\documentclass[12pt]{article} 

\usepackage{enumerate}

\usepackage{algorithm}

\usepackage{algorithmic}

\usepackage{amsmath}
\usepackage{amsfonts}
\usepackage{amssymb}
\usepackage{url}

\DeclareMathAlphabet{\eufrak}{U}{}{}{} 
\SetMathAlphabet\eufrak{normal}{U}{euf}{m}{n}
\SetMathAlphabet\eufrak{bold}{U}{euf}{b}{n}
\numberwithin{equation}{section}

\newenvironment{Proof}{\removelastskip\par\medskip
\noindent{\em Proof.} \rm}{\penalty-20\null\hfill$\square$\par\medbreak}

\usepackage{color}
\definecolor{lightblue}{rgb}{0,0.2,0.5}

\usepackage[colorlinks=true, urlcolor=blue,linkcolor=blue, citecolor=lightblue]{hyperref}

\allowdisplaybreaks

 \def\ind{{\bf 1}}
 
 \def\inte{{\mathord{\mathbb N}}}
 
 \def\qu{{\mathord{\mathbb Z}}}

 \def\inte{{\mathord{{\rm I\kern-3pt N}}}}
 \def\sZZ{{\rm Z\kern-.45em{}Z}}
 \def\sQQ{{\kern 0.27em \vrule height1.45ex width0.03em depth0em
           \kern-0.30em \rm Q}}
 \def\qu{{\mathchoice
         {\sQQ}
         {\sQQ}
   {\kern 0.225em \vrule height1.05ex width0.025em depth0em \kern-0.25em \rm Q}
   {\kern 0.180em \vrule height0.78ex width0.020em depth0em \kern-0.20em \rm Q}
         }}
 \def\sGG{{\kern 0.27em \vrule height1.45ex width0.03em depth0em
           \kern-0.30em \rm G}}
 \def\gg{{\mathchoice
         {\sGG}
         {\sGG}
   {\kern 0.225em \vrule height1.05ex width0.025em depth0em \kern-0.25em \rm G}
   {\kern 0.180em \vrule height0.78ex width0.020em depth0em \kern-0.20em \rm G}
         }}

 \newtheorem{prop}{Proposition}[section]
 \newtheorem{lemma}[prop]{Lemma}
 
 \newtheorem{corollary}[prop]{Corollary}

\def\E{\mathop{\hbox{\rm I\kern-0.20em E}}\nolimits}

\def\inte{{\mathord{\mathbb N}}}

\def\qu{{\mathord{\mathbb Q}}}

\def\PP{{\mathord{{\rm I\kern-2.8pt P}}}}

\def\E{\mathop{\hbox{\rm I\kern-0.20em E}}\nolimits}

\def\Det{\mathop{\hbox{\rm Det}}\nolimits}

\def\inte{{\mathord{{\rm I\kern-2.8pt N}}}}
\def\PP{{\mathord{{\rm I\kern-2.8pt P}}}}

\def\real{{\mathord{\mathbb R}}}

\def\inte{{\mathord{\mathbb N}}}

\def\qu{{\mathord{\mathbb Z}}}

 \newcounter{hyp}
\usepackage{fixltx2e,amsmath}
\MakeRobust{\eqref}

\title{\huge
 Second-order multi-object filtering with target interaction using determinantal point processes 
}

\author{
  \large
  Nicolas Privault\thanks{nprivault@ntu.edu.sg}
  ~and Timothy Teoh\thanks{teoh0094@e.ntu.edu.sg} 
\\ 
\normalsize 
Division of Mathematical Sciences 
\\ 
\normalsize 
School of Physical and Mathematical Sciences 
\\ 
\normalsize 
Nanyang Technological University 
\\ 
\normalsize 
21 Nanyang Link 
\\ 
\normalsize 
Singapore 637371
}

\usepackage{xcolor,graphicx}

\usepackage{mathrsfs}

\usepackage{amssymb,amsmath}
\usepackage{comment}
\usepackage{enumerate}

\usepackage{multicol}

\usepackage{subcaption}

\usepackage[round]{natbib} 

\bibpunct{\textcolor{lightblue}{(}}{\textcolor{lightblue}{)}}{,}{a}{}{;}

\let\oldcitet=\citet
\let\oldcitep=\citep 
\renewcommand{\cite}[1]{\textcolor[rgb]{0,0,1}{\oldcitet{#1}}}
\renewcommand{\citet}[1]{\textcolor[rgb]{0,0,1}{\oldcitet{#1}}}
\renewcommand{\citep}[1]{\textcolor[rgb]{0,0,1}{\oldcitep{#1}}}

\newcommand{\underlying}{{\Lambda}} 

\newcommand{\underlyingsigma}{\mathcal{B}(\underlying )}

\def\espaceconfig{\mathbf N_\sigma}

\def\topoconfig{{\mathscr N_\sigma}}

\def\topoconfig{{\mathscr N_\sigma}}

\def\espaceconfiglambda{\mathbf N_\sigma(\Lambda)}

\newcommand{\pp}{{\mathbb P}}

\renewcommand{\d}{d} 

\newcommand{\X}{{\espaceconfig}}

\newcommand{\radon}{{\nu}}

\usepackage{dsfont}

\usepackage{float}

\def\espaceconfig{\mathbf N_\sigma}

\newcommand{\Tr}{\mathrm{Tr}}

\begin{document}
 
\maketitle

\baselineskip0.6cm
 
\begin{abstract}
 The Probability Hypothesis Density (PHD) filter, 
 which is used for multi-target tracking based on sensor measurements, 
 relies on the propagation of the first-order moment,
 or intensity function, of a point process. 
 This algorithm assumes that
 targets behave independently, 
 an hypothesis which may not hold in practice
 due to potential target interactions. 
 In this paper, 
 we construct a second-order PHD filter based on Determinantal Point Processes (DPPs)
 which are able to model repulsion between targets. 
 Such processes are characterized by their first and second-order moments,
 which allows the algorithm 
 to propagate variance and covariance information
 in addition to first-order target count estimates. 
 Our approach relies on posterior moment formulas for the estimation of a general hidden point process after a thinning operation and a superposition with a Poisson Point Process (PPP), and on suitable approximation formulas in the determinantal point process setting.
 The repulsive properties of determinantal point processes apply 
 to the modeling of negative correlation between distinct
 measurement domains. Monte Carlo simulations with correlation estimates are provided. 
\end{abstract}

\noindent {\bf Keywords}: 
Probability hypothesis density (PHD) filter; higher-order statistics; correlation; second-order moment; determinantal point processes; multi-object filtering; multi-target tracking.

\noindent
{\em Mathematics Subject Classification (2010):} 60G35; 60G55; 62M30; 62L12. 
\baselineskip0.7cm

\section{Introduction} 
Probability Hypothesis Density (PHD) filters
have been introduced in \cite{mahler2003}
for multi-target tracking in cluttered environments.
The construction of the prediction point process $\Phi$
therein uses multiplicative point processes, see e.g. \cite{moyal}, \cite{moyal2},
by thinning and shifting a prior point process $\Psi$,
and superposition with a birth point process.
The posterior point process $\Phi | \Xi$ is obtained by conditioning
$\Phi$ given a measurement point process of targets $\Xi$, also constructed by
thinning, shifting and superposition. 
This step relies on Bayesian estimation with a Poisson point process prior,
see e.g. \cite{vanlieshout2}, \cite{mori}, \cite{portenko}. 
PHD filters have low complexity, and they allow for explicit update formulas
see e.g. \cite{cdh} for a review. 

\bigskip

While the PHD filter of \cite{mahler2003} is based on Poisson point processes,
several extensions of the PHD filter to non-Poisson prior distributions
have been proposed.
Cardinalized Probability Hypothesis Density (CPHD) filters
have been introduced in \cite{mahler2007} as
a generalization in which the target count is allowed to have
an arbitrary distribution.
In \cite{maskell}, discretized Gamma distributions
are used to design an efficient approximation
of the CPHD filter cardinality distribution. 
Other generalizations include the 
Gauss-Poisson point processes that generalize the Poisson point process
by allowing for two-point clusters,
and have been used in \cite{singh}.
The PHD filter has been implemented using the Sequential Monte Carlo (SMC) method
in \cite{vo-singh}, and using Gaussian mixtures in \cite{vo-ma}.

\bigskip

PHD filters approximate the distribution of the number of targets
by a Poisson distribution estimated
by a single mean (or variance) parameter, which can result
into high variance estimates when the estimated mean is high.
Second-order PHD filters that can propagate
distinct information on mean and variance parameters have been
recently proposed in \cite{sdhc}, based on the Panjer point process
defined therein, where the Panjer cardinality distribution encompasses the
binomial, Poisson and negative binomial distributions. 
Other multi-target filters propagating second-order moment information have
also been recently proposed, see e.g. 
\cite{demelo} for a filter that propagates second-order point process 
factorial cumulants. 

\bigskip

A common feature of cardinalized filters is to assume that
target locations are distributed
as $n$ independent random samples according
to a reference intensity measure,
given that the observation window contains $n$ points. 
While this hypothesis is natural and facilitates an explicit
derivation of prediction formulas, it does not
reflect potential interaction between targets.
In addition, as observed in the simulations of Section~\ref{s7},
the presence of repulsion between targets
can degrade the performance of the Poisson PHD filter. 

\bigskip

As a response, 
we propose to construct a PHD filter based on determinantal point
processes introduced in \cite{macchi},
which are able to model repulsion among configuration points 
on a target domain $\Lambda \subset \real^d$,
see also \cite{soshnikov} and \cite{shirai}.
Taking into account correlation via 
more general point process-based PHD filters poses
several challenges linked to the derivation of closed form
filtering formulas. 
In addition, the distribution of general point processes
relies on Janossy densities which may not be characterized by
the knowledge of a finite number of moments. 
Determinantal point processes,
on the other hand, are characterized by their kernel functions
$(K(x,y))_{(x,y) \in \Lambda^2}$, and their Janossy densities can be recovered from 
their first and second-order moments.
In this setting, the knowledge of first and second-order moments
can be used to update the Janossy densities that
characterize the underlying determinantal point process.

\bigskip

Discrete determinantal point processes
have also been recently used for the pruning of Gaussian components
in the Gaussian Mixture (GM) PHD filter in \cite{jorquera2}, 
see also \cite{jorquera} for other applications to multi-target tracking.
Permanental processes have been used in \cite{mahlerpermanental}
to propagate a joint Poisson distribution, however this approach
is distinct from the determinantal setting. 
See \cite{koch} for the use the exclusion principle in multi-target 
 tracking by a fermionic filtering update 
 of anti-symmetric components in the joint probability density functions
 of states. 
 Note also that determinantal point processes have been originally
 introduced in \cite{macchi}
 to represent configurations of fermions.

\bigskip

After recalling general facts and notation on point processes in
Section~\ref{s2}, we derive general formulas for the
distribution, and for the first and second-order 
moments, of a posterior point process $\Phi | \Xi$
in Section~\ref{s2.1}, see also \cite{lund}.
In Section~\ref{s2.2} we review the construction of determinantal point processes,
and 
in Section~\ref{s5} we present a second-order PHD filtering algorithm based
on determinantal point processes, with the computation 
of the prediction kernel $K_\Phi (x,y )$ and of the updated kernel
$K_{\Phi \mid \Xi} (x,y )$. 
An implementation of the Poisson PHD filter that allows for
performance evaluation using measurement-estimate associations
is presented in Section~\ref{s7} with numerical illustrations. 
This simulation is based on
the sequential Monte Carlo (or particle filtering) method 
with a nearly constant turn-rate motion dynamics, see 
\cite{cantoni}, \cite{litiancheng}, to which we add a repulsion term.
In Section~\ref{s7.1} 
we implement the determinantal PHD filter using the sequential Monte Carlo method. 
The implementation of the algorithm relies on closed-form filter update expressions
obtained from approximation formulas for corrector terms and Janossy densities
presented in appendix. 

\section{Preliminaries on point processes}
\label{s2}
In this section we review the properties of
point processes; see, e.g. \cite{daley}, 
\cite{dfpt2}, and references therein. 
For any subset
$A \subseteq \real^d$, let $| A|$ denote the cardinality of $A$,
setting $| A|=\infty$ if $A$ is not finite, and
let
\begin{equation*}
\espaceconfig:=\{\xi\subseteq \real^d \ :\ |\xi \cap
A |<\infty\quad\text{ for all relatively compact sets }A \subset\real^d \}
\end{equation*}
 denote the set of locally finite point configurations
 on $\real^d$, which is identified with the set of all 
 nonnegative integer-valued Radon measures $\xi$ on $\real^d$
 such that $\xi (\{x\}) \in\{0,1\}$
 for all $x \in \real^d$.
 We denote by $\topoconfig$ the Borel $\sigma$-field
 generated by the weakest topology that makes the mappings 
\begin{equation*}
	\xi \mapsto \langle f, \xi \rangle := \sum_{y \in \xi} f(y)
\end{equation*}
continuous for all continuous and compactly supported functions $f$ on $\real^d$.
Given $\Lambda$ a relatively compact subset of $\real^d$, we 
let $\espaceconfiglambda$ be the space of finite configurations on $\Lambda$. 

\bigskip 

We consider a simple and locally finite point process $\Phi$ on $\Lambda$, 
defined as a random element on a probability space $(\Omega,\topoconfig )$
with values in $\espaceconfiglambda$, and denote its distribution by
$\pp$.
The point process $\Phi$ 
is characterized by its Laplace transform $\mathcal{L}_\Phi$ which is defined,
for any measurable nonnegative function $f$ on $\underlying$, by
\begin{equation}
  \label{lapl}
  \mathcal{L}_\Phi (f) = \int_\X e^{- \langle f, \xi \rangle} \,\pp(\d \xi)	.
\end{equation}
We denote the expectation of an integrable random variable
$F$ defined on $(\espaceconfig,\topoconfig,\pp)$ by
\begin{equation*}
\E[F(\Phi) ] := \int_{\X} F(\xi) \,\pp(\d \xi).
\end{equation*}

\subsubsection*{Janossy densities}

For any relatively compact subset $A \subseteq \underlying$, the Janossy densities of $\Phi$ {\em w.r.t.} a reference Radon measure $\radon$ on $\underlying$ 
are symmetric measurable functions $j^{(n)}_{\Lambda}:\Lambda^n\rightarrow [0,\infty)$
  satisfying 
\begin{equation}
\nonumber 
\E\left[ F(\bold \Phi ) \right] =
 F(\emptyset )\,j_{\Lambda}^{(0)}
 + \sum_{n \ge 1} \frac{1}{n!} \int_{\Lambda^n} F(\{x_1,\ldots,x_n\})\,j_{\Lambda}^n\left( x_{1}, \ldots, x_{n}\right)\,\radon(\d x_1)\cdots\radon(\d x_n), 
\end{equation}
for all measurable functions $F:\espaceconfiglambda \rightarrow [0,\infty)$; 
 see, e.g., \cite{georgiiyoo}.

\bigskip 

For $n \ge 1$, the Janossy density $j^{(n)}_\Lambda (x_1,\ldots,x_n)$
is proportional, up to a multiplicative constant, to 
the joint density of the $n$ points
of the point process, given that it has exactly $n$ points.
 For $n=0$, 
 $j^{(0)}_\Lambda (\emptyset)$ is the probability that there are no points in $\Lambda$.
\subsubsection*{Correlation functions}
 The correlation functions of $\Phi$ {\em w.r.t.} the 
 reference measure $\radon$ on $\underlying$ are measurable symmetric functions
 $\rho^{(k)}_\Phi :\underlying^k \longrightarrow [0,\infty)$ such that
\begin{equation}
\label{coor-in}
\E\left[\prod_{i=1}^{k}\bold \Phi (B_i)\right]=\int_{B_1\times\cdots\times
B_k}\rho^{(k)}_\Phi (x_1,\ldots,x_k)\,\radon(\d x_1)\cdots\radon(\d x_k),
\end{equation}
 for any
 family of mutually disjoint bounded subsets $B_1,\ldots,B_k$ of $\underlying$,
 $k\geq 1$. More generally,
 if $B_1,\hdots,B_n$ are disjoint bounded Borel subsets of $\underlying$
 and $k_1,\hdots,k_n$ are integers such that $\sum_{i=1}^n k_i = N$, 
 we have 
	\begin{equation*}
		\E\left[\prod_{i=1}^n \frac{\Phi (B_i)!}{(\Phi (B_i)-k_i)!}\right] = \int_{B_1^{k_1}\times \cdots \times B_n^{k_n}} \rho^{(N)}_\Phi (x_1,\hdots,x_N ) \,\radon(\d x_1)\cdots \radon(\d x_N)	.
	\end{equation*}
 In addition, we let 
 $\rho^{(n)}_\Phi (x_1,\ldots,x_n) =0$
 whenever $x_i=x_j$ for some $1 \leq i\neq j \leq n$.
In other words, the factorial moment density
$\rho^{(n)}_\Phi (x_1,\ldots , x_n)$ of $\Phi$, 
$x_1,\ldots , x_n \in \Lambda$, $x_i \not= x_j$, $1\leq i<j \leq n$,
is defined from the relation
$$   \int_{B_1\times \cdots \times B_n} \rho^{(n)}_\Phi (x_1,\ldots , x_n)
  \nu (dx_1) \cdots \nu (dx_n )
   = 
\E\left[
  \sum_{x_1,\ldots , x_n\in \Phi}
  \ind_{B_1}(x_1)\cdots 
  \ind_{B_n}(x_n)
  \right],
$$ 
for mutually disjoint measurable subsets
$B_1,\ldots, B_n \subset \Lambda$,
where $\ind_{B_i}$ denotes the indicator function over $B_i$,
$i=1,\ldots , n$. 
Heuristically, 
$	\rho^{(n)}_\Phi (x_1,\ldots,x_n) \, \radon(\d x_1)\cdots \radon(\d x_n)
$ 
 represents
 the probability of finding a particle in the vicinity of each $x_i$, $i=1,\ldots,n$.

\bigskip 

We also recall that the Janossy densities $j^{(n)}_\Lambda$
can be recovered from the correlation functions
$\rho^{(m)}_\Phi$ via the relation
\begin{equation*}
j^{(n)}_\Lambda (x_1,\ldots,x_n) = \sum_{m \ge 0} \frac{(-1)^m}{m!} \int_{\Lambda^m} \rho^{(n+m)}_\Phi (x_1,\ldots,x_n,y_1,\ldots,y_m)\,\radon(\d y_1)\cdots \radon(\d y_m),
\end{equation*}
and vice versa using the equality
\begin{equation*}
\rho^{(n)}_\Phi (x_1,\ldots,x_n) = \sum_{m \ge 0} \frac{1}{m!} \int_{\Lambda^m} j^{(m+n)}_\Lambda (x_1,\ldots,x_n,y_1,\ldots,y_m)\,\radon(\d y_1)\cdots \radon(\d y_m),
\end{equation*}
see Theorem~5.4.II of \cite{daley}.

\subsubsection*{Probability generating functionals} 
The Probability Generating Functional (PGFl)
of the point process $\Phi$, see \cite{moyal}, is defined by 
\begin{eqnarray*} 
  {
    h \mapsto {\cal G}_\Phi (h)
  }
  & := & 
  {
    \E\left[
  \prod_{i=1}^{\Phi ( \Lambda )}
  h(X_i)
  \right]
  }
  \\
& = &
  {
    j^{(0)}_\Phi 
+
\sum_{n\geq 1}
\frac{1}{n!}
\int_{\Lambda^n} j^{(n)}_\Phi (x_1, \ldots , x_n ) \prod_{i=1}^n h(x_i)
\ \!
\nu (dx_1 )
\cdots \nu (dx_n), 
  }
\end{eqnarray*} 
 for $h \in L^\infty (\Lambda )$ a bounded measurable function on $\Lambda$.
 Given ${\cal F}$ a functional on $L^\infty ( \Lambda )$,
 we will use the functional derivative
 $\partial_g / \partial h$
 of ${\cal F}(h)$ in the direction of
 $g \in L^\infty ( \Lambda )$, defined as 
 $$
 \frac{\partial_g }{\partial h} {\cal F}(h) 
 := \lim_{\varepsilon \to 0} \frac{{\cal F}(h+ \varepsilon g) - {\cal F}(h)}{\varepsilon}.
 $$  
  Given $x\in \Lambda$, we also let 
  \begin{equation}
    \label{dirac} 
 \frac{\partial_{\delta_x} }{\partial h} {\cal F}(h) 
 := \lim_{n \to \infty}
 \frac{\partial_{g_n} }{\partial h} {\cal F}(h) 
 ,
\end{equation} 
 where $(g_n)_{n \geq 1}$ is a sequence of bounded functions
 converging weakly to the Dirac distribution $\delta_x$ at $x\in \Lambda$.

\bigskip 

 This construction allows one to recover
 the Janossy densities $j^{(n)}_\Phi (x_1,\ldots , x_n)$
and factorial moment densities
$\rho^{(n)}_\Phi (x_1,\ldots , x_n)$ 
of $\Phi$ from the PGFl ${\cal G}_\Phi (h)$ as
\begin{equation}
  \label{jnsy} 
{
  j^{(n)}_\Phi (x_1,\ldots , x_n)
=
\frac{\partial_{\delta_{x_1}} }{\partial h}
\cdots
\frac{\partial_{\delta_{x_n}} }{\partial h}
     {\cal G}_{\Phi } (h)_{\mid h=0},
     \quad
     x_1,\ldots , x_n \in \Lambda, 
}
\end{equation} 
see e.g. \S~2.4 of \cite{cdh},
and as 
\begin{equation}
  \label{**} 
     {
       \rho^{(n)}_\Phi (x_1,\ldots , x_n)
=
\frac{\partial_{\delta_{x_1}} }{\partial h}
\cdots
\frac{\partial_{\delta_{x_n}} }{\partial h}
     {\cal G}_{\Phi } (h)_{\mid h=1},
     \quad
x_1,\ldots , x_n \in \Lambda, 
     }
\end{equation} 
 with $x_i \not= x_j$, $1\leq i<j \leq n$; see, e.g., \S~3.4 of \cite{cdh}. 
\section{Posterior point process distribution}
\label{s2.1}
In this section we compute the Janossy densities, and 
the first and second-order moments, of 
a posterior point process of targets $\Phi$ given the
point process $\Xi$ of sensor measurements.
The case of a Poisson prior was treated in Chapter~5 of \cite{vanlieshout2},
see Theorem~29 therein and also \cite{mori} for early
sensor fusion applications, or Theorems~6.1-6.2 \cite{portenko} for related derivations
based on Laplace transforms. 
\\

In Propositions~\ref{p1} and \ref{p2} below 
 we will use extensions of the corrector terms $l^{(1)}_{z_{1:m}}$, $l^{(2)}_{z_{1:m}}$ 
introduced in \cite{duhc} for the cardinalized PHD filter,
see \cite{batuongvo}.
We start with a review of the thinning and shifting of point processes; 
 see, e.g., \cite{cdh} and references therein for details.
\subsubsection*{Thinning and shifting of point processes} 
 The point process $\Xi$ of
 sensor measurements is constructed via the following steps. 
\begin{enumerate}[(i)]
\item
 Thinning and shifting. 
 Every target point $x \in \Phi$ is kept with probability 
 $p_d(x)\in (0,1]$ 
 and shifted according to the probability density function
 $l_d (\cdot | x)$, by branching the hidden point process $\Phi$ 
 with a Bernoulli point process $\Xi_s$ with PGFl 
 \begin{eqnarray}
   \nonumber 
  g \mapsto {\cal G}_{\Xi_s} ( g \mid x ) & : = &
  q_d(x) + p_d(x) \int_\Lambda g(z) l_d (z|x) \nu (dz)
  \\
  \label{abdfasfd}
  & = &
  q_d(x) + \int_\Lambda g(z) \tilde{l}_d (z|x) \nu (dz),
\end{eqnarray} 
 where for compactness of notation we take
\begin{equation}
\label{qd}
 q_d (x):= 1 - p_d (x) \quad
\mbox{and}
\quad 
\tilde{l}_d (z|x) : = p_d(x) l_d (z|x), \quad x\in \Lambda.
\end{equation} 
\item The point process $\Xi$ is obtained by
  superposing a Poisson point process $\Xi_c$ with intensity
  function $l_c (\cdot )$, representing clutter,
  to the above thinning and shifting of $\Phi$. 
\end{enumerate}
 In the sequel we use the shorthand notation 
 $$
 x_{1:n} = (x_1,\ldots ,x_n)\in \Lambda^n, 
 \quad
 \mbox{and}
 \quad 
 \nu (dx_{1:n} ): = \nu (dx_1 )\cdots \nu (dx_n ),
 $$
 with $x_{1:0}=\emptyset$; see \cite{duhc}, \cite{sdhc}. 
 The joint PGFl of the point process $( \Phi , \Xi )$ is given by 
\begin{equation}
  \label{mspt} 
 (h,g) \mapsto
{\cal G}_{ \Phi , \Xi } (h,g) 
 := {\cal G}_{\Xi_c} (g)
 {\cal G}_\Phi \big( h (\cdot ) {\cal G}_{\Xi_s} (g \mid \cdot \ \! ) \big),
\end{equation} 
see, e.g., Theorem~1.1 of \cite{moyal2}
where, taking $p_d(x):=p_d$ for simplicity, $x\in \Lambda$, we have 
   \begin{align*} 
       {\cal G}_\Phi \big( h (\cdot ) {\cal G}_{\Xi_s} (g \mid \cdot \ \! ) \big)
& = 
       \sum_{n\geq 0}
\frac{1}{n!}
\int_{\Lambda^n}
j^{(n)}_\Phi (x_{1:n})
 {
 \prod_{i=1}^n
\left( h(x_i) \left(
q_d + \int_\Lambda g(z) \tilde{l}_d (z|x_i ) \nu (dz) \right)
\right) 
\nu (dx_{1:n})
}
    \\
 & =  
    {
      \sum_{n\geq 0}
\frac{1}{n!}
\int_{\Lambda^n}
j^{(n)}_\Phi (x_{1:n})
\prod_{j=1}^n
h(x_j)
    }
    {
 \sum_{k=0}^n
    {n \choose k}
    q_d^{n-k}
    \prod_{i=1}^k 
    \int_\Lambda g(z) \tilde{l}_d (z|x_i ) \nu (dz) 
\nu (dx_{1:n})
    }
    \\
 & =   
    {
      \sum_{k\geq 0}
    \int_{\Lambda^k}
    \prod_{i=1}^k 
    \int_\Lambda g(z) \tilde{l}_d (z|x_i ) \nu (dz) 
    }
     {
 \sum_{n\geq 0}
    \frac{    q_d^n}{n!k!}
    \int_{\Lambda^n}
    j^{(k+n)}_\Phi (x_{1:k+n})
\prod_{j=1}^{k+n}  h(x_j) 
    \nu (dx_{1:k+n}). 
    }
   \end{align*} 
 The marginal PGFl of the point process $\Xi$ is given by 
 \begin{eqnarray} 
\nonumber
   g \mapsto
{\cal G}_\Xi (g) 
& = & 
{\cal G}_{  \Phi, \Xi } (\ind,g) 
\\
\nonumber
& = & {\cal G}_{\Xi_c} (g)
 {\cal G}_\Phi \big( {\cal G}_{\Xi_s} (g \mid \cdot \ \! ) \big)
 \\
   \label{mgnal}
  & = & 
  \mathcal{G}_{\Xi_c}(g)
  \mathcal{G}_\Phi\left(
  q_d(\cdot)+\int_\Lambda g (y) \tilde{l}_d (y\ \! |\ \! \cdot)\nu (dy) \right). 
\end{eqnarray} 
\subsubsection*{Marginal moments of $\Xi$} 
The first derivative of $  \mathcal{G}_\Xi (g)$
 is given by
\begin{eqnarray*}
    \frac{\partial_{\delta_x}}{\partial g}
  \mathcal{G}_\Xi (g)
 & = & 
  \frac{\partial_{\delta_x}}{\partial g}\big(\mathcal{G}_\Phi(\mathcal{G}_{\Xi_s}(g \mid \cdot))\mathcal{G}_{\Xi_c}(g)\big)
    \\
  &= & 
  \mathcal{G}_{\Xi_c}(g)
  \frac{\partial_{\delta_x}}{\partial g}\mathcal{G}_\Phi(\mathcal{G}_{\Xi_s}(g \mid \cdot)
  +\mathcal{G}_\Phi(\mathcal{G}_{\Xi_s}(g \mid \cdot))
  \frac{\partial_{\delta_x}}{\partial g}\mathcal{G}_{\Xi_c}(g)
  \\
  &= & 
  \mathcal{G}_{\Xi_c}(g)
  \frac{\partial_{\delta_x}}{\partial g}\mathcal{G}_\Phi
  \left(
  q_d(\cdot)+\int_\Lambda g(y)\tilde{l}_d (y\ \! |\ \! \cdot)\nu(dy))
  \right)
  \\
  &  & +\mathcal{G}_\Phi \left(
  q_d(\cdot)+\int_\Lambda g(y)\tilde{l}_d (y\ \! |\ \! \cdot)\nu(dy)
  \right)
  \frac{\partial_{\delta_x}}{\partial g}\mathcal{G}_{\Xi_c}(g), 
\end{eqnarray*}
from which the first-order moment density of $\Xi$ can be computed
after setting $g=1$ as 
\begin{equation} 
 \label{1stfm}
 {
   \mu^{(1)}_\Xi(x)=
  \frac{\partial_{\delta_{x}} }{\partial g}
  \mathcal{G}_\Xi(g)_{\mid g=1}
= \mu^{(1)}_{\Xi_c}(x) + \int_\Lambda \tilde{l}_d (x|y)\mu^{(1)}_\Phi(y)\nu (dy).
 }
\end{equation} 
Similarly, the second derivative of $\mathcal{G}_\Xi (g)$ is given by 
\begin{align*}
    {
      \frac{\partial_{\delta_x}}{\partial g}\frac{\partial_{\delta_y}}{\partial g}\bigg(\mathcal{G}_\Phi(\mathcal{G}_{\Xi_s}(g \mid \cdot))\mathcal{G}_{\Xi_c}(g)\bigg)
  }
  & = 
  {
    \mathcal{G}_{\Xi_c}(g)  \frac{\partial_{\delta_x}}{\partial g}\frac{\partial_{\delta_y}}{\partial g}\mathcal{G}_\Phi(\mathcal{G}_{\Xi_s}(g \mid \cdot))
+
  \frac{\partial_{\delta_x}}{\partial g}\mathcal{G}_\Phi(\mathcal{G}_{\Xi_s}(g \mid \cdot))
  \frac{\partial_{\delta_y}}{\partial g}\mathcal{G}_{\Xi_c}(g)
  }
  \\
  & 
  {
    \quad
    + \frac{\partial_{\delta_y}}{\partial g}\mathcal{G}_\Phi(\mathcal{G}_{\Xi_s}(g \mid \cdot))
  \frac{\partial_{\delta_x}}{\partial g}\mathcal{G}_{\Xi_c}(g)
  +
  \mathcal{G}_\Phi(\mathcal{G}_{\Xi_s}(g \mid \cdot))
  \frac{\partial_{\delta_x}}{\partial g}\frac{\partial_{\delta_y}}{\partial g}(\mathcal{G}_{\Xi_c}(g)),
  }
\end{align*}
 from which the second-order factorial moment density of $\Xi$ can
 be computed after setting $g=1$ as 
 \begin{eqnarray}
   \label{2ndfm} 
   \rho^{(2)}_{\Xi}(x,y) & = & \int_{\Lambda^2}\tilde{l}_d (x|u)
   \tilde{l}_d (y|v)\rho^{(2)}_{\Phi}(u,v)\nu(du)\nu(dv)
       +\rho^{(2)}_{\Xi_c}(x,y)
   \\
     \nonumber
   & &
   {
     +\mu^{(1)}_{\Xi_c}(y)\int_{\Lambda}
     \tilde{l}_d (x|u)\mu^{(1)}_{\Phi}(u)\nu(du)
}
     {
       +\mu^{(1)}_{\Xi_c}(x)\int_{\Lambda}
       \tilde{l}_d (y|v)\mu^{(1)}_{\Phi}(v)\nu (dv)
 , 
   }
 \end{eqnarray} 
$x,y\in \Lambda$, $x\not= y$. 
\subsubsection*{Posterior distribution} 
In Lemma~\ref{ldjkldsf} below, we derive the general expression
of the Janossy densities of the posterior point process $\Phi | \Xi$
given the sensor measurements $\Xi$.
In the sequel, we let $|S|$ denote the cardinality of subsets
$S\subset \{1,\ldots , m\}$, and we use the notation
$z_{1:m} = (z_1, \ldots , z_n)$,
while 
$z_{1:m} \setminus z$ denotes the sequence
$(z_1, \ldots , z_n)$ with the omission of $z$
if $z\in \{z_1,\ldots , z_n\}$. 
\begin{lemma}
  \label{ldjkldsf} 
The $n$-th conditional Janossy density
of $\Phi$ given that $\Xi = z_{1:m} = (z_1,\ldots , z_m)$ 
satisfies 
\begin{equation} 
\label{x012} 
  j^{(n)}_{\Phi \mid \Xi = z_{1:m}} (x_1, \ldots , x_n)
   = 
   \frac{  j^{(n,m)}_{\Phi , \Xi = z_{1:m}} (x_1, \ldots , x_n)}{j^{(m)}_\Xi (z_{1:m})},
   \qquad
   x_1, \ldots , x_n \in \Lambda, 
\end{equation} 
$m,n\geq 0$, where
\\
$(i)$ 
  the $(n,m)$-th joint Janossy density of $(\Phi , \Xi)$ is given by 
\begin{equation} 
  \label{x1} 
      j^{(n,m)}_{\Phi , \Xi = z_{1:m}} (x_1, \ldots , x_n)
      =
      {
    j^{(n)}_\Phi (x_1,\ldots , x_n)
  }
  {
 \sum_{S\subset \{1,\ldots , m\} \atop |S| \leq n } 
  \frac{n!   q_d^{n-|S|} }{(n-|S|)!}
 \prod_{j \notin S} l_c(z_j) 
 \sum_{\pi : S \rightarrow \{1,\ldots , n\} } 
    \prod_{i \in S} \tilde{l}_d (z_i | x_{\pi (i)}) 
, }
\end{equation} 
the above sum being over injective
mappings $\pi : S \rightarrow \{1,\ldots , n\}$, and 
\\
$(ii)$ the Janossy densities of the measurement point process $\Xi$ are given by 
\begin{equation} 
  \label{jd}
    {
      j^{(m)}_\Xi (z_1,\ldots , z_m)
      =
      \sum_{n\geq 0}
      \sum_{S\subset \{1,\ldots , m\} \atop
|S|\leq n } 
      }
  {
     \prod_{j \notin S} l_c(z_j) 
     \frac{    q_d^{n-|S|}}{(n-|S|)!} 
    \int_{\Lambda^n} 
    j^{(n)}_\Phi (x_{1:n})
    \prod_{i \in S} \tilde{l}_d (z_i | x_i) 
    \nu (dx_{1:n}), 
  }
\end{equation} 
$m\geq 0$. 
\end{lemma}
\begin{Proof}
  In order to derive 
  the $(n,m)$-th joint Janossy density of $(\Phi,\Xi)$
  as in \eqref{jnsy}, we need to compute   
  $$
  \frac{\partial_{\eta_1}}{\partial h}\cdots\frac{\partial_{\eta_n}}{\partial h}\frac{\partial_{f_1}}{\partial g}\cdots\frac{\partial_{f_m}}{\partial g}\mathcal{G}_{\Phi,\Xi}(h,g)
    $$
  in the directions of the functions $\eta_1,\ldots ,\eta_n$, $f_1,\ldots ,f_m
  \in L^\infty ( \Lambda )$. 
  For a given set $S\subseteq\{1,\ldots ,m\}$
  we let $\pi:=\left|S\right|$ and denote the elements of $S$
  as $S(1),\ldots ,S(\pi)$ in increasing order,
  where $S$ is identified to the
  mapping $S:\{1,\ldots ,\pi\}\rightarrow\{1,\ldots ,m\}$. 
  By the Fa\`a di Bruno's formula, see, e.g., \cite{Faa},
  and \eqref{mspt}, we have 
    \begin{align}
&\frac{\partial_{f_1}}{\partial g}\cdots\frac{\partial_{f_m}}{\partial g}\mathcal{G}_{\Phi,\Xi}(h,g)=\frac{\partial_{f_1}}{\partial g}\cdots\frac{\partial_{f_m}}{\partial g}\big(\mathcal{G}_{\Xi_c}(g)\mathcal{G}_\Phi(h(\cdot)\mathcal{G}_{\Xi_s}(g\ \! |\ \! \cdot))\big)\nonumber\\
      &=\sum_{\substack{S\subseteq\{1,\ldots ,m\}\\ \pi=\left|S\right|\\ Q \subset \{1,\ldots ,m\}\backslash S}}
      \left(
      \frac{\partial_{f_{S(1)}}}{\partial g}\cdots\frac{\partial_{f_{S(\pi)}}}{\partial g}\mathcal{G}_\Phi(h(\cdot)\mathcal{G}_{\Xi_s}(g\ \! |\ \! \cdot))\right)
      \left( \frac{\partial_{f_{Q(1)}}}{\partial g}\cdots\frac{\partial_{f_{Q(m-\pi)}}}{\partial g}\mathcal{G}_{\Xi_c}(g)
      \right),\label{bruno1}
\end{align}\noindent
    where the index set $S\subseteq\{1,\ldots ,m\}$ runs through
    the collection of $2^m$ subsets of $\{1,\ldots ,m\}$.
 \noindent
 Next,
 the $m^{\text{th}}$ derivative of $\mathcal{G}_\Phi(h(\cdot)\mathcal{G}_{\Xi_s}(g\ \! |\ \! \cdot))$
 can be computed by a standard induction argument as 
 \begin{align}
   \label{8.33} 
  &\frac{\partial_{f_1}}{\partial g}\cdots\frac{\partial_{f_m}}{\partial g}\mathcal{G}_\Phi(h(\cdot)\mathcal{G}_{\Xi_s}(g\ \! |\ \! \cdot))
  \\
  \nonumber
  &=\sum_{a=m}^\infty\frac{1}{a!}\int_{\Lambda^a}j_\Phi^{(a)}(u_{1:a})
  \hskip-0.2cm
  \sum_{\substack{s_1,\ldots ,s_m=1\\s_1\neq\cdots\neq s_m}}^a \prod_{l=1}^m\Big(h(u_{s_l})\int_\Lambda f_l(z)\tilde{l}_d(z|u_{s_l})\nu(dz)\Big)
  \hskip-0.5cm
  \prod_{s\in\{1,\ldots ,a\}\backslash\{s_1,\ldots ,s_m\}}^a
  \hskip-1cm
  \big(h(u_s)\mathcal{G}_{\Xi_s}(g|u_s)\big)\nu(du_{1:a}).
 \end{align}
 Substituting $f_1,\ldots ,f_m$ with Dirac delta functions 
$\delta_{z_1},\ldots ,\delta_{z_m}$ at the distinct configuration
 points $z_{1:m}\in\Lambda$ as in \eqref{dirac}
 and setting $g=0$ 
in \eqref{8.33}, we find 
\begin{align}
  \nonumber
  \frac{\partial_{\delta_{z_1}}}{\partial g}\cdots\frac{\partial_{\delta_{z_m}}}{\partial g}\mathcal{G}_\Phi(h(\cdot)\mathcal{G}_{\Xi_s}(g\ \! |\ \! \cdot))_{|g=0}
  &=\sum_{a=m}^\infty\frac{q_d^{a-m}}{a!}\int_{\Lambda^a}j_\Phi^{(a)}(u_{1:a})
  \prod_{k=1}^a h(u_k) \sum_{\substack{s_1,\ldots ,s_m=1\\s_1\neq\cdots\neq s_m}}^a \prod_{l=1}^m\tilde{l}_d(z_l|u_{s_l})\nu(du_{1:a})
  \\
  \label{spec}
  &=\sum_{a=m}^\infty\frac{q_d^{a-m}}{(a-m)!}\int_{\Lambda^a}j_\Phi^{(a)}(u_{1:a})
  \prod_{k=1}^a h(u_k)
  \prod_{l=1}^m\tilde{l}_d(z_l|u_l)\nu(du_{1:a}).
\end{align}
Next, substituting \eqref{spec} and the relation
$$
  \frac{\partial_{\delta_{z_{Q(1)}}}}{\partial g}\cdots\frac{\partial_{\delta_{z_{Q(m-\pi)}}}}{\partial g}\mathcal{G}_{\Xi_c}(g)_{|g=0}
   =
  j_{\Xi_c}^{(m-\pi)}(z_{Q(1)},...,z_{Q(m-\pi)})
  =\prod_{k=1}^{m-\pi}l_c(z_{Q(k)})
  $$ 
 into \eqref{bruno1}, we obtain 
 \begin{align}
   \label{eqeq} 
  & \frac{\partial_{\delta_{z_1}}}{\partial g}\cdots\frac{\partial_{\delta_{z_m}}}{\partial g}\mathcal{G}_{\Phi,\Xi}(h,g)_{|g=0}
  \\
  \nonumber
  &=\sum_{a=0}^\infty\sum_{\substack{S\subset\{1,\ldots ,m\}
      \\ \left|S\right|\leq a}}\prod_{j\notin S}l_c(z_j)\frac{q_d^{a-\left|S\right|}}{(a-\left|S\right|)!} \int_{\Lambda^a}j_\Phi^{(a)}(u_{1:a})\prod_{k=1}^a h(u_k)\prod_{i\in S}\tilde{l}_d(z_i|u_i)\nu(du_{1:a}). 
\end{align}
\noindent
 Hence we find 
\begin{align*}
  & \frac{\partial_{\eta_1}}{\partial h}\cdots\frac{\partial_{\eta_n}}{\partial h}\frac{\partial_{\delta_{z_1}}}{\partial g}\cdots\frac{\partial_{\delta_{z_m}}}{\partial g}\mathcal{G}_{\Phi,\Xi}(h,g)_{|h=0,g=0}
  \\
  &=\frac{\partial_{\eta_1}}{\partial h}\cdots\frac{\partial_{\eta_n}}{\partial h}\sum_{a=0}^\infty\sum_{\substack{S\subset\{1,\ldots ,m\}
      \\ \left|S\right|\leq a}}\prod_{j\notin S}l_c(z_j)\frac{q_d^{a-\left|S\right|}}{(a-\left|S\right|)!} \int_{\Lambda^a}j_\Phi^{(a)}(u_{1:a})\prod_{k=1}^a h(u_k)\prod_{i\in S}\tilde{l}_d(z_i|u_i)\nu(du_{1:a})_{|h=0}
\\
&=
\sum_{\substack{S\subset\{1,\ldots ,m\}\\ \left|S\right|\leq n}}\prod_{j\notin S}l_c(z_j)\frac{q_d^{n-\left|S\right|}}{(n-\left|S\right|)!}\int_{\Lambda^n}j_\Phi^{(n)}(u_{1:n})
\sum_{\substack{k_1,\ldots ,k_n=1\\k_1\neq\cdots \neq k_n}}^n\eta_1(u_{k_1})\cdots\eta_n(u_{k_n})
\prod_{i\in S}\tilde{l}_d(z_i|u_i)\nu(du_{1:n}), 
\end{align*}\noindent
after setting $h=0$. 
By substituting $\eta_1,\ldots ,\eta_n$ with Dirac delta functions $\delta_{x_1},\ldots ,\delta_{x_n}$ at distinct configuration points $x_{1:n}\in\Lambda$, 
 the $(n,m)$-th joint Janossy density of $(\Phi,\Xi)$ is then given by
\begin{align*}
j_{\Phi,\Xi=z_{1:m}}^{(n,m)}(x_{1:n}) & =\frac{\partial_{\delta_{x_1}}}{\partial h}\cdots\frac{\partial_{\delta_{x_n}}}{\partial h}\frac{\partial_{\delta_{z_1}}}{\partial g}\cdots\frac{\partial_{\delta_{z_m}}}{\partial g}\mathcal{G}_{\Phi,\Xi}(h,g)_{|h=0,g=0}
\\
&=j_\Phi^{(n)}(x_{1:n})\sum_{\substack{S\subset\{1,\ldots ,m\}\\ \left|S\right|\leq n}}\frac{n!q_d^{n-\left|S\right|}}{(n-\left|S\right|)!}\prod_{j\notin S}l_c(z_j)\sum_{\pi:S\rightarrow\{1,\ldots ,n\}}\prod_{i\in S}\tilde{l}_d(z_i|x_{\pi(i)}), 
\end{align*}
 which shows \eqref{x1}. 
By \eqref{jnsy}, \eqref{mspt}, \eqref{mgnal} and \eqref{eqeq}, we have 
\begin{align*}
 j_\Xi^{(m)}(z_{1:m}) & =\frac{\partial_{\delta_{z_1}}}{\partial g}\cdots\frac{\partial_{\delta_{z_m}}}{\partial g}\mathcal{G}_{\Xi}(g)_{|g=0}\\
&=\frac{\partial_{\delta_{z_1}}}{\partial g}\cdots\frac{\partial_{\delta_{z_m}}}{\partial g}\mathcal{G}_{\Phi,\Xi}(\mathbf{1},g)_{|g=0}\\
&=\sum_{n=0}^\infty\sum_{\substack{S\subset\{1,\ldots ,m\}\\ \left|S\right|\leq n}}\prod_{j\notin S}l_c(z_j)\frac{q_d^{n-\left|S\right|}}{(n-\left|S\right|)!}\int_{\Lambda^n}j_\Phi^{(n)}(u_{1:n})\prod_{i\in S}\tilde{l}_d(z_i|u_i)\nu(du_{1:n}), 
\end{align*}
 which shows \eqref{jd}.
Finally, \eqref{x012} follows from the Bayes formula.
\end{Proof}
 The combinatorics of Lemma~\ref{ldjkldsf} is similar to
 Theorem~1 of \cite{lund}, which instead computes the conditional likelihood
 $j^{(m)}_{\Xi \mid \Phi = x_{1:n}} (z_1, \ldots , z_m)$
 of the observed point process $\Xi$ given a Poisson point process $\Phi$.
\\
 
Note that
\eqref{x1}
and 
\eqref{jd}
admit natural combinatorial interpretations by
identifying $S^c = \{ 1,\ldots ,m \} \setminus S$ to the set of
  points created according to the Poisson point process with intensity
  function $l_c(z)$, and by letting $n-|S|$ denote the number of 
  points in $\Phi$ deleted with probability $q_d$
  by the Bernoulli point process $\Xi_s$.
\subsubsection*{Poisson case}
In case $\Phi$ is the Poisson point process with intensity measure $\nu (dx)$
we have $j^{(n)}_\Phi=e^{ - \nu ( \Lambda) }$, $n\geq 0$, hence \eqref{jd} recovers
the classical expression 
\begin{align*} 
\nonumber 
  j^{(m)}_\Xi (z_1,\ldots , z_m)
& =  
e^{ - \nu ( \Lambda) }
\sum_{n\geq 0}
\sum_{S\subset \{1,\ldots , m\} \atop |S|\leq n } 
  \prod_{j \notin S} l_c(z_j) 
     \frac{  (  q_d \nu ( \Lambda ) )^{n-|S|}}{(n-|S|)!} 
    \prod_{i \in S} \int_\Lambda \tilde{l}_d (z_i | x_i) \nu (dx_i) 
\\
&
{
  =  
e^{ - \nu ( \Lambda) }
\sum_{S\subset \{1,\ldots , m\} } 
  \prod_{j \notin S} l_c(z_j) 
\prod_{i \in S} \int_\Lambda \tilde{l}_d (z_i | x_i) \nu (dx_i)
\sum_{n\geq |S|}
     \frac{  (  q_d \nu ( \Lambda ) )^{n-|S|}}{(n-|S|)!} 
}
\\
&
{
  =  
e^{ - p_d \nu ( \Lambda) }
\sum_{S\subset \{1,\ldots , m\} } 
\prod_{j \notin S} l_c(z_j) 
\prod_{i \in S} \int_\Lambda \tilde{l}_d (z_i | x_i) \nu (dx_i)
}
\\
& {
  =  
e^{ - p_d \nu ( \Lambda) }
\prod_{j=1}^m
\left(
l_c(z_j) 
+ \int_\Lambda \tilde{l}_d (z_j | x ) \nu (dx)
\right),
\qquad
 m\geq 0.
}
\end{align*} 
\subsubsection*{First-order posterior moment}
 In the next proposition we express the first-order conditional moment 
 of $\Phi$ given the sensor measurements $\Xi =z_{1:m}=(z_1,\ldots , z_m)$, 
 using extensions of the corrector terms
 $l^{(1)}_{z_{1:m}}$ introduced in 
 \cite{duhc}
 for the cardinalized PHD filter, see Equation~(19) in Lemma~1 therein, 
 and also Equation~(41) in Theorem~IV.7 of \cite{sdhc}
 for the Panjer-based PHD filter. 
\begin{prop}
\label{p1}
The first-order conditional moment of $\Phi$ given that $\Xi =(z_1,\ldots , z_m)$ 
is given by its density 
\begin{equation}
  \label{rho11} 
  \mu^{(1)}_{\Phi \mid \Xi =z_{1:m}} (x) 
  = 
  q_d l^{(1)}_{z_{1:m}} ( x )
  +
  \sum_{z\in z_{1:m}}
  \tilde{l}_d(x\mid z)
  l^{(1)}_{z_{1:m}} (x;z)
  ,
 \end{equation} 
with respect to $\nu (dx)$, $m\geq 0$, where
\begin{equation}
  \label{l1corr} 
    l^{(1)}_{z_{1:m}} ( x )
: =     \frac{\Upsilon^{(1)}_{z_{1:m}} ( x)
    }{j^{(m)}_\Xi (z_{1:m})},
\qquad
      l^{(1)}_{z_{1:m}} (x;z): = 
 \frac{
   \Upsilon^{(1)}_{z_{1:m}\! \setminus z} (x)
 }{j^{(m)}_\Xi (z_{1:m})}
 ,
\end{equation} 
 are corrector terms,
 $j^{(m)}_\Xi (z_{1:m}) = j^{(m)}_\Xi (z_1,\ldots , z_m)$ is given by
  \eqref{jd},  and
\begin{equation} 
  \label{djksds1}
      \Upsilon^{(1)}_{z_{1:m}} (x)
 : =
    \sum_{p\geq 0}
     \displaystyle \sum_{S\subset \{1,\ldots , m\} \atop |S| \leq p} 
     \frac{q_d^{p-|S|}}{(p-|S|)!}
  \prod_{j \notin S} l_c(z_j) 
     \int_{\Lambda^p}
    j^{(p+1)}_\Phi (x_{1:p} , x ) 
    \prod_{i \in S} \tilde{l}_d (z_i | x_i ) 
    \nu ( dx_{1:p} ), 
\end{equation} 
$m\geq 0$. 
\end{prop}
\begin{Proof}
The first-order joint moment density of
$\Phi$ with $\Xi =(z_1,\ldots , z_m)$
  can be obtained from the PGFl \eqref{mspt} as 
  $$
      \mu^{(1)}_{\Phi , \Xi =z_{1:m}} (x) = 
  \frac{\partial_{\delta_x}}{\partial h} {\cal G}_{\Phi , \Xi =z_{1:m}} (h)_{\mid h=1} 
    , 
    $$ 
  or, using the joint Janossy densities \eqref{x1} and
  denoting by $d\hat{x}_r$ the absence of $dx_r$, as 
\begin{align} 
\nonumber 
&   \mu^{(1)}_{\Phi, \Xi =z_{1:m}} (x) 
=
     \sum_{p\geq 1}     \frac{1}{p!}
     \sum_{r=1}^p
     \int_{\Lambda^{p-1}}
     j^{(p,m)}_{\Phi, \Xi = z_{1:m} } (x_{1:p})_{\mid x_r = x}
     \nu ( dx_1) \cdots \nu ( d\hat{x}_r) \cdots \nu ( dx_p)
\\
     \nonumber
     &
     {
      = 
     \sum_{p\geq 1}
     \sum_{r=1}^p
 \sum_{S\subset \{1,\ldots , m\} \atop |S| \leq p} 
 \frac{  q_d^{p-|S|}}{(p-|S|)!}
     \prod_{j \notin S} l_c(z_j) 
 \int_{\Lambda^{p-1}}
   j^{(p)}_\Phi (x_{1:p})_{\mid x_r = x}
 \sum_{\pi : S \rightarrow \{1,\ldots , p\} } 
    \prod_{i \in S} \tilde{l}_d (z_i | x_{\pi (i) } )_{\mid x_r = x} 
     \nu ( dx_{1:p} \! \setminus dx_r) 
     }
     \\
     \nonumber
     & {
       = 
     \sum_{p\geq 1}
     \sum_{r=1}^p
     \sum_{S\subset \{1,\ldots , m\} \atop |S| \leq p-1} 
     \frac{ q_d^{p-|S|} }{(p-|S|)!}
  \prod_{j \notin S} l_c(z_j) 
        \sum_{\pi : S \rightarrow \{1,\ldots , p\} \setminus \{r\}} 
   \int_{\Lambda^{p-1}}
   j^{(p)}_\Phi (x_{1:p})_{\mid x_r = x}
    \prod_{i \in S} \tilde{l}_d (z_i | x_{\pi (i)} )
     \nu ( dx_{1:p} \! \setminus dx_r) 
     }
     \\
     \nonumber
     &  {
       + 
     \sum_{p\geq 1}
     \sum_{r=1}^m
    \sum_{S\subset \{1,\ldots , m\} \atop |S| \leq p, r\in S} 
     \frac{  q_d^{p-|S|}}{(p-|S|)!} \prod_{j \notin S} l_c(z_j) 
  \sum_{\pi : S \rightarrow \{1,\ldots , p\} } 
   \int_{\Lambda^{p-1}}
   \prod_{i \in S} \tilde{l}_d (z_i | x_i)_{\mid x_{\pi (r)} = x}
    j^{(p)}_\Phi (x_{1:p})_{\mid x_{\pi (r)} = x}
     \nu ( dx_{1:p} \! \setminus dx_{\pi (r) }) 
     }
     \\
     \nonumber
     & {
       = 
     \sum_{p\geq 1}
     \sum_{r=1}^p
    \sum_{S\subset \{1,\ldots , m\} \atop |S| \leq p-1} 
     \frac{ q_d^{p-|S|}}{(p-|S|)!}   \prod_{j \notin S} l_c(z_j) 
        \sum_{\pi : S \rightarrow \{1,\ldots , p\} \setminus \{r\}} 
   \int_{\Lambda^{p-1}}
   \prod_{i \in S} \tilde{l}_d (z_i | x_{\pi (i)} )
    j^{(p)}_\Phi (x_{1:p})_{\mid x_r = x}
     \nu ( dx_{1:p} \! \setminus dx_r) 
     }
     \\
     \nonumber
     &  {
       + 
     \sum_{p\geq 1}
     \sum_{r=1}^m
     \tilde{l}_d (z_r|x)
     \sum_{S\subset \{1,\ldots , m\} \atop |S| \leq p, r\in S} 
     \frac{     q_d^{p-|S|}}{(p-|S|)!} \prod_{j \notin S} l_c(z_j) 
        }
     \\
     & \qquad
       {
         \times \sum_{\pi : S \rightarrow \{1,\ldots , p\} } 
   \int_{\Lambda^{p-1}}
   \prod_{i \in S \setminus \{ r\}} \tilde{l}_d (z_i | x_i)
    j^{(p)}_\Phi (x_{1:p})_{\mid x_{\pi (r)} = x}
     \nu ( dx_{1:p} \! \setminus dx_{\pi (r) } ) 
     }
     \\
     \nonumber
     & {
       = 
     q_d
     \sum_{p\geq 1}
     \sum_{S\subset \{1,\ldots , m\} \atop |S| \leq p-1} 
     \frac{ q_d^{p-1-|S|}}{(p-|S|-1)!} \prod_{j \notin S} l_c(z_j) 
   \sum_{\pi : S \rightarrow \{1,\ldots , p-1\} } 
   \int_{\Lambda^{p-1}}
   \prod_{i \in S} \tilde{l}_d (z_i | x_{\pi (i)} )
    j^{(p)}_\Phi (x_{1:p-1},x)
     \nu ( dx_{1:p-1})
     }
     \\
     \nonumber
     &  {
       + 
          \sum_{r=1}^m
          \sum_{p\geq 1}
          \tilde{l}_d (z_r|x)
     \sum_{S\subset \{1,\ldots , m\} \atop |S \setminus \{ r \} | \leq p-1, r\in S} 
     \frac{  q_d^{p+1-| S \setminus \{r\}|}}{(p-|S|-1)!} \prod_{j \notin S} l_c(z_j) 
        }
     \\
     \nonumber
      & \qquad \times \sum_{\pi : S \setminus \{r\} \rightarrow \{1,\ldots , p-1\} } 
   \int_{\Lambda^{p-1}}
   \prod_{i \in S \setminus \{ r\}} \tilde{l}_d (z_i | x_{\pi (i)} )
    j^{(p)}_\Phi (x_{1:p-1},x) \nu ( dx_{1:p-1} )
     \\
     \nonumber
     & {
       = 
          q_d    \Upsilon^{(1)}_{z_{1:m}} (x)
  +
  \sum_{z\in z_{1:m}}
  \tilde{l}_d(x\mid z)
\Upsilon^{(1)}_{z_{1:m}\! \setminus z} (x)
,
     }
\end{align} 
and it remains to divide by $j^{(m)}_\Xi (z_{1:m})$.
\end{Proof}
\subsubsection*{Second-order posterior moment} 
 Similarly, the second partial moment of 
the first-order integral of
$\Phi$ when $\Xi = z_{1:m} = (z_1,\ldots , z_m)$ 
is obtained in the next
proposition,
which uses an extension of the corrector terms
 $l^{(2)}_{z_{1:m}}$ introduced in 
\cite{duhc} for the cardinalized PHD filter, see
Equation~(29) in Lemma~2 therein, and also  
 Equation~(42) in Theorem~IV.8 of \cite{sdhc}
 for the Panjer-based PHD filter. 
 \begin{prop}
\label{p2}
The second-order conditional factorial moment of 
$\Phi$ given that $\Xi = (z_1,\ldots , z_m)$ 
is given by its density 
\begin{eqnarray}
  \label{rho22} 
    \rho^{(2)}_{\Phi \mid \Xi =z_{1:m}} (x,y)  
 & = &       q_d^2
l^{(2)}_{z_{1:m}} (x,y) 
  + q_d 
 \sum_{z\in z_{1:m}} 
 \big( \tilde{l}_d (z | x) + \tilde{l}_d (z | y) \big) 
     l^{(2)}_{z_{1:m}} (x,y;z)
\\
\nonumber 
  & &
  {
    + 
 \sum_{z,z'\in z_{1:m} \atop z \not= z'} 
 \tilde{l}_d (z | x)
 \tilde{l}_d (z' | y)
 l^{(2)}_{z_{1:m}} (x,y;z,z')
 ,
 \quad x,y\in \Lambda, \quad x\not= y,
  }
\end{eqnarray} 
 with respect to $\nu (dx)\nu (dy)$, with the corrector terms 
\begin{equation}
  \label{l2corr} 
  {
    l^{(2)}_{z_{1:m}} (x,y) 
  : =   \frac{\Upsilon^{(2)}_{z_{1:m}} (x,y)
    }{j^{(m)}_\Xi (z_{1:m})}
,
\qquad
     l^{(2)}_{z_{1:m}} (x,y;z)
: = 
 \frac{
   \Upsilon^{(2)}_{z_{1:m} \! \setminus z} (x,y)
 }{j^{(m)}_\Xi (z_{1:m})}
 ,
  }
  \end{equation}
 and 
\begin{equation}
  \label{l2corr.2} 
  {
 l^{(2)}_{z_{1:m}} (x,y;z,z')
: = 
 \frac{
      \Upsilon^{(2)}_{z_{1:m} \! \setminus \{z,z'\} } (x,y)
    }{j^{(m)}_\Xi (z_{1:m})}
 ,
  }
  \end{equation}
 where $j^{(m)}_\Xi (z_1,\ldots , z_m)$ is as in \eqref{jd},
 and 
\begin{equation} 
  \label{djksds}
      \Upsilon^{(2)}_{z_{1:m}} (x,y) 
 : =
    \sum_{p\geq 0} 
 \displaystyle \sum_{S\subset \{1,\ldots , m\} \atop |S| \leq p} 
 \frac{q_d^{p-|S|}}{(p-|S|)!}
  \prod_{j \notin S} l_c(z_j)
      \int_{\Lambda^p}
           j^{(p+2)}_\Phi (x_{1:p}, x , y ) 
    \prod_{i \in S } \tilde{l}_d (z_i | x_i ) 
\nu ( dx_{1:p})
,
\end{equation} 
$x,y\in \Lambda$, $m \geq 0$.
\end{prop}
\begin{Proof}
  Factorial moments can be computed
  using the second derivative of the conditional PGFl
 \eqref{mspt},
 see \eqref{**}, 
 or equivalently using the joint Janossy densities \eqref{x1} as 
 in the proof of Proposition~\ref{p1}.
 We have 
\begin{align} 
\nonumber 
  &
{
  \rho^{(2)}_{\Phi , \Xi =z_{1:m}} (x,y)
  = 
     \sum_{p\geq 0}
     \frac{1}{p!}
 \sum_{r,u=1 \atop r\not= u}^p 
  \int_{\Lambda^{p-2}}
     j^{(p,m)}_{\Phi, \Xi = z_{1:m} } (x_{1:p} )_{\mid x_r = x \atop \mid x_u = y}
     \nu ( dx_1) \cdots \nu ( d\hat{x}_r) \cdots \nu ( d\hat{x}_u) \cdots \nu ( dx_p)
}
\\
\nonumber
& {
  = 
     \sum_{p\geq 2}
    \sum_{r,u=1\atop r\not= u}^p 
  \displaystyle \sum_{S\subset \{1,\ldots , m\} \atop |S| \leq p} 
     \frac{ q_d^{p-|S|}}{(p-|S|)!}  \prod_{j \notin S} l_c(z_j) 
 j^{(p)}_\Phi (x_{1:p})_{\mid x_r = x\atop \mid x_u = y} 
}
\\
\nonumber
& 
{
  \qquad
  \times
       \!\!     \!\!     \!\!
       \sum_{\pi : S \rightarrow \{1,\ldots , p \} } 
          \int_{\Lambda^{p-1}}
          \prod_{i \in S} \tilde{l}_d (z_i | x_{\pi (i) })_{\mid x_r = x\atop \mid x_u = y}
          \nu ( dx_1) \cdots \nu ( d\hat{x}_r) \cdots \nu ( d\hat{x}_u) \cdots \nu ( dx_p)
}
\\
     \nonumber
     &
     {
       = 
 q_d^2 \sum_{p\geq 0}
 \displaystyle \sum_{S\subset \{1,\ldots , m\} \atop |S| \leq p} 
 \frac{  q_d^{p-|S|}}{(p-|S|)!}
 \prod_{j \notin S} l_c(z_j) 
         \sum_{\pi : S \rightarrow \{1,\ldots , p \} } 
          \int_{\Lambda^p}
          j^{(p+2)}_\Phi (x_{1:p}, x, y ) 
    \prod_{i \in S} \tilde{l}_d (z_i | x_{\pi (i) }) 
\nu ( dx_{1:p})
     }
     \\
     \nonumber
     &
     {
       \quad + q_d 
 \sum_{r=1}^m 
 \tilde{l}_d (z_r | x)
    \sum_{p\geq 0}
    \displaystyle \sum_{S\subset \{1,\ldots , m\}
      \atop |S| \leq p+1, r\in S} 
    \frac{  q_d^{p-|S|}}{(p-|S|)!}
  \prod_{j \notin S} l_c(z_j) 
     }
     \\
  \nonumber
  & {
    \quad \quad
    \times
         \!\!         \!\!
        \sum_{\pi : S \setminus \{ r \} \rightarrow \{1,\ldots , p \} } 
          \int_{\Lambda^p}
          j^{(p+2)}_\Phi (x_{1:p} , x , y ) 
    \prod_{i \in S \setminus \{ r \} } \tilde{l}_d (z_i | x_{\pi (i) })  \nu ( dx_{1:p}) 
  }
  \\
     \nonumber
     &
     {
       \quad + q_d 
 \sum_{r=1}^m 
 \tilde{l}_d (z_r | y)
    \sum_{p\geq 0}
  \displaystyle \sum_{S\subset \{1,\ldots , m\} \atop |S| \leq p+1, r \in S} 
    \frac{  q_d^{p-|S|}}{(p-|S|)!}
    \prod_{j \notin S} l_c(z_j)
     }
     \\
  \nonumber
  & {
    \quad \quad\times
\!\!
\!\!
\!\!
         \sum_{\pi : S \setminus \{ r \} \rightarrow \{1,\ldots , p \} } 
  \int_{\Lambda^p}
  j^{(p+2)}_\Phi (x_{1:p} , x , y ) 
    \prod_{i \in S \setminus \{ r \} } \tilde{l}_d (z_i | x_{\pi (i) })  \nu ( dx_{1:p}) 
  }
  \\
     \nonumber
     &
     {
       \quad + 
 \sum_{r,u=1 \atop r \not= u}^m 
 \tilde{l}_d (z_r | x)
 \tilde{l}_d (z_u | y)
    \sum_{p\geq 0}
    \displaystyle \sum_{S\subset \{1,\ldots , m\} \atop |S| \leq p+2, r,u\in S} 
    \frac{  q_d^{p-|S|}}{(p-|S|)!}
  \prod_{j \notin S} l_c(z_j)
     }
     \\
  \nonumber
  &
  {
    \quad \quad\times
   \!\!   \!\!   \!\!
        \sum_{\pi : S \setminus \{r,u\} \rightarrow \{1,\ldots , p \} }
    \int_{\Lambda^p}
    j^{(p+2)}_\Phi (x_{1:p} , x , y ) 
    \prod_{i \in S \setminus \{ r ,u\} } \tilde{l}_d (z_i | x_{\pi (i) }) 
\nu ( dx_{1:p}). 
  }
\\
     \nonumber
     &
     {
 =       q_d^2
 \Upsilon^{(2)}_{z_{1:m}} (x,y)
  + q_d 
 \sum_{z\in z_{1:m}} 
 \big( \tilde{l}_d (z | x) + \tilde{l}_d (z | y) \big) 
 \Upsilon^{(2)}_{z_{1:m}\setminus z} (x,y) 
    + 
 \sum_{z,z'\in z_{1:m} \atop z \not= z'} 
 \tilde{l}_d (z | x)
 \tilde{l}_d (z' | y)
 \Upsilon^{(2)}_{z_{1:m}\setminus \{z,z'\} } (x,y) 
 ,
  }
\end{align} 
 $x,y\in \Lambda$, $x\not= y$, 
 and it remains to divide by $j^{(m)}_\Xi (z_{1:m})$.
\end{Proof}
\subsubsection*{Poisson case}
In the case of a Poisson point process 
with $j^{(n)}_\Phi=e^{ - \nu ( \Lambda) }$, $n\geq 0$, 
\eqref{djksds1} reads 
\begin{eqnarray*} 
      \Upsilon^{(1)}_{z_{1:m}} (x)
  & = & 
  {
    e^{ - \nu ( \Lambda) }
     \sum_{p\geq 0}     \displaystyle \sum_{S\subset \{1,\ldots , m\} \atop |S| \leq p} 
     \frac{q_d^{p-|S|}}{(p-|S|)!}
   \prod_{j \notin S} l_c(z_j) 
   \nu ( \Lambda)^{p-|S|} 
     \prod_{i \in S}
   \int_{\Lambda}
    \tilde{l}_d (z_i | u ) 
    \nu ( du )
  }
  \\
 & = & 
  {
    e^{ - \nu ( \Lambda) }
     \displaystyle \sum_{S\subset \{1,\ldots , m\}} 
     \prod_{j \notin S} l_c(z_j)
     \prod_{i \in S}
   \int_{\Lambda}
    \tilde{l}_d (z_i | u ) 
    \nu ( du )
     \sum_{p\geq |S|}
     \frac{q_d^{p-|S|}}{(p-|S|)!}
   \nu ( \Lambda)^{p-|S|} 
  }
  \\
 & = & 
  {
    e^{ - p_d \nu ( \Lambda) }
    \displaystyle \sum_{S\subset \{1,\ldots , m\}} 
     \prod_{j \notin S} l_c(z_j)
     \prod_{i \in S}
   \int_{\Lambda}
    \tilde{l}_d (z_i | u ) 
    \nu ( du )
  }
  \\
    & = &
  {
    e^{ - p_d \nu ( \Lambda) }
\prod_{j=1}^m
\left(
l_c(z_j) 
+ \int_\Lambda \tilde{l}_d (z_j | u ) \nu (du)
\right)
  }
  \\
& = &
  {
    j^{(m)}_\Xi (z_{1:m}),
\qquad
m\geq 0, 
  }
\end{eqnarray*} 
 and similarly from \eqref{djksds} we find 
$$ 
\Upsilon^{(2)}_{z_{1:m}} (x,y) 
=
 \Upsilon^{(1)}_{z_{1:m}} (x)
  =
j^{(m)}_\Xi (z_{1:m})
. 
$$   
Hence, in the Poisson case 
the corrector terms are given by
$l^{(1)}_{z_{1:m}} ( x ) = l^{(2)}_{z_{1:m}} (x,y) = 1$ and 
$$ 
       l^{(1)}_{z_{1:m}} (x;z) = 
     l^{(2)}_{z_{1:m}} (x,y;z)
 = 
 \frac{1}{
l_c(z) 
+ \int_\Lambda \tilde{l}_d (z | u ) \nu (du)
 }
 ,
 $$
 with 
\begin{eqnarray*} 
 l^{(2)}_{z_{1:m}} (x,y;z,z')
 & = &  
  l^{(1)}_{z_{1:m}} (x;z)
  l^{(1)}_{z_{1:m}} (y,z')
  \\
   & = & 
  \frac{1}{
    \big(
    l_c(z) 
+ \int_\Lambda \tilde{l}_d (z | u ) \nu (du)
\big)
  \big( 
l_c(z') 
+ \int_\Lambda \tilde{l}_d (z' | u ) \nu (du)
\big)
  }
  ,
\end{eqnarray*} 
and the first and second (factorial) moment densities
 \eqref{rho11}, \eqref{rho22}
 of $\Phi$ with respect to $\nu (dx)$
 given the point process $\Xi$ recover the classical expressions  
$$
 \mu^{(1)}_{\Phi \mid \Xi =z_{1:m}} (x) 
 = q_d + \sum_{z\in z_{1:m}}
 \frac{ \tilde{l}_d (z | x)}{ l_c (z ) + \int_{\Lambda} \tilde{l}_d (z | u ) \nu ( du )  }, 
$$ 
 of first order moment density,
 see Relation~(2.87) in \cite{cdh}, and
 the second-order moment density 
 \begin{eqnarray*} 
          \rho^{(2)}_{\Phi \mid \Xi =z_{1:m}} (x,y)  
      & = &
    q_d^2
  + q_d 
 \sum_{z\in z_{1:m}} 
 \frac{\tilde{l}_d (z | x)+\tilde{l}_d (z | y)}{
        l_c(z)
   +
  \int_\Lambda \tilde{l}_d (z | u) \nu ( du)
 }
              \\
     \nonumber
     & &
     { 
     + 
 \sum_{r,p=1 \atop r \not= p}^m 
 \frac{  \tilde{l}_d (z_r | x) \tilde{l}_d (z_p | y)
}{
     \big(
  l_c(z_r)
  +
  \int_\Lambda \tilde{l}_d (z_r | u) \nu ( du)
  \big)
  \big(
  l_c(z_p)
  +
  \int_\Lambda \tilde{l}_d (z_p | v) \nu ( dv)
  \big)
 },
 }
\end{eqnarray*} 
    $x,y\in \Lambda$, $x\not= y$,
    $m\geq 0$. 
    See, e.g., Proposition~V.1(a) of \cite{sdhc}
     and Exercise~4.3.4 in \cite{cdh}. 
\subsubsection*{Posterior covariance}
 For $A, B$ measurable subsets of $\real^d$, let
$$ 
{
  c^{(2)}_{\Phi \mid \Xi =(z_{1:m})} \! (A,B)  
  := 
  \mu^{(2)}_{\Phi \mid \Xi =z_{1:m}} \! (A,B)
  -
   \mu^{(1)}_{\Phi \mid \Xi =z_{1:m}} \! (A) 
   \mu^{(1)}_{\Phi \mid \Xi =z_{1:m}} \! (B),
}
$$
   denote the posterior covariance,
   where $\mu^{(1)}_{\Phi \mid \Xi =z_{1:m}} (A)$ is the posterior first order moment 
   $$
{
     \mu^{(1)}_{\Phi \mid \Xi =z_{1:m}} (A) 
:=
 \int_A \mu^{(1)}_{\Phi \mid \Xi =z_{1:m}} (x) \nu (dx),
   }
   $$
 and $\mu^{(2)}_{\Phi \mid \Xi =z_{1:m}} (A,B)$ is the posterior second-order moment 
 $$
 \mu^{(2)}_{\Phi \mid \Xi =z_{1:m}} (A,B)
 = 
      \int_{A\cap B} \mu^{(1)}_{\Phi \mid \Xi =z_{1:m}} (x) \nu (dx ) 
      + 
 \int_{A\times B} \rho^{(2)}_{\Phi \mid \Xi =z_{1:m}} (x,y) \nu (dx) \nu (dy). 
$$ 
Using Relations~\eqref{rho11} and \eqref{rho22}
in Propositions~\ref{p1}-\ref{p2}, 
we obtain the following 
representation of the posterior covariance.
\begin{prop}
  \label{djlksd}
  The posterior covariance
  $c^{(2)}_{\Phi \mid \Xi =z_{1:m}} (A,B)$
  of $\Phi$ given that $\Xi =(z_1,\ldots , z_m)$ 
is given by 
\begin{align} 
  \nonumber
  &
  {
    c^{(2)}_{\Phi \mid \Xi =z_{1:m}} (A,B)  
 = 
q_d
\int_{A\cap B}
         l^{(1)}_{z_{1:m}} ( x )
     \nu (dx)
    + q_d^2 \int_{A\times B}
    \left(
    l^{(2)}_{z_{1:m}} (x,y )
    -
    l^{(1)}_{z_{1:m}} ( x )
    l^{(1)}_{z_{1:m}} ( y )
    \right)
    \nu (dx) \nu (dy)
  }
  \\
\nonumber
&
  {
+ q_d \sum_{z\in z_{1:m}}
 \int_{A\times B} \tilde{l}_d (z | x) 
 \left(
 l^{(2)}_{z_{1:m}} (x,y;z)
 -
 l^{(2)}_{z_{1:m}} (x,y)
 l^{(1)}_{z_{1:m}} (x;z)
 \right) \nu ( dx ) \nu (dy) 
  }
  \\
\nonumber
&
{
 + q_d \sum_{z\in z_{1:m}} 
 \int_{A\times B}
 \tilde{l}_d (z | y) 
 \left(
 l^{(2)}_{z_{1:m}} (x,y;z)
 -
 l^{(1)}_{z_{1:m}} (y;z)
 l^{(2)}_{z_{1:m}} (x,y)
 \right)
    \nu ( dx ) \nu (dy)
}
\\
\nonumber
&  {
   + 
\sum_{z\in z_{1:m}} \left(
\int_{A\cap B}\tilde{l}_d(z|x)
l^{(1)}_{z_{1:m}} (x;z)
\nu ( dx)
    - 
    \int_A \tilde{l}_d(z|x)
    l^{(1)}_{z_{1:m}} (x;z)
    \nu (dx)
    \int_B \tilde{l}_d(z|y)
    l^{(1)}_{z_{1:m}} (y;z)
    \nu (dy)
    \right)
}
\\
\label{fklds323} 
&
{
 + 
    \sum_{z,z'\in z_{1:m} \atop z\not= z'} 
    \int_{A\times B}
    \tilde{l}_d(z|x) \tilde{l}_d(z'|y)
    \left(
     l^{(2)}_{z_{1:m}} (x,y;z,z')
    -
     l^{(1)}_{z_{1:m}} (x;z)
     l^{(1)}_{z_{1:m}} (y;z')
\right)
    \nu (dx) \nu ( dy) ,
}
\end{align} 
$m\geq 0$.
\end{prop}
When $A=B$, Relation~\eqref{fklds323} becomes the variance identity 
\begin{eqnarray*} 
  \nonumber
    c^{(2)}_{\Phi \mid \Xi =z_{1:m}} (A,A)  
& = & 
q_d
\int_A
         l^{(1)}_{z_{1:m}} ( x )
     \nu (dx)
   + 
\sum_{z\in z_{1:m}} 
\int_A \tilde{l}_d(z|x)
l^{(1)}_{z_{1:m}} (x;z)
\nu ( dx)
  \\
\nonumber
& & 
     + q_d^2
     \left(
     \int_{A^2} 
    l^{(2)}_{z_{1:m}} (x,y )
    \nu (dx) \nu (dy)
    -
    \left(
    \int_A l^{(1)}_{z_{1:m}} ( x )
    \nu (dx) 
    \right)^2 
    \right)
\\
\nonumber
& & 
+ 2 q_d \sum_{z\in z_{1:m}}
 \int_{A^2} \tilde{l}_d (z | x) 
 \left(
 l^{(2)}_{z_{1:m}} (x,y;z)
 -
 l^{(2)}_{z_{1:m}} (x,y)
 l^{(1)}_{z_{1:m}} (x;z)
 \right) \nu ( dx ) \nu (dy) 
\\
\nonumber 
& & 
 + \sum_{z,z'\in z_{1:m} \atop z\not= z'} 
    \int_{A^2}
    \tilde{l}_d(z|x) \tilde{l}_d(z'|y)
     l^{(2)}_{z_{1:m}} (x,y;z,z')
    \nu (dx) \nu ( dy)
\\
\nonumber 
& &  - \sum_{z,z'\in z_{1:m}} 
    \int_{A} \tilde{l}_d(z|x) l^{(1)}_{z_{1:m}} (x;z)
    \nu (dx) 
    \int_{A} \tilde{l}_d(z|x) l^{(1)}_{z_{1:m}} (x;z')
    \nu (dx),  
\end{eqnarray*} 
 which takes a form similar to the variance update formula
obtained for the Panjer-based PHD filter, 
see Equations~(41)-(42) of Theorem~IV.8 in \cite{sdhc}. 
\subsubsection*{Poisson case}
In the case of a Poisson point process 
 with $j^{(n)}_\Phi=e^{ - \nu ( \Lambda) }$, $n\geq 0$,
 Proposition~\ref{djlksd} recovers the covariance
\begin{eqnarray*} 
    {
      c^{(2)}_{\Phi \mid \Xi =z_{1:m}} (A,B)  
  }
 & = &
  {
    q_d \nu ( A\cap B ) 
 + \sum_{z\in z_{1:m}}
    \frac{ \int_{A\cap B} \tilde{l}_d (z | x) \nu ( dx ) }{
     l_c (z ) + \int_{\Lambda} \tilde{l}_d (z | u ) \nu ( du )  }
  }
  \\
  & &
  {
    - \sum_{z\in z_{1:m}} 
 \frac{  \int_A \tilde{l}_d (z | x) \nu ( dx ) \int_B \tilde{l}_d (z | y) \nu ( dy ) }{
  \big(
  l_c(z)
  +
  \int_\Lambda \tilde{l}_d (z | u) \nu ( du)
  \big)^2
}
  }
  ,
\end{eqnarray*} 
see, e.g., Equation~(41) and Proposition~V.1(a) of \cite{sdhc}, 
and the variance 
$$       c^{(2)}_{\Phi \mid \Xi =z_{1:m}} ( A , A )  
=
q_d \nu ( A ) 
    + \sum_{z\in z_{1:m}}
    \frac{ \int_A \tilde{l}_d (z | x) \nu ( dx ) }{
     l_c (z ) + \int_A \tilde{l}_d (z | u ) \nu ( du )  }
 \left(1 - 
 \frac{  \int_A \tilde{l}_d (z | x) \nu ( dx ) }{
   l_c(z) +
  \int_\Lambda \tilde{l}_d (z | u) \nu ( du)
  }\right),
$$ 
 see also Exercise~4.3.4 in \cite{cdh}. 
\section{Determinantal point processes} 
\label{s2.2}
In this section we review the properties of
determinantal point processes;
 see, e.g., \cite{dfpt2} and references therein for
additional background.  
\subsubsection*{Kernels and integral operators}
For any compact set $\Lambda\subseteq \real^d$, we denote by $L^2(\Lambda,\radon)$ the
Hilbert space of square-integrable functions {\em
w.r.t.} the restriction of the Radon measure $\radon$ on $\Lambda$,
equipped with the inner product
\begin{equation*}
\langle f,g\rangle_{L^2(\Lambda,\radon)}:=\int_{\Lambda}f(x) g(x) \,\radon(\d x),\quad\text{$f,g\in
L^2(\Lambda,\radon)$}. 
\end{equation*}
 By definition, an integral operator $\mathcal{K}:L^2(\underlying,\radon)\to L^2(\underlying,\radon)$ with kernel
$K:\underlying^2\to \real$ is a bounded operator
defined by 
\begin{equation*}
\mathcal{K}f(x):=\int_{\underlying}K(x,y)f(y)\,\radon(\d y),\quad\text{for
$\radon$-almost all $x\in \underlying$}.
\end{equation*}
It can be shown that $\mathcal{K}$ is a compact operator, 
which is self-adjoint if its kernel verifies
\begin{equation}
  \nonumber 
K(x,y)=K(y,x),\quad\text{for $\radon^{\otimes 2}$-almost
all $(x,y)\in \underlying^2$.}
\end{equation}
Equivalently, this means that the integral operator 
$\mathcal{K}$ is self-adjoint for any compact set $\Lambda\subseteq \real^d$.
If $\mathcal{K}$ is self-adjoint, by the spectral theorem
we have that $L^2(\Lambda,\radon)$ has
an orthonormal basis $(\varphi_n )_{n\geq 1}$ of eigenfunctions
of $\mathcal{K}$ with corresponding eigenvalues
$( \radon_n )_{n\geq 1}$, and the kernel $K$ of $\mathcal{K}$ can be written as
\begin{equation}
\label{eq:kdecomp}
K (x,y) = \sum_{n \ge 1} \radon_n \varphi_n (x) \varphi_n(y),
\qquad x, y \in \Lambda.
\end{equation}
For $\mathcal{K}$ a self-adjoint integral operator of trace class,
i.e. 
\begin{equation*}
\sum_{n\geq 1}|\radon_n|<\infty,
\end{equation*}
we
define the trace of $\mathcal{K}$ as $\mathrm{Tr\,} \mathcal{K} = \sum_{n\geq 1}\radon_n$.
Let also $\mathrm{Id}$ denote the identity operator on $L^2(\underlying,\radon)$ and let $\mathcal{K}$ be a trace class operator on $L^2(\underlying,\radon)$.
We define the Fredholm determinant of $\mathrm{Id}+\mathcal{K}$ as
\begin{equation}
\nonumber 
\Det(\mathrm{Id}+\mathcal{K}) = \mathrm{exp}\left(\sum_{n \ge 1}\frac{(-1)^{n-1}}{n} \Tr(\mathcal{K}^n)\right), 
\end{equation}
 with the relation 
\begin{equation}
\nonumber 
\Det(\mathrm{Id}+\mathcal{K})=\sum_{n \ge 0} \frac{1}{n!}
 \int_{\underlying^n} \mathrm{det}\big( K(x_i,x_j)_{1\le i,j \le n} \big)\, \radon(\d x_1)\cdots \radon(\d x_n) ,
\end{equation}
where 
$\mathrm{det}\big( K(x_i,x_j)_{1\le i,j \le n} \big)$ is the determinant of the $n\times n$ matrix
$(K(x_i,x_j))_{1\le i,j \le n}$, see Theorem~$2.4$ of \cite{shirai}, and 
also \cite{brezis2} for more details on Fredholm determinants. 

\subsubsection*{Determinantal point processes}

In the sequel we consider
a self-adjoint trace class operator $\mathcal{K}_\Psi$ 
on $L^2(\Lambda , \nu)$ with spectrum contained in $[0,1)$, and denote by
  $K_\Psi : \Lambda \times \Lambda \longrightarrow \real$ 
   the kernel of ${\mathcal K}_\Psi$.

\bigskip 

By the results in \cite{macchi} and
\cite{soshnikov} (see also Lemma 4.2.6 and Theorem 4.5.5 in
\cite{hough}) 
the determinantal point process
$\Psi$ on $\underlying$,
with integral operator ${\mathcal K}_\Psi $ is defined as in
\eqref{coor-in} by its correlation functions 
\begin{equation*}
\rho^{(n)}_\Psi (x_1,\ldots,x_n)=\mathrm{det}\big(K_\Psi (x_i,x_j)_{1\leq i,j\leq n}\big),
\end{equation*}
    {\em w.r.t.} the measure $\radon$ on $(\underlying,\underlyingsigma)$,
    $x_1,\ldots , x_n \in \Lambda$, with $x_i \not= x_j$, $1\leq i<j \leq n$, 
see also Lemma~3.3 of \cite{shirai}. 
In particular, we have 
\begin{equation}
  \label{kx}
  \mu^{(1)}_\Psi (x) =
  \rho^{(1)}_\Psi (x) = K_\Psi(x,x),
  \qquad x\in \Lambda, 
\end{equation} 
 and 
\begin{equation}
  \label{kx2}
 \rho^{(2)}_\Psi (x , y) = K_\Psi(x,x)K_\Psi(y,y)- ( K_\Psi (x,y) )^2, 
\end{equation} 
 $x,y\in \Lambda$, $x\not=y$, i.e. 
\begin{equation}
  \label{kxy} 
  \rho^{(2)}_\Psi (x , y) - \mu^{(1)}_\Psi (x) \mu^{(1)}_\Psi ( y) 
  = - ( K_\Psi(x,y) )^2 \leq 0, 
\quad x,y\in \Lambda, \quad x\not=y, 
\end{equation} 
 with $\rho^{(2)}_\Psi (x , x) := 0$, $x\in \Lambda$.  
The covariance of the determinantal point process $\Psi$ is then given by 
\begin{eqnarray} 
  \nonumber
      {
    c^{(2)}_\Psi (A,B)  
    }
    & = &
    { 
    \int_{A\cap B} \mu^{(1)}_\Psi (x) \mu (dx ) 
    + 
 \int_{A\times B}
 \big(
 \rho^{(2)}_\Psi (x,y)  
 -
 \mu^{(1)}_\Psi (x)
 \mu^{(1)}_\Psi (y) \big)
 \nu (dx) \nu (dy)
  }
  \\
  \label{fksdlfds0}
  & = &
  {
    \int_{A\cap B} K_\Psi (x,x) \nu (dx ) - \int_{A\times B} ( K_\Psi (x,y) )^2 \nu (dx) \nu (dy), 
  }
\end{eqnarray} 
which shows that the determinantal point process $\Psi$ is {\em negatively correlated}, 
since when $A\cap B = \emptyset$ we have
\begin{equation} 
  \label{fksdlfds}
  c^{(2)}_\Psi (A,B)  
 = - \int_{A\times B} ( K_\Psi (x,y) )^2 \nu (dx) \nu (dy) \leq 0.
\end{equation}
  
The interaction operator ${\mathcal J}_\Psi$ on $L^2(\underlying,\radon)$
is defined as 
\begin{equation}
  \label{idk}
  {\mathcal J}_\Psi : = ( \mathrm{Id} - {\mathcal K}_\Psi )^{-1} {\mathcal K}_\Psi, 
\end{equation} 
 and has the kernel 
\[
J_\Psi(x,y) =\sum_{n\geq 1}\frac{\mu_n}{1-\mu_n}\varphi_n (x)\varphi_n (y),
\qquad
x,y\in\Lambda, 
\]
by \eqref{eq:kdecomp}.
 For $\alpha=\{x_1,\ldots,x_n\}\in\espaceconfiglambda$,
we denote by $\mathrm{det}\,J_\Psi (\alpha)$
the determinant
$$
(x_1,\ldots,x_n)\mapsto
\mathrm{det}\,J_\Psi (x_1,\ldots,x_n) := 
\mathrm{det}\,\big(J_\Psi (x_i,x_j)_{1\leq i,j\leq n}\big),
$$ 
which is $\radon^{\otimes n}(x_1,\ldots,x_n)$-a.e. nonnegative; 
 see, e.g., the appendix of \cite{georgiiyoo}. 

\bigskip 

By Lemma 3.3 in \cite{shirai}
the determinantal point process $\Psi$
on $\Lambda $ with kernel
$K_\Psi(x,y)$, $x,y \in \Lambda$,
admits the Janossy densities
\begin{equation}
\label{eq:janossy}
j_\Psi^n(x_1,\ldots,x_n) = \mathrm{Det } (\mathrm{Id} -\mathcal{K}_\Psi  )\,
\mathrm{det}\,\big(J_\Psi (x_i,x_j)_{1\leq i,j\leq n}\big), 
\qquad x_1,\ldots,x_n \in \Lambda.  
\end{equation}
In addition, from e.g. \cite{shirai} (see Theorem 3.6 therein)
the Laplace transform \eqref{lapl} of $\Psi$ is given by
$$
\mathcal{L}_\Psi (f)= \Det \left( \mathrm{Id}-\mathcal{K}_\varphi  \right),
$$
for each nonnegative $f$ on $\underlying$ with compact support, where $\varphi = 1-e^{-f}$ and $\mathcal{K}_\varphi$ is the trace class integral operator with kernel
$$
K_\varphi (x,y) = \sqrt{\varphi(x)} K_\Psi (x,y)  \sqrt{\varphi(y)}, \quad x,y \in \underlying.
$$
\section{Determinantal PHD filter}

In this section we construct a second-order PHD filter
based on determinantal point processes.
We show that approximate closed-form filter
update expressions can be derived 
using approximation formulas stated in appendix 
for the corrector terms
$l^{(1)}_{z_{1:m}}$, $l^{(2)}_{z_{1:m}}$ 
and Janossy densities $j^{(n)}_\Phi$, when
the underlying point process has low cross-correlations.

\bigskip 

 In the sequel we will restrict the class of determinantal kernels
considered to a class of finite range interaction point
processes, by enforcing the condition
\begin{equation}
  \label{enf} 
J(x,y) = 0 \mbox{~for all~} x,y\in \Lambda  
\mbox{~such that~} |x-y|>\eta d(\Lambda),
\end{equation} 
as in e.g. Proposition~3.9 in \cite{georgiiyoo},
where $d(\Lambda)$ is the diameter of $\Lambda$ and $\eta \in (0,1)$. 

\subsubsection*{Prediction step} 
\label{s5}

The prediction point process $\Phi$
is constructed by branching
the prior point process $\Psi$ with a Bernoulli point process
$\Phi_s$ with
probability of survival $p_s(x)$ at the point $x \in \Lambda$, 
spatial likelihood density $l_s(\cdot |x)$ from state $x$, and 
characterized by the PGFl 
\begin{equation}
  \label{pgfl1}
        {\cal G}_{\Phi_s} ( g  \mid  x )
= 
1-p_s(x)+p_s(x)\int_\Lambda g (u) l_s(u|x) \nu (du).
\end{equation} 
According to \eqref{mspt}, the PGFl of
the prediction point process $\Phi$
is given by
\begin{align*}
  \mathcal{G}_\Phi (h)&=\mathcal{G}_{\Phi_b}(h)
  \mathcal{G}_\Psi (\mathcal{G}_{\Phi_s} (h  \mid  \cdot))
  \\
   & =
  \mathcal{G}_{\Phi_b}(h)
  \mathcal{G}_\Psi\left(
  1-p_s(\cdot)+p_s(\cdot)\int_\Lambda h (u) l_s(u\ \! |\ \! \cdot) \nu (du) \right), 
\end{align*}
 where $\mathcal{G}_{\Phi_b}$ is the PGFl of the Poisson birth
 point process $\Phi_b$ of new targets. 
 In the sequel we use the notation convention \eqref{qd}, i.e. 
 $\tilde{l}_s (x|u) : = p_s(u) l_s (x|u)$, for compactness of notation. 
\begin{prop}
 \label{pfs}
  Assume that the prior point process $\Psi$ is
  a determinantal point process with kernel
  $K_\Psi(x,y)$.
  Then, the prediction first and second-order moment densities of
  $\Phi$ are given by
  \begin{equation}
    \label{prediction_fom}
\mu^{(1)}_\Phi(x) =
\mu^{(1)}_{\Phi_b}(x)
+
\int_\Lambda \tilde{l}_s (x|u)K_\Psi(u,u)\nu (du),
\quad x\in \Lambda, 
\end{equation}  
 and 
 \begin{align} 
   \label{prediction_som}
      \rho^{(2)}_{\Phi}(x,y)
      = 
      &
      \int_{\Lambda^2}\tilde{l}_s (x|u)\tilde{l}_s (y|v)
      (K_\Psi(u,u)K_\Psi(v,v)-(K_\Psi(u,v))^2)\nu(du)\nu(dv)
  \\
  \nonumber
  &
  {
    +\mu^{(1)}_{\Phi_b}(y)\int_{\Lambda}\tilde{l}_s (x|u)K_{\Psi}(u,u)\nu(du)
    +\mu^{(1)}_{\Phi_b}(x)\int_{\Lambda}\tilde{l}_s (y|v)K_{\Psi}(v,v)\nu(dv)
+\rho^{(2)}_{\Phi_b}(x,y),
    }
\end{align} 
$x,y \in \Lambda$, $x\not= y$. 
\end{prop}
\begin{Proof}
  The expressions \eqref{prediction_fom}-\eqref{prediction_som} of 
  the prediction first and second-order (factorial) moment densities
  are obtained from \eqref{1stfm} and \eqref{2ndfm}
   as   
$$ 
  \mu^{(1)}_\Phi(x) = 
  \frac{\partial_{\delta_{x}} }{\partial h}
  \mathcal{G}_\Phi(h)_{\mid h=1}
= 
  \mu^{(1)}_{\Phi_b}(x) + \int_\Lambda p_s(u)l_s(x|u)\mu^{(1)}_\Psi(u)\nu (du)
$$ 
 and
\begin{align*}
  {
    \rho^{(2)}_{\Phi}(x,y)}
  &
  {
    =\int_{\Lambda^2}p_s(u)l_s(x|u)p_s(v)l_s(y|v)\rho^{(2)}_{\Psi}(u,v)\nu(du)\nu(dv)
  }
  {
    +\mu^{(1)}_{\Phi_b}(y)\int_{\Lambda}p_s(u)l_s(x|u)\mu^{(1)}_{\Psi}(u)\nu(du)
  }
  \\
& \quad 
  {
    +\mu^{(1)}_{\Phi_b}(x)\int_{\Lambda}p_s(v)l_s(y|v)\mu^{(1)}_{\Psi}(v)\nu(dv)
+\rho^{(2)}_{\Phi_b}(x,y), 
  }
  \qquad x,y\in \Lambda, \quad x\not= y.
\end{align*}
\end{Proof}
From Proposition~\ref{pfs}
we can model the prediction point process $\Phi$ as
a determinantal process with prediction kernel $K_{\Phi}$, whose
diagonal entries are given by 
$$ 
K_{\Phi}(x,x) = \mu^{(1)}_\Phi(x) =
\mu^{(1)}_{\Phi_b}(x)
+
\int_\Lambda \tilde{l}_s (x|u)K_\Psi(u,u)\nu (du), 
$$ 
and whose nondiagonal entries satisfy 
$$ 
K_\Phi(x,y)=\sqrt{K_\Phi(x,x)K_\Phi(y,y)-\rho^{(2)}_{\Phi}(x,y)},
\quad
x,y \in \Lambda, 
$$
from \eqref{kx2}, where $\rho^{(2)}_{\Phi}(x,y)$
is given by \eqref{prediction_som}
when $x\not= y$, and $\rho^{(2)}_{\Phi}(x,x):=0$, $x\in \Lambda$. 
 The prediction Janossy kernel $J_\Phi (x,y)$
 of the operator ${\mathcal J}_\Phi$ is then
 computed by the formula
$$
{\mathcal J}_\Phi = ( \mathrm{Id} - {\mathcal K}_\Phi )^{-1} {\mathcal K}_\Phi, 
$$
 see \eqref{idk}. 
\subsubsection*{Update step} 
\begin{enumerate}[1.] 
  \item First order moment update. 
Proposition~\ref{p1} and \eqref{kx} show that 
the diagonal values 
of the posterior kernel $K_{\Phi \mid \Xi =z_{1:m}}$ are given by 
\begin{eqnarray*} 
    K_{\Phi \mid \Xi =z_{1:m}} (x,x) 
  & = & \mu^{(1)}_{\Phi \mid \Xi =z_{1:m}} (x)
  \\
  & = &   q_d
l^{(1)}_{z_{1:m}} ( x)
     + 
    \sum_{z\in z_{1:m}} \tilde{l}_d(z|x)
    l^{(1)}_{z_{1:m}} (x;z),
    \quad x \in \Lambda,
    \quad
m\geq 0. 
\end{eqnarray*} 
Based on the approximation 
of the corrector term $l^{(1)}_{z_{1:m}}$ stated in Proposition~\ref{l1}, 
we will estimate the first-order posterior moment density as 
\begin{eqnarray}
 \nonumber 
  K_{\Phi \mid \Xi =z_{1:m}} (x,x)
  & = & \mu^{(1)}_{\Phi \mid \Xi =z_{1:m}} (x) 
  \\
  \label{update_fom}
  & \simeq &  
q_d K_\Phi ( x , x )
 + \sum_{z\in z_{1:m}}
    \frac{ J_\Phi ( x , x ) \tilde{l}_d (z | x)}{
     l_c (z ) + \int_{\Lambda} \tilde{l}_d (z | u ) J_\Phi ( u,u) \nu ( du )  },  
\end{eqnarray} 
 where we choose to approximate
$l^{(1)}_{z_{1:m}} ( x ) \simeq
\mu^{(1)}_\Phi (x)
= K_\Phi ( x , x )$
for consistency with the standard Poisson PHD filter. 
The estimate \eqref{update_fom} allows one to locate the targets by maximizing 
$K_{\Phi \mid \Xi =z_{1:m}} (x,x) = \mu^{(1)}_{\Phi \mid \Xi =z_{1:m}} (x)$
over $x$,
and to estimate the number of targets as 
 \begin{equation}\label{mu}
   \gamma_{\Phi|\Xi=z_{1:m}} := \int_\Lambda K_{\Phi|\Xi=z_{1:m}}(x,x) \nu (dx), 
   \quad x\in \Lambda. 
 \end{equation}  
\item Cross-diagonal kernel update. 
 As a consequence of
 \eqref{kxy}, i.e. 
 $$
 { 
  ( K_{\Phi \mid \Xi =z_{1:m}} (x,y) )^2 
  = \mu^{(1)}_{\Phi | \Xi =z_{1:m}} (x) 
  \mu^{(1)}_{\Phi | \Xi =z_{1:m}} (y) 
-
  \rho^{(2)}_{\Phi \mid \Xi =z_{1:m}} (x,y), 
 }
 $$ 
$x,y \in \Lambda$,  
 the cross-diagonal entries of the posterior kernel
 $K_{\Phi \mid \Xi =z_{1:m}} (x,y)$ can be estimated from 
 \eqref{kx2} as 
 \begin{equation} 
   \label{kx3} 
     K_{\Phi \mid \Xi =z_{1:m}} (x,y)
 = \sqrt{K_{\Phi \mid \Xi =z_{1:m}} (x,x)K_{\Phi \mid \Xi =z_{1:m}} (y,y)-\rho^{(2)}_{\Phi \mid \Xi =z_{1:m}} (x,y)},
\end{equation} 
 $x,y \in \Lambda$. 
 The above relation \eqref{kx3} can be rewritten
  from Propositions~\ref{p1}-\ref{p2} as 
    \begin{eqnarray} 
         \label{cjkdf}
        ( K_{\Phi \mid \Xi =z_{1:m}} (x,y) )^2 
 & = & q_d^2 \left(
      l^{(1)}_{z_{1:m}} ( x)
      l^{(1)}_{z_{1:m}} ( y)
      -
            l^{(2)}_{z_{1:m}} ( x,y)
\right)
      \\
      \nonumber 
       & & 
      {
        - q_d \sum_{z\in z_{1:m}} 
 \big( \tilde{l}_d (z | x) + \tilde{l}_d (z | y) \big) 
            l^{(2)}_{z_{1:m}} (x,y;z)
      }
      \\
\nonumber
& &
{
  + q_d
\sum_{z\in z_{1:m}}
\big(
\tilde{l}_d(z|y)
            l^{(1)}_{z_{1:m}} ( x)
            l^{(1)}_{z_{1:m}} (y;z)
 +
 \tilde{l}_d(z|x)
            l^{(1)}_{z_{1:m}} ( y)
            l^{(1)}_{z_{1:m}} (x;z)
 \big)
}
\\
\nonumber
& &
{
  + 
    \sum_{z\in z_{1:m}} \tilde{l}_d(z|x) \tilde{l}_d(z|y)
            l^{(1)}_{z_{1:m}} (x;z)
            l^{(1)}_{z_{1:m}} (y;z)
}
\\
\nonumber 
      & &
      {
        + 
        \sum_{z,z'\in z_{1:m} \atop z\not= z' }
        \! \! \! \!
        \tilde{l}_d(z|x) \tilde{l}_d(z'|y)
    \left(
            l^{(1)}_{z_{1:m}} (x;z)
            l^{(1)}_{z_{1:m}} (y;z')
    - 
            l^{(2)}_{z_{1:m}} ( x,y,z,z')
    \right)
    , \qquad 
      }
    \end{eqnarray} 
    $x,y \in \Lambda$, $m\geq 0$. 
 In practice we will estimate the posterior kernel
  $K_{\Phi \mid \Xi =z_{1:m}} (x,y)$ in \eqref{kx3} using 
  \eqref{update_fom} and 
  the approximation
  of the corrector term $l^{(2)}_{z_{1:m}}$
 in Proposition~\ref{l2}, to obtain 
 \begin{eqnarray}
    \label{update_som}
      \rho^{(2)}_{\Phi \mid \Xi =z_{1:m}} (x,y)  
  & \simeq &  
    q_d^2 (J_\Phi (x,x)J_\Phi (y,y)-J_\Phi (x,y)^2 )
  \\
    \nonumber
    &   &
    {
      +     q_d  \sum_{z\in z_{1:m}} 
    \frac{
          (J_\Phi (x,x)J_\Phi (y,y)-J_\Phi (x,y)^2 )
 \big( \tilde{l}_d (z | x)
      +
 \tilde{l}_d (z | y) \big)}{
   s_c(z)
    }
    }
        \\
    \nonumber 
  &   &
    {
      + 
     \sum_{z,z'\in z_{1:m} \atop z \not= z'}
     \frac{  
       (J_\Phi (x,x)J_\Phi (y,y)-J_\Phi (x,y)^2 )
 \tilde{l}_d (z | x)
 \tilde{l}_d (z' | y)
}{
         s_c(z) s_c(z')
         -
  \int_{\Lambda^2} J_\Phi (u,v)^2 \tilde{l}_d (z | u) \tilde{l}_d (z' | v) \nu (du) \nu (dv) 
 }
    }
\end{eqnarray} 
with $\rho^{(2)}_{\Phi \mid \Xi =z_{1:m}} (x,x):=0$, $x\in \Lambda$,
see \eqref{djkls2} which yields the expression of
$K_{\Phi \mid \Xi =z_{1:m}} (x,y)^2$ in Proposition~\ref{thasd}, where
$s_c(z)$ is defined in \eqref{scz}.  
 After completing the update step, we move to the next prediction step
 by taking $K_\Psi (x,y) := K_{\Phi \mid \Xi =z_{1:m}} (x,y)$,
 $x,y \in \Lambda$. 
\end{enumerate} 
\section{Implementation} 
\label{s7}
We implement the Determinantal Point Process (DPP)
and Poisson Point Process (PPP) PHD filters using
the sequential Monte Carlo (or particle filtering) method 
as in \cite{litiancheng},
together with the roughening method of \cite{sattar},
which allow us to estimate otherwise intractable integrals
using discretized particle summations. 
Our ground truth dynamics
follows the \textit{nearly constant turn-rate} motion dynamics
of \cite{cantoni}, \cite{litiancheng},
with the addition of a repulsion term. 
The state of each target at time $t$ is given by
$\mathbf{x}_t=\left(x_t, \dot{x}_t, y_t, \dot{y}_t, \theta_t\right)^\top$, where $x_t$, $y_t$ are the cartesian coordinates, $\dot{x}_t$, $\dot{y}_t$ are the respective velocities, and $\theta_t$ is the turn rate.
At time $t+1$, the location of every target $i$ for $i\in \{1,\ldots , n\}$
is given by 
\begin{equation} 
  \label{motion_dynamics}
  \mathbf{x}^i_{t+1}=\mathbf{F}(\theta_t)\mathbf{x}^i_t+\mathbf{G}\mathbf{v}^i_t
  +
  \mathbf{s}^i_t,
\end{equation}
where
$$
\mathbf{s}^i_t=\left(\begin{array}{c}
\zeta_x\sum_{j=1}^N\frac{x^i_t-x^j_t}{\left\vert \mathbf{x}^i_t-\mathbf{x}^j_t\right\vert}\\
0\\
\zeta_y\sum_{j=1}^N\frac{y^i_t-y^j_t}{\left\vert \mathbf{x}^i_t-\mathbf{x}^j_t\right\vert}\\
0\\
0 
\end{array}\right)
$$
 is a term which models repulsion among targets. 

\bigskip 

 Here, $\mathbf{v}_t=(v_x,v_y,v_{\theta})^\top$
is a zero-mean acceleration noise
distributed according to the zero-mean Gaussian noise
\begin{align}
  \label{vt} 
\mathbf{v}_t \sim \mathcal{N}\left(\left(\begin{array}{c}
0 \\
0 \\
0
\end{array}\right),\left(\begin{array}{ccc}
\sigma^2_{v_x} & 0 & 0\\
0 & \sigma^2_{v_y} & 0 \\
0 & 0 & \sigma^2_{v_{\theta}}
\end{array}\right)\right), 
\end{align}
 and 
$$\mathbf{F}=\left(\begin{array}{ccccc}
1 & {\sin(\tau\theta_t)}/{\theta_t} & 0 & {(\cos(\tau\theta_t)-1)}/{\theta_t}& 0\\
0 & \cos(\tau\theta_t) & 0 & -\sin(\tau\theta_t)& 0\\
0 & {(\cos(\tau\theta_t)-1)}/{\theta_t}& 1 & {\sin(\tau\theta_t)}/{\theta_t} & 0\\
0 & \sin(\tau\theta_t) & 0 & \cos(\tau\theta_t) & 0\\
0 & 0 & 0 & 0 & 1
\end{array}\right),
 $$
 $$ \mathbf{G}=\left(\begin{array}{ccc}
\tau^2/2 & 0 & 0\\
\tau & 0 & 0\\
0 & \tau^2/2 & 0\\
0 & \tau & 0\\
0 & 0 & \tau
\end{array}\right),
$$
with $\tau > 0$ the time sampling period.
 When $\mathbf{F}=0$ and $\mathbf{G}=I_d$,
the repulsive interaction motion dynamics
\eqref{motion_dynamics} 
has the law the Ginibre DPP for stationary
 distribution; see, for example, Equation~(2.19) in \S~2.2 of \cite{osada2}.

\bigskip 

 The measurement vector of each target at time $t$ is
written as $\mathbf{m}_t =(m_{r_t},m_{\omega_t})^\top$ 
using the range and bearing components 
$m_{r_t}$ and $m_{\omega_t}$. 
The measurement generated by every target at time $t+1$ is then given by 
\begin{align}\label{measurement_eqn}
\mathbf{m}_{t+1}=\mathbf{p}_{t+1}+\mathbf{w}_{t+1},
\end{align}
where
 $\mathbf{p}_{t+1}=\left(\begin{array}{c}
\sqrt{x^2_{t+1}+y^2_{t+1}}\\
\arctan( y_{t+1} / x_{t+1} )
\end{array}\right)$ and
 the measurement noise vector
$\mathbf{w}_{t+1}=(w_{r_{t+1}},w_{\omega_{t+1}})^\top$ 
is distributed according to the 
zero-mean Gaussian noise 
\begin{align}
  \label{wt} 
\mathbf{w}_{t+1} \sim \mathcal{N}\left(\left(\begin{array}{c}
0 \\
0
\end{array}\right),\left(\begin{array}{cc}
\sigma^2_{w_r} & 0\\
0 & \sigma^2_{w_{\omega}}
\end{array}\right)\right).
\end{align}
The spatial likelihood densities
$l_s( \mathbf{z} |\mathbf{x})$ and $\tilde{l}_d(\mathbf{z}|\mathbf{x})$ 
from a target state $\mathbf{x}$ 
to a measurement $\mathbf{z}$
in \eqref{abdfasfd} and \eqref{pgfl1}
 respectively follow the zero-mean 
multivariate Gaussian distribution \eqref{vt} of $\mathbf{v}_t$
and the multivariate Gaussian distribution \eqref{wt} of $\mathbf{w}_t$.   
In addition, the model generates measurement information from every
target with a constant probability of detection $p_d$,
and the spatially distributed clutter measurement points are generated
according to a Poisson point process 
with constant density $l_c(z)$. 

\bigskip 

\begin{figure}[H]
\centering
\hskip-0.2cm
\begin{subfigure}{.48\textwidth}
\centering
\includegraphics[width=1.\textwidth]{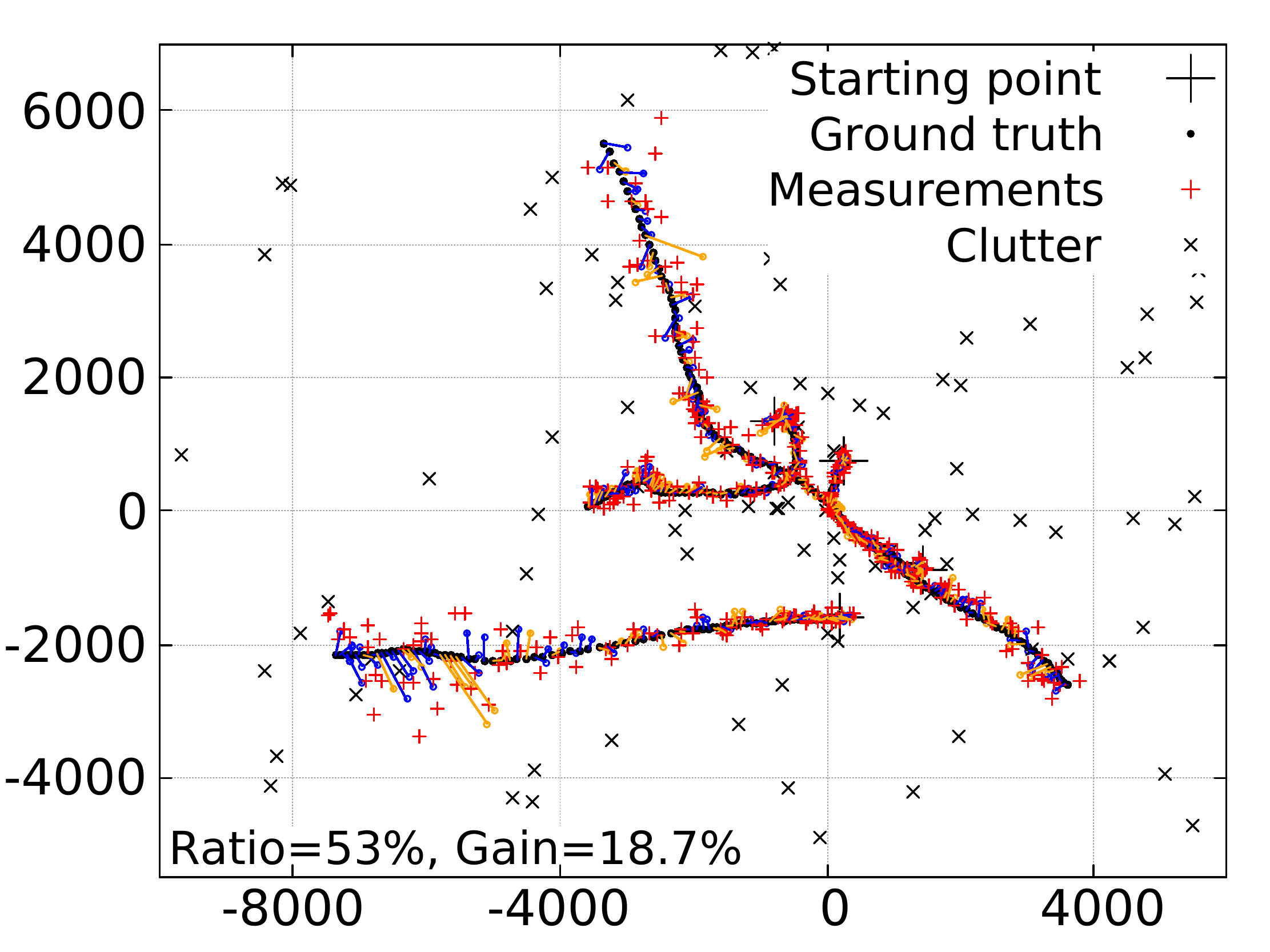}
\caption{\small PHD filter with no repulsion} 
\label{fig1.0}
\end{subfigure}
\hskip0.5cm
\begin{subfigure}{.48\textwidth}
\centering
\includegraphics[width=1.\textwidth]{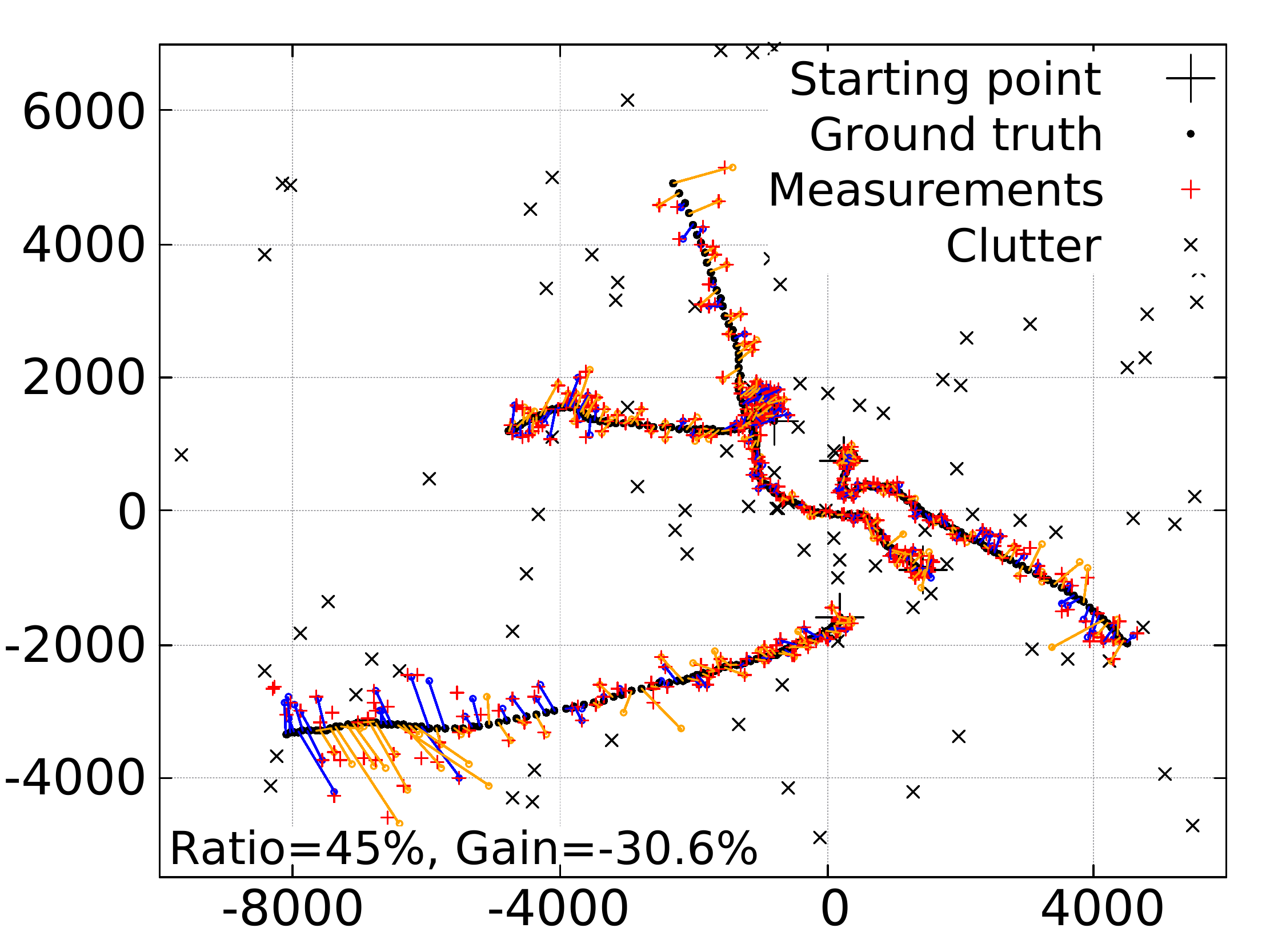}
\caption{\small PHD filter with repulsion} 
\label{fig2.0}
\end{subfigure}
\caption{\small PHD filter simulations with and without repulsion} 
\label{fig2.00}
\end{figure}

\begin{figure}[H]
\centering
\hskip-0.2cm
\begin{subfigure}{.48\textwidth}
\centering
\includegraphics[width=1.\textwidth]{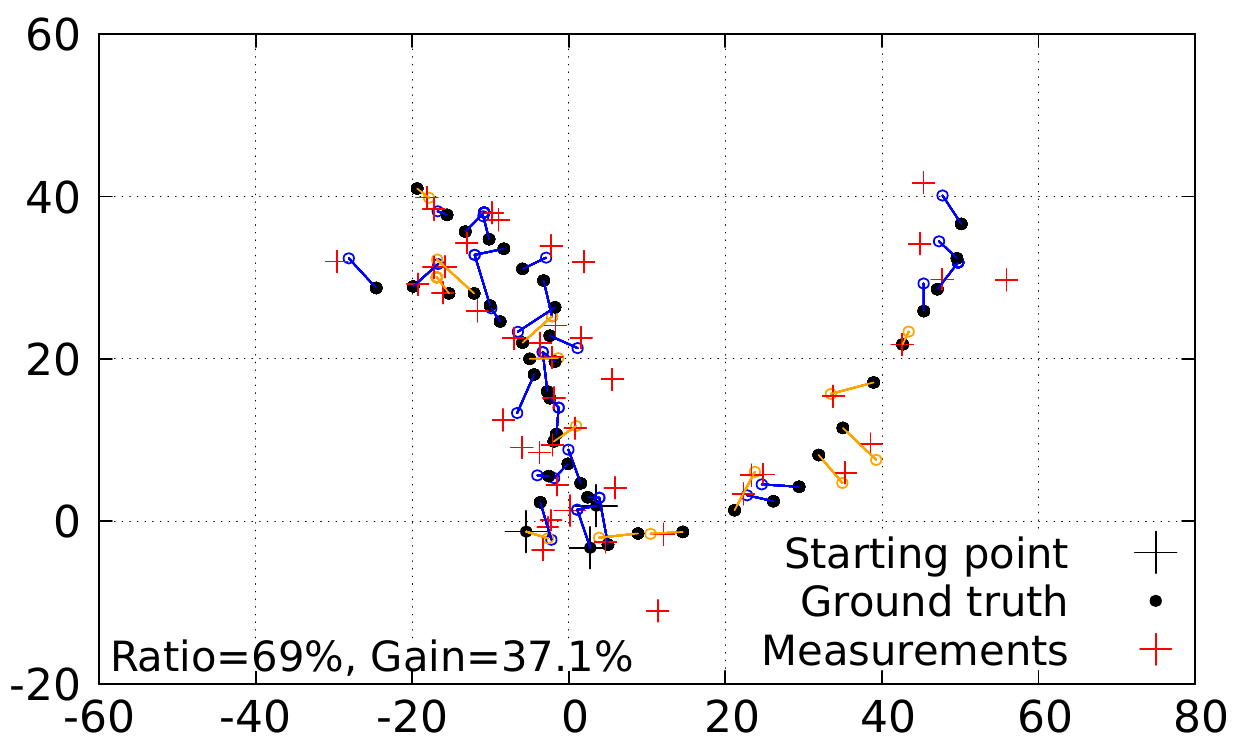}
\caption{\small PHD filter with data association for $\zeta = 0$.}
\label{fig6}
\end{subfigure}
\begin{subfigure}{.48\textwidth}
\centering
\includegraphics[width=1.\textwidth]{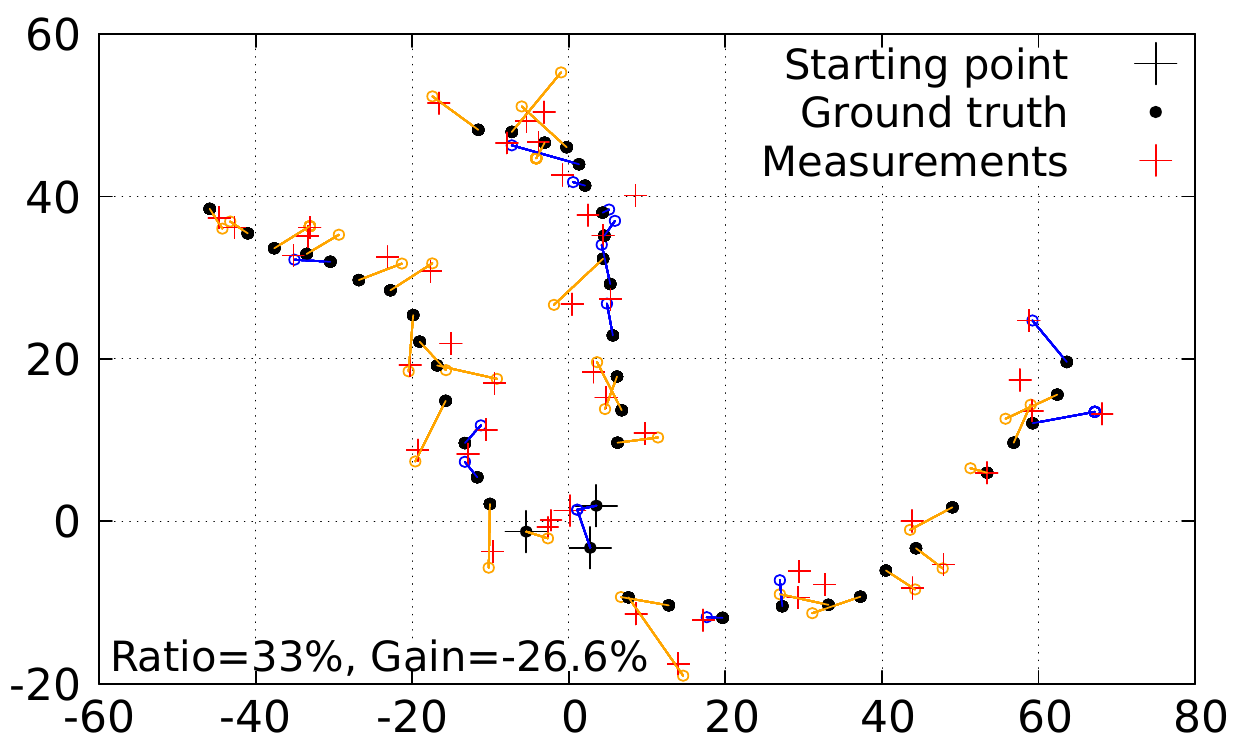}
\caption{\small PHD filter with data association for $\zeta =30$.}
\label{fig14}
\end{subfigure}
\caption{\small PHD filter with data association for $\zeta =0$ and $\zeta =30$.}
\label{fig14.00}
\end{figure}

\begin{figure}[H]
\centering
\begin{subfigure}{.5\textwidth}
\centering
\includegraphics[width=1.\textwidth]{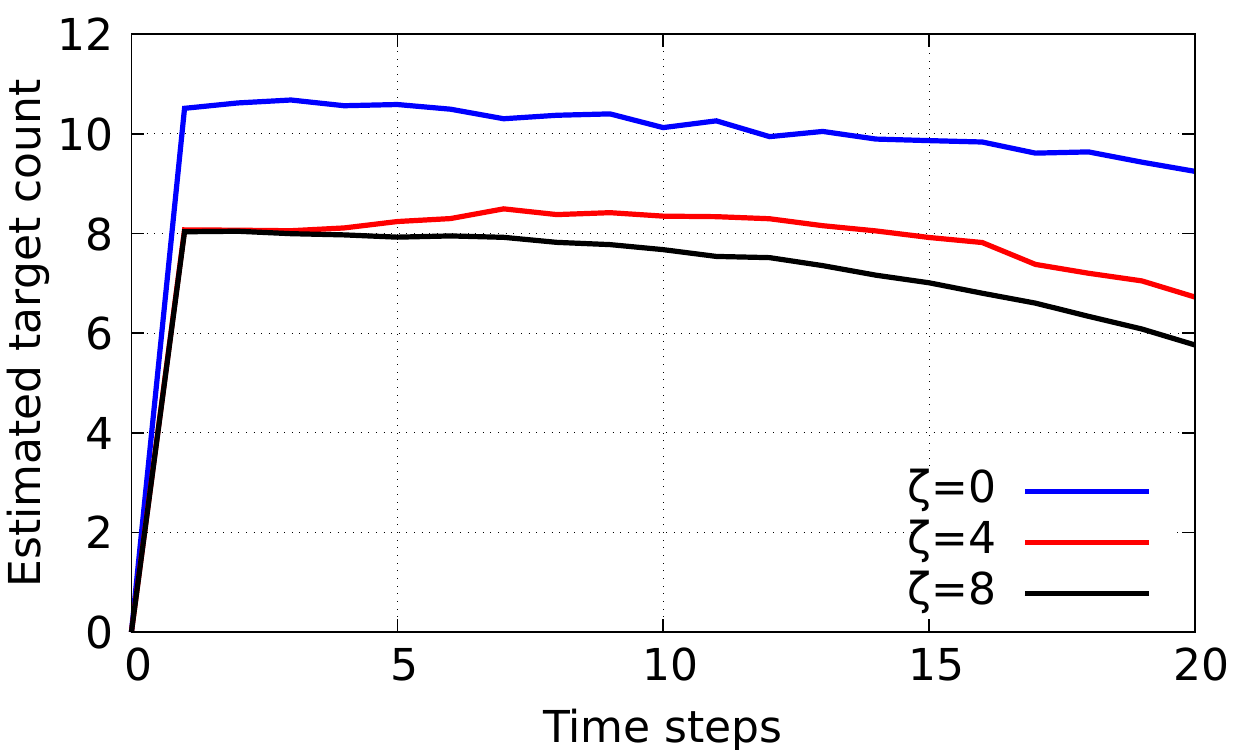}
\caption{\small First moment estimation of target count.}
\label{fig3.0}
\end{subfigure}
\hskip-0.15cm
\begin{subfigure}{.5\textwidth}
\centering
\includegraphics[width=1.\textwidth]{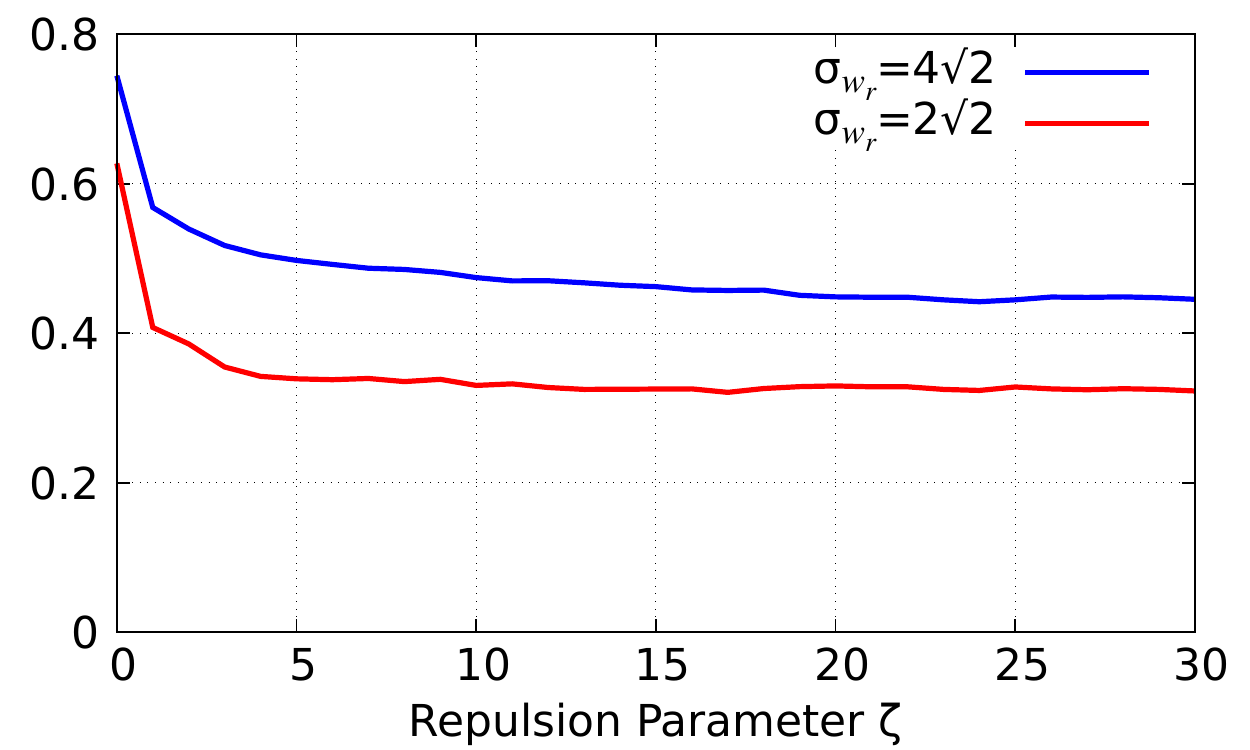}
\caption{\small Good estimate ratio with $p_d=1$ and no clutter.}
\label{fig4}
\end{subfigure}
\caption{\small Graphs of first moment target counts and good estimate ratios.}
\label{fig4.00}
\end{figure}

The implementation of Figures~\ref{fig2.00} to \ref{fig4.00}
use $1000$ particles at initialization,
 $100$ resampling particles per expected target,
 and $100$ new particles per expected target birth,
 which follows a time-dependent Poisson birth process.
 The starting locations are uniformly distributed within a
  square domain, we take the spatial standard deviations (s.d.) $\sigma_{v_x}=\sigma_{v_y}= 1 \ \mathrm{m}/\mathrm{s}^2$, the turn-rate noise s.d. $\sigma_{v_{\theta}}=\pi = 180 \ \mathrm{rad}/\mathrm{s}$, bearing distribution s.d. $\sigma_{w_{\omega}}= \pi = 180 \ \mathrm{rad}$, and range distribution s.d. $\sigma_{w_r}=2\sqrt{2}\ \mathrm{m}$. 

 \bigskip

 For illustration and performance assessment purposes, our code
 displays the association  
 between target-originated measurements and posterior state estimates.
 For this, given a measurement we select its associated estimate
 by minimizing the distance between the measurement and all candidate
 estimates.
 In addition, the blue edges show the estimates which
 improve over the corresponding measurements
 in terms of Euclidean distances
 to the ground truth,
 while the orange edges show the estimates which perform worse than measurements. 
 
 \bigskip

 Our simulations also display a ratio of good estimate counts
 against total measurement counts,
 as well as a gain metric which
 measures the relative improvement in distance between
 estimates and measurements.
 Positive gain correspond to a good estimate ratio above $50\%$, and
 negative gain is realized when the ratio falls below $50\%$.
 
\vspace{-0.2cm}

Figure~\ref{fig2.00}-(\subref{fig2.0}) 
shows that when using a nonzero value for
the repulsion parameter $\zeta=\zeta_x=\zeta_y$, the good estimate ratio
with repulsive interaction becomes 
lower as compared with the non repulsive setting
of Figure~\ref{fig2.00}-(\subref{fig1.0}), with $50$ time steps. 
In Figure~\ref{fig2.00} 
 the Poisson clutter rate is $l_c=1$, the probability of detection $p_d=0.9$ and the probability of survival is $p_s=1$,
with four targets. 

\bigskip

In Figures~\ref{fig14.00}-(\subref{fig6})-(\subref{fig14}) we provide
further illustrations of three-target interaction and PHD filter output of a single trial at different repulsion values $\zeta=0$ and $30$, with $15$ time steps
and $p_d=p_s=1$. 

Figure~\ref{fig14.00}-(\subref{fig6}) shows the PHD filter output
with $\zeta =0$, where the targets can become closer to each other
without repulsion, and with positive gain. 
For the repulsion value $\zeta =30$ as in Figure~\ref{fig14.00}-(\subref{fig14}) 
the repulsion effect among the three targets become much more evident,
the good estimate ratio becomes lower, 
and the gain becomes negative. 

\bigskip

The above results are summarized in Figure~\ref{fig4.00}. 
 Figure~\ref{fig4.00}-(\subref{fig3.0}) 
 presents the SMC-PHD filter first-order moment
 output with $200$ Monte Carlo runs and
 $10$ targets across $20$ time steps, with different repulsion
 parameter values $\zeta =0,4,8$. 
 Figure~\ref{fig4.00}-(\subref{fig4})
 presents the good estimate ratio for
various values of the repulsion coefficient $\zeta$, 
with $100$ Monte Carlo runs on three targets across $15$ time steps,
with $p_d=1$ and no clutter. 

We note that the PHD filter is correctly estimating 
the target count when 
the repulsion coefficient $\zeta$ vanishes;  
however, the estimation falls short
for nonzero values of $\zeta$,
showing the performance degradation of the
PPP-PHD filter in the presence of target interaction. 

\section{Determinantal PHD filtering algorithm} 
\label{s7.1}
\subsubsection*{Initialization $(t=0)$} 
The state dynamics of the initial set of $N_{\Phi,0}$ particles
is sampled according to a uniform distribution on
the state space $\Lambda$.
The diagonal entries
of the prior discretized determinantal kernel $K_{\Phi,0}$
at time $t=0$ are initialized to $\gamma_{\Phi,0} / N_{\Phi,0}$
with $\gamma_{\Phi,0}$ a prior intensity value. 
The nondiagonal entries
are initialized to $\alpha\gamma_{\Phi,0} / N_{\Phi,0}$ where $\alpha \geq 0$, 
except for those which are set to zero according to Condition~\eqref{enf}
with the matrix index threshold $\eta =10\%$.

\bigskip

Using \eqref{idk} we then compute the discretized Janossy kernel $J_{\Phi,0}$ 
which is needed for the evaluation of the posterior determinantal kernel
$\tilde{K}_{\Phi|\Xi,0}$. 
Letting
$P_p$ denote the number of resampled particles per target
and 
$\gamma_0 :=\sum_{i=1}^{N_{\Phi,0}}K_{\Phi|\Xi,0}(\mathbf{x}_i,\mathbf{x}_i)$,
we resample $\tilde{N}_0 :=P_p\times\lfloor\gamma_0\rfloor$
particles $\{\mathbf{x}_i\}_{i=1}^{\tilde{N}_0}$
that better describe the target locations 
 as in \cite{sattar} 
 by maximizing the diagonal entries of 
 $\tilde{K}_{\Phi|\Xi,0}$ over $\Lambda$,
 where $\lfloor \cdot\rfloor$ denotes the
 integer floor function. 
Those particles are then used to initialize the
post-resampling determinantal kernel 
$K_0$ and to compute the post-resampling Janossy kernel $J_0$
in order 
to estimate the updated kernel $\tilde{K}_{\Phi|\Xi}$
 from \eqref{update_fom}, \eqref{kx3} and \eqref{update_som}.
The updated kernel $\tilde{K}_{\Phi|\Xi}$ is then set
as the prior kernel $K_\Psi$ of the next time step.
 
\bigskip

\floatname{algorithm}{Initialization (time $t=0$)}
\renewcommand{\thealgorithm}{}

\begin{algorithm}[H]
  \footnotesize 
  \caption{} 
  \vspace{-0.4cm}
\begin{multicols}{2}
  \begin{algorithmic}
    \STATE {Set} $\gamma_{\Phi,0}\in\mathbb{R}^+$, $P_p\in\mathbb{N}$,
    $\eta \in (0,1)$, 
    and $\alpha\in\mathbb{R}^+$.
	\STATE {Sample} $N_{\Phi,0}$ particle state dynamics $\{\mathbf{x}_i\}_{i=1}^{N_{\Phi,0}}$ for initial birth process $\Phi$ uniformly distributed within state space $\Lambda$.
	\STATE {Initialize} the (prior) determinantal kernel $K_{\Phi,0}$.
	\FOR{$1\leq i, j \leq N_{\Phi,0}$}
	\STATE $K_{\Phi,0}(\mathbf{x}_i,\mathbf{x}_i) := $ \\ $ \gamma_{\Phi,0} / N_{\Phi,0}$
	\IF{$i\neq j$ \AND $|i-j|\leq \eta P_b$}
	\STATE $K_{\Phi,0}(\mathbf{x}_i,\mathbf{x}_j) := $ \\ $ \alpha\gamma_{\Phi,0} / N_{\Phi,0}$
	\ENDIF
	\ENDFOR
	\STATE {Compute} the Janossy kernel $J_{\Phi,0} :=(I-K_{\Phi,0})^{-1}K_{\Phi,0}$.
	\STATE {Compute} the posterior determinantal kernel $K_{\Phi|\Xi,0}$ using \eqref{update_fom}, \eqref{kx3} and \eqref{update_som}.
	\STATE {Perform} resampling as in \cite{sattar} to obtain the particle state dynamics $\{\mathbf{x}_i\}_{i=1}^{\tilde{N}_0}$ where $\tilde{N}_0 :=P_p\times\lfloor\gamma_0\rfloor$ and $\gamma_0 :=\sum_{i=1}^{N_{\Phi,0}}K_{\Phi|\Xi,0}(\mathbf{x}_i,\mathbf{x}_i)$.
    \STATE {Initialize} the post-resampling determinantal kernel $K_0$ as follows: 
	\FOR{$1\leq i, j \leq \tilde{N}_0$}
	\STATE $K_0(\mathbf{x}_i,\mathbf{x}_i) := \gamma_0 / \tilde{N}_0$
	\IF{$i\neq j$ \AND $|i-j|\leq \eta P_p$}
	\STATE $K_0(\mathbf{x}_i,\mathbf{x}_j) := \alpha\gamma_0 / \tilde{N}_0$
	\ENDIF
	\ENDFOR
	\STATE {Compute} the post-resampling Janossy kernel \\ $J_0=(I-K_0)^{-1}K_0$.
	\STATE {Compute} the posterior determinantal kernel $\tilde{K}_{\Phi|\Xi,0}$ using \eqref{update_fom}, \eqref{kx3} and \eqref{update_som}.
	
	\STATE {Estimate} the number of targets as
        $\gamma_{\Phi|\Xi,0} := \sum_{i=1}^{\tilde{N}_{0}} \tilde{K}_{\Phi|\Xi,0}(\mathbf{x}_i,\mathbf{x}_i)$.
	\end{algorithmic}
\end{multicols}
\vspace{-0.2cm}
\end{algorithm}

\subsubsection*{Algorithm $(t\geq 1)$} 

The general algorithm proceeds to compute the prediction state transition dynamics $\{\mathbf{x}_{t+1|t}^{(i)}\}_{i=1}^{\tilde{N}_t}$ using \eqref{motion_dynamics}, 
followed by the computation of
the prediction determinantal state transition kernel $K_{\Phi,t+1|t}$ using \eqref{prediction_fom} and \eqref{prediction_som}.
Letting $P_b$ denote the number of particles per birth target
and 
$\gamma_{\Phi,t+1}:=\sum_{i=1}^{\tilde{N}_t}K_{\Psi,t+1}(\mathbf{x}_i,\mathbf{x}_i) \nu (\{x_i \})$,
 we sample the state dynamics of
$N_{\Phi,t+1}:=P_b\times\lfloor\gamma_{\Phi,t+1}\rfloor$ particles for the
target birth process $\Phi$.  
The discretized prediction determinantal kernel $K_{\Phi,t+1|t}$ is then extended to incorporate
the set of additional $N_{\Phi,t+1}$ particles by assigning the diagonal entries
corresponding to these new particles to $\gamma_{\Phi,t+1} / N_{\Phi,t+1}$ and 
the nondiagonal entries to
$\alpha\gamma_{\Phi,t+1} / N_{\Phi,t+1}$,
and by setting all other new entries to $0$
according to Condition~\eqref{enf} with the matrix index threshold $\eta :=10\%$. 
Thereafter, we compute the discretized Janossy kernel $J_{\Phi,t+1|t}$ using \eqref{idk}
and then the discretized posterior determinantal kernel $K_{\Phi|\Xi,t+1}$ using \eqref{update_fom}, \eqref{kx3} and \eqref{update_som}.
Next, letting
$\gamma_{t+1}:=\sum_{i=1}^{N_{t+1}}K_{\Phi|\Xi,t+1}(\mathbf{x}_i,\mathbf{x}_i)$
we resample
$\tilde{N}_{t+1}:=P_p\times \lfloor\gamma_{t+1}\rfloor$
particles with state dynamics $\{\mathbf{x}_i\}_{i=1}^{\tilde{N}_{t+1}}$ 
as in \cite{sattar}, 
by maximizing the diagonal entries of $K_{\Phi|\Xi,t+1}$
over $\Lambda$. 
Those particles are then used to initialize the
post-resampling determinantal kernel $K_{t+1}$ 
by setting diagonal entries to $\gamma_{t+1} / \tilde{N}_{t+1}$ and  
nondiagonal entries to $\alpha\gamma_{t+1} / \tilde{N}_{t+1}$,
except for those which are set to zero according to Condition~\eqref{enf}
with $\eta =10\%$. Finally, 
we recompute the post-resampling Janossy kernel $J_{t+1}$ and
the posterior determinantal kernel $\tilde{K}_{\Phi|\Xi,t+1}$ using \eqref{update_fom}, \eqref{kx3} and \eqref{update_som}.  

\floatname{algorithm}{DPP-PHD Filter (time $t+1\geq 1$)}

\begin{algorithm}[h]
  \footnotesize
  \caption{} 
  \vspace{-0.4cm}
  \begin{multicols}{2}
    \begin{algorithmic}
	\STATE {Compute} the (prediction) state transition dynamics $\{\mathbf{x}_{t+1|t}^{(i)}\}_{i=1}^{\tilde{N}_t}$ from $\{\mathbf{x}_{t}^{(i)}\}_{i=1}^{\tilde{N}_t}$ based on \eqref{motion_dynamics}.
	\STATE {Compute} the (prediction) determinantal state transition kernel $K_{\Phi, t+1|t}$ from $K_{\Psi,t+1} :=\tilde{K}_{\Phi|\Xi,t}$ using \eqref{prediction_fom} and \eqref{prediction_som}.
	\STATE {Sample} $N_{\Phi,t+1}$ new particle state dynamics for birth process $\Phi$ at time $t+1$ uniformly distributed within state space $\Lambda$ to generate $\{\mathbf{x}^{(i)}\}_{i=1}^{N_{\Phi,t+1}}$ where $N_{\Phi,t+1}:=P_{b}\times \lfloor\gamma_{\Phi,t+1}\rfloor$ and $\gamma_{\Phi,t+1}:=\sum_{i=1}^{\tilde{N}_t}K_{\Psi,t+1}(\mathbf{x}_i,\mathbf{x}_i)$.
	\STATE {Extend} the (prediction) determinantal kernel $K_{\Phi,t+1|t}$ to dimension $N_{t+1}:=\tilde{N}_t+N_{\Phi,t+1}$ with state dynamics $\{\mathbf{x}_i\}_{i=1}^{N_{t+1}}:=\{\mathbf{x}_{t+1|t}^{(i)}\}_{i=1}^{\tilde{N}_t}\cup\{\mathbf{x}^{(i)}\}_{i=1}^{N_{\Phi,t+1}}$, where the following indexes are allocated to new particles.
	\FOR{$\tilde{N}_t+1\leq i,j \leq N_{t+1}$}
	\STATE $K_{\Phi,t+1|t}(\mathbf{x}_i,\mathbf{x}_i):= \gamma_{\Phi,t+1} / N_{\Phi,t+1}$
	\IF{$i\neq j$ \AND $|i-j|\leq \eta P_b$}
	\STATE $K_{\Phi,t+1|t}(\mathbf{x}_i,\mathbf{x}_j):= \alpha\gamma_{\Phi,t+1} / N_{\Phi,t+1}$
	\ENDIF
	\ENDFOR
	\FOR{$1\leq i,j \leq N_{t+1}$}
	\IF{$i\leq \tilde{N}_t \And j\geq \tilde{N}_t+1$}
	\STATE $K_{\Phi,t+1|t}(\mathbf{x}_i,\mathbf{x}_j):=0$
	\STATE $K_{\Phi,t+1|t}(\mathbf{x}_j,\mathbf{x}_i):=0$
	\ENDIF
	\ENDFOR
	\STATE {Compute} the Janossy kernel \\ $J_{\Phi,t+1|t}:=$ \\ $(I-K_{\Phi,t+1|t})^{-1}{K_{\Phi,t+1|t}}$.
	\STATE {Compute} the posterior determinantal kernel $K_{\Phi|\Xi,t+1}$ using \eqref{update_fom}, \eqref{kx3} and \eqref{update_som}.
	\STATE {Perform} resampling as in \cite{sattar} to obtain the particle state dynamics $\{\mathbf{x}_i\}_{i=1}^{\tilde{N}_{t+1}}$ where $\tilde{N}_{t+1}:=P_p\times\lfloor\gamma_{t+1}\rfloor$ (capped at $1000$ particles) and $\gamma_{t+1}:=\sum_{i=1}^{N_{t+1}}K_{\Phi|\Xi,t+1}(\mathbf{x}_i,\mathbf{x}_i)$. 
	\STATE {Initialize} the post-resampling determinantal kernel $K_{t+1}$ by 
	\FOR{$1\leq i, j \leq \tilde{N}_{t+1}$}
	\STATE $K_{t+1}(\mathbf{x}_i,\mathbf{x}_i):= \gamma_{t+1} / \tilde{N}_{t+1}$
	\IF{$i\neq j$ \AND $|i-j|\leq \eta P_p$}
	\STATE $K_{t+1}(\mathbf{x}_i,\mathbf{x}_j):= $ \\ $ \alpha\gamma_{t+1} / \tilde{N}_{t+1}$
	\ENDIF
	\ENDFOR
	\STATE {Compute} the post-resampling Janossy kernel \\ $J_{t+1}=(I-K_{t+1})^{-1}{K_{t+1}}.$
	\STATE {Compute} the posterior determinantal kernel $\tilde{K}_{\Phi|\Xi,t+1}$ using \eqref{update_fom}, \eqref{kx3} and \eqref{update_som}.
	\STATE {Estimate} the number of targets as \\  
        $\gamma_{\Phi|\Xi,t+1} := $ \\ $\sum_{i=1}^{\tilde{N}_{t+1}} \tilde{K}_{\Phi|\Xi,t+1}(\mathbf{x}_i,\mathbf{x}_i)$.
	\end{algorithmic}
    \end{multicols}
  \vspace{-0.2cm}
\end{algorithm}

The complexity of this DPP-PHD filter is cubic in the number of discretization steps
due to the presence of matrix inversions in the algorithm. 

\subsubsection*{Numerical results} 
\begin{figure}[H]
\begin{subfigure}{.4\textwidth}
\centering
\includegraphics[width=7.8cm,height=4.5cm]{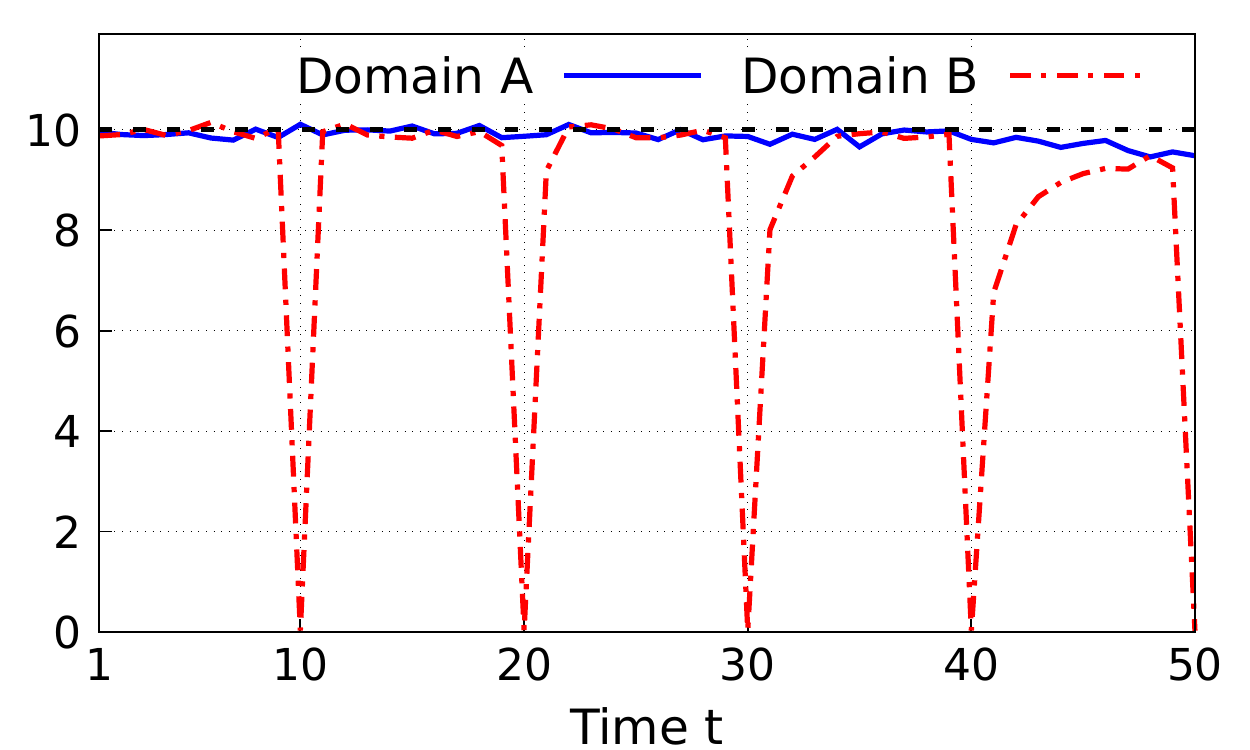}
\caption{Target counts estimates.} 
\label{fig1a}
\end{subfigure}
\hskip1.5cm
\begin{subfigure}{.4\textwidth}
\centering
\includegraphics[width=7.8cm,height=4.5cm]{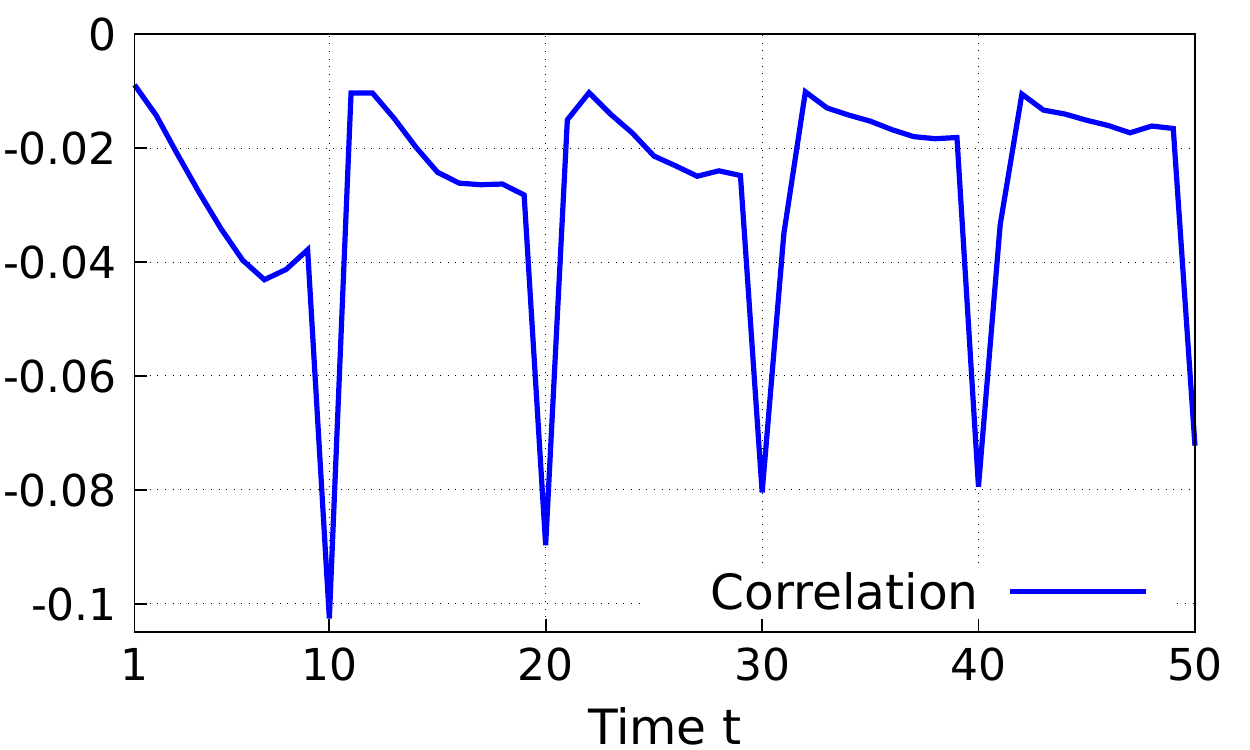}

\caption{Correlation estimates.} 
\label{fig1b}
\end{subfigure}
\caption{MC estimates with $p_d=0.9$ and $5$ clutter points per domain with $\alpha = 4.0$.} 
\label{fig1} 
\end{figure}

In Figure~\ref{fig1} we assess the spooky effect (see \cite{franken})
of our DPP-PHD filter, 
following the approach applied in \cite{sdhc} to second-order PHD filters. 
 Our tracking scenario consists of
two disjoint square domains $A$ and $B$ of size $150$~m by $150$~m,
which are located 
150m diagonally apart.
In each domain, $10$ targets are initialized and their state dynamics are centrally distributed at the first time step.
The runtime of the experiment is set at $50$
with $100$ Monte Carlo (MC) runs, 
the targets survive throughout ($p_s=1$),
and their trajectories remain within the observation domains. 
In Figure~\ref{fig1} we
take the spatial standard deviations (s.d.) $\sigma_{v_x} = \sigma_{v_y} = 1.0$~m/s$^2$, 
turn-rate noise s.d. $\sigma_{v_{\theta}} = \pi = 180$~rad/s,
bearing distribution s.d. $\sigma_{w_{\omega}} = \pi = 180$~rad,
and range distribution s.d. $\sigma_{w_r} = \sqrt{2}$~m. 

\bigskip

In a similar setting to \cite{sdhc},
all targets in domain $B$ are compelled to be misdetected in every cycle of $10$ time steps. We use a constant probability of detection $p_d=0.9$ and mean clutter count at $5$ in each measurement space.
 
\bigskip

At initialization in Figure~\ref{fig1} we set $N_{\Phi,0}=800$,
$\gamma_{\Phi,0}=2$
and
$\alpha = 4$.
The DPP-PHD
filter implementation uses $P_p=30$ resampled particles per target, and $P_b=10$ particles per birth target.

\bigskip

Figure~\ref{fig1}-(\subref{fig1a}) 
shows the estimated intensities in domains $A$ and $B$, where domain $A$ is unaffected by the rapid drop in the intensity of domain $B$.
 The posterior correlation estimates in Figure~\ref{fig1}-(\subref{fig1b}) 
 are computed by rescaling the covariance expression \eqref{fksdlfds}
 written as 
$$
 c^{(2)}_{\Phi \mid \Xi =z_{1:m}} (A,B)  
 = - \int_{A\times B} ( K_{\Phi \mid \Xi =z_{1:m}} (x,y) )^2 \nu (dx) \nu (dy), 
$$
 as in Corollary~\ref{fdsfsfd}, where $K_{\Phi \mid \Xi =z_{1:m}} (x,y)$ is estimated
 as in Proposition~\ref{thasd} from 
 \eqref{update_fom} and \eqref{kx3}.
 Figure~\ref{fig1}-(\subref{fig1b}) shows negative correlations
 due to the determinantal point process nature,
 which leads to a drop in negative correlation during
 the compelled misdetection at each $10$-steps cycle.

\begin{figure}[H]
\begin{subfigure}{.4\textwidth}
\centering
\includegraphics[width=7.8cm,height=4.5cm]{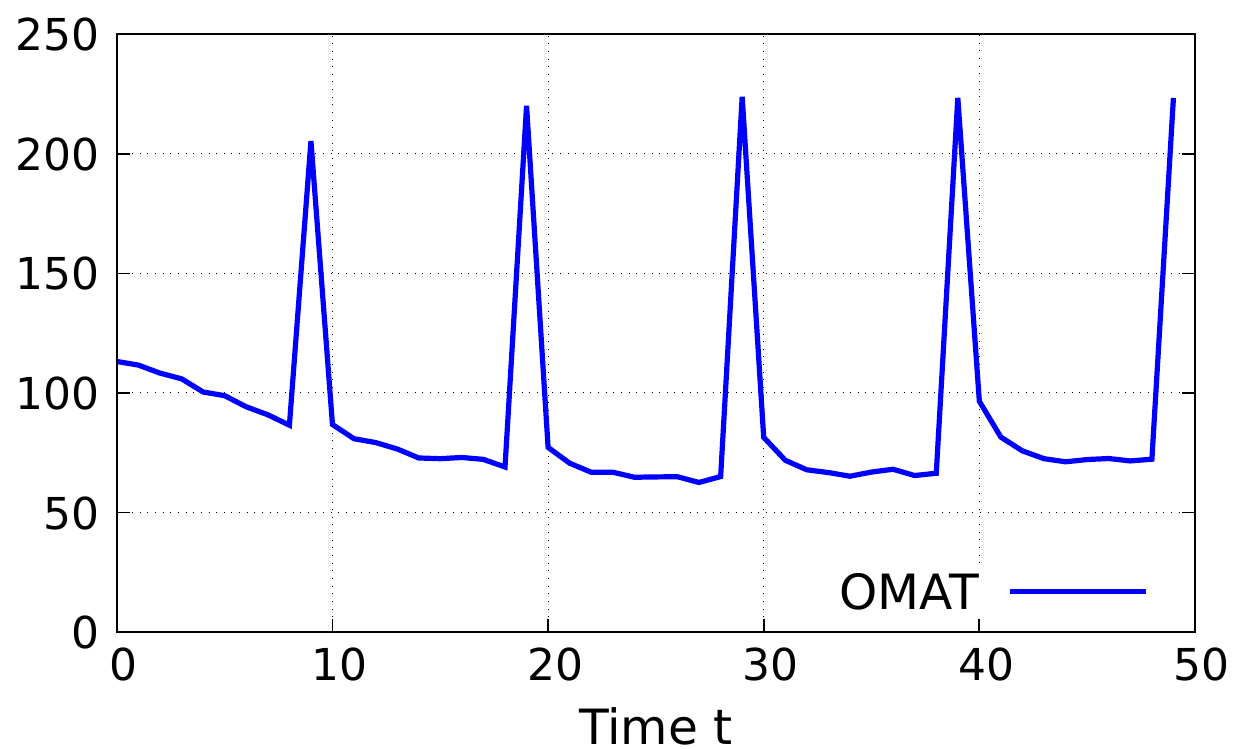}
\caption{OMAT distance estimates.}
\label{fig1a-2}
\end{subfigure}
\hskip1.5cm
\begin{subfigure}{.4\textwidth}
\centering
\includegraphics[width=7.8cm,height=4.5cm]{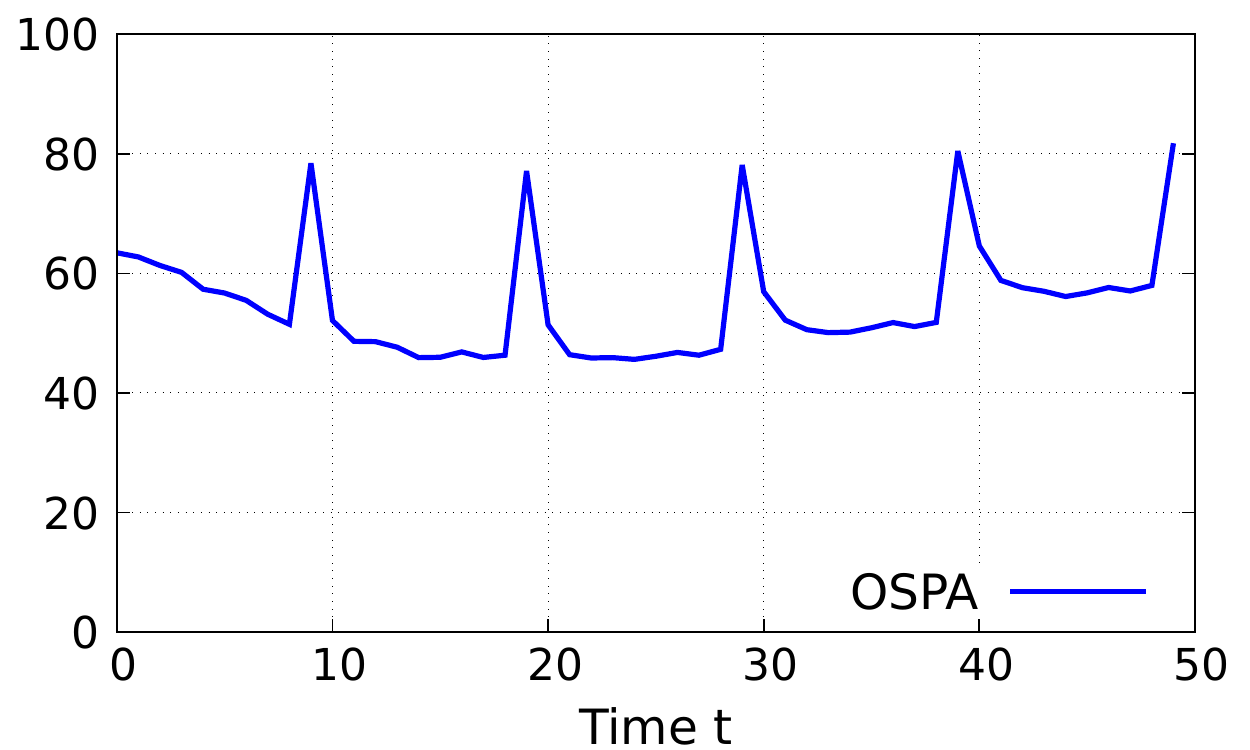}

\caption{OSPA distance estimates.} 
\label{fig1b-2}
\end{subfigure}
\caption{Miss-distance performance evaluation.} 
\label{fig1-2} 
\end{figure}

Figure~\ref{fig1-2} presents miss-distance performance estimates
 for the experiment of Figure~\ref{fig1}, using the 
$L^2$-Optimal Mass Transfer (OMAT, \cite{hoffman}) metric,
and the $L^2$-Optimal Subpattern Assignment (OSPA, \cite{schumacher})
 metric with threshould $c=100$, which solves the inconsistencies
 encountered with the OMAT metric and 
 takes into account differences in cardinalities.  

\begin{figure}[H]
\begin{subfigure}{.4\textwidth}
\centering
\includegraphics[width=7.8cm,height=4.5cm]{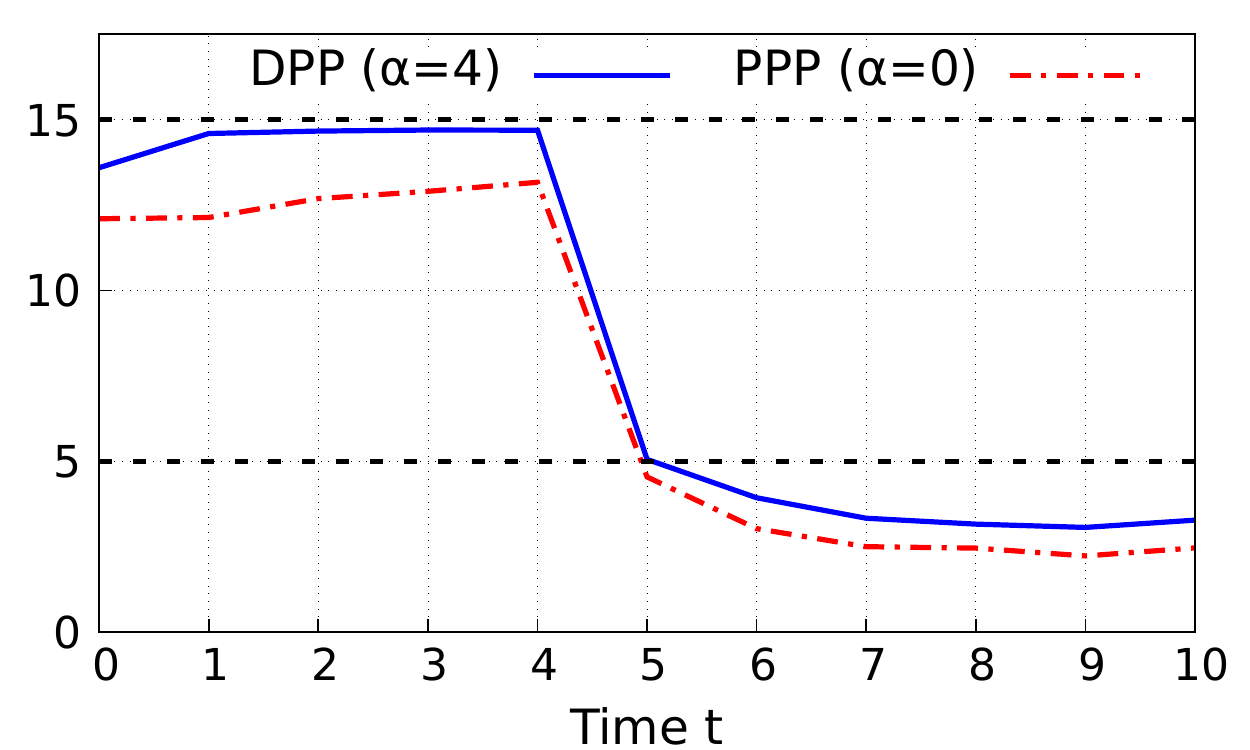}
\caption{Graph with $P_p=50$ and $P_b=40$.}
\label{fig2b}
\end{subfigure}  
\hskip1.5cm
\begin{subfigure}{.4\textwidth}
  \centering
  \includegraphics[width=7.8cm,height=4.5cm]{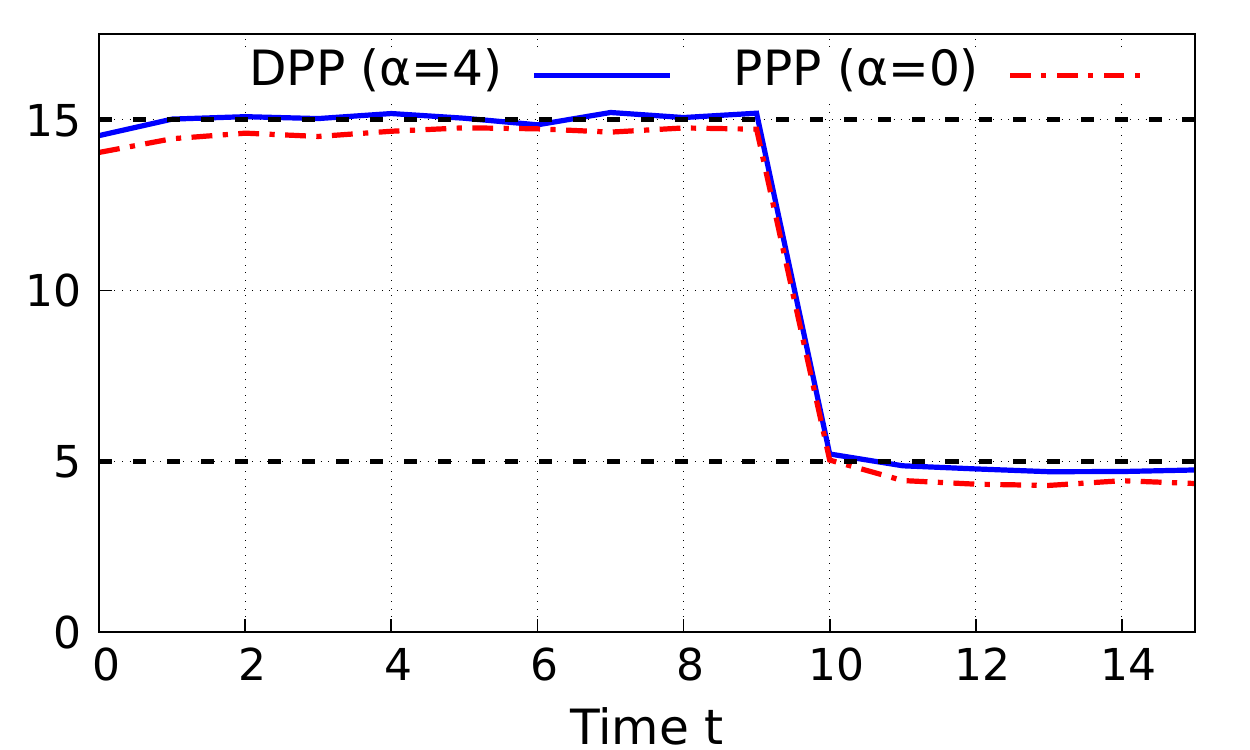}
  \caption{Graph with $P_p=50$ and $P_b=60$.}
\label{fig2a}
\end{subfigure}
\caption{Target count estimates from $15$ to $5$ targets, $p_d=0.9$ and $1$ to $0.3$ clutter points.}
\label{fig2}
\end{figure}

In Figure~\ref{fig2} we compare the robustness of the DPP and PPP-PHD filters
when both filters are subjected to sudden death in the number of targets in a single domain of size $100$m by $100$m, beginning with $15$ targets at the first time step.

\bigskip

Figure~\ref{fig2}-(\subref{fig2b}) uses $300$ Monte Carlo runs, while
Figure~\ref{fig2}-(\subref{fig2a}) 
relies on $200$ Monte Carlo runs. The runtime of each Monte Carlo run spans from time $t=0$ to time $t=15$, 
and the probability of survival is $p_s:=1$.
The initial $15$ targets are maintained until time $t=9$
when $10$ random targets are forced to die and the remaining $5$ targets survive until the end of the time interval.
In Figure~\ref{fig2} we
take the spatial standard deviations (s.d.) $\sigma_{v_x} = \sigma_{v_y} = 1.0$~m/s$^2$, 
turn-rate noise s.d. $\sigma_{v_{\theta}} = \pi = 180$~rad/s,
with bearing and range distribution s.d. 
$\sigma_{w_{\omega}} = \pi = 180$~rad, 
$\sigma_{w_r} = \sqrt{2}$~m as in Figure~\ref{fig1}, 
with probability of detection $p_d=0.95$, mean clutter count at $1$ up to time
$t=9$ and then at $0.06$ afterwards for Figure~\ref{fig2}-(\subref{fig2a}),
and mean clutter count at $1$ up to time
$t=9$ and then at $0.3$ afterwards for Figure~\ref{fig2}-(\subref{fig2b}).
 At initialization in Figure~\ref{fig2}, we set $N_{\Phi,0}=6000$ and $\gamma_{\Phi,0}=0.2$. Both our
DPP and PPP-PHD filter
implementations use $P_p=50$ resampled particles per target in Figure~\ref{fig2}, $P_b=60$ particles per birth target in Figure~\ref{fig2}-(\subref{fig2a}),
 and $P_b=40$ particles per birth target in Figure~\ref{fig2}-(\subref{fig2b}). 

\begin{figure}[H]
  \begin{subfigure}{.4\textwidth}
  \centering
\includegraphics[width=7.8cm,height=4.5cm]{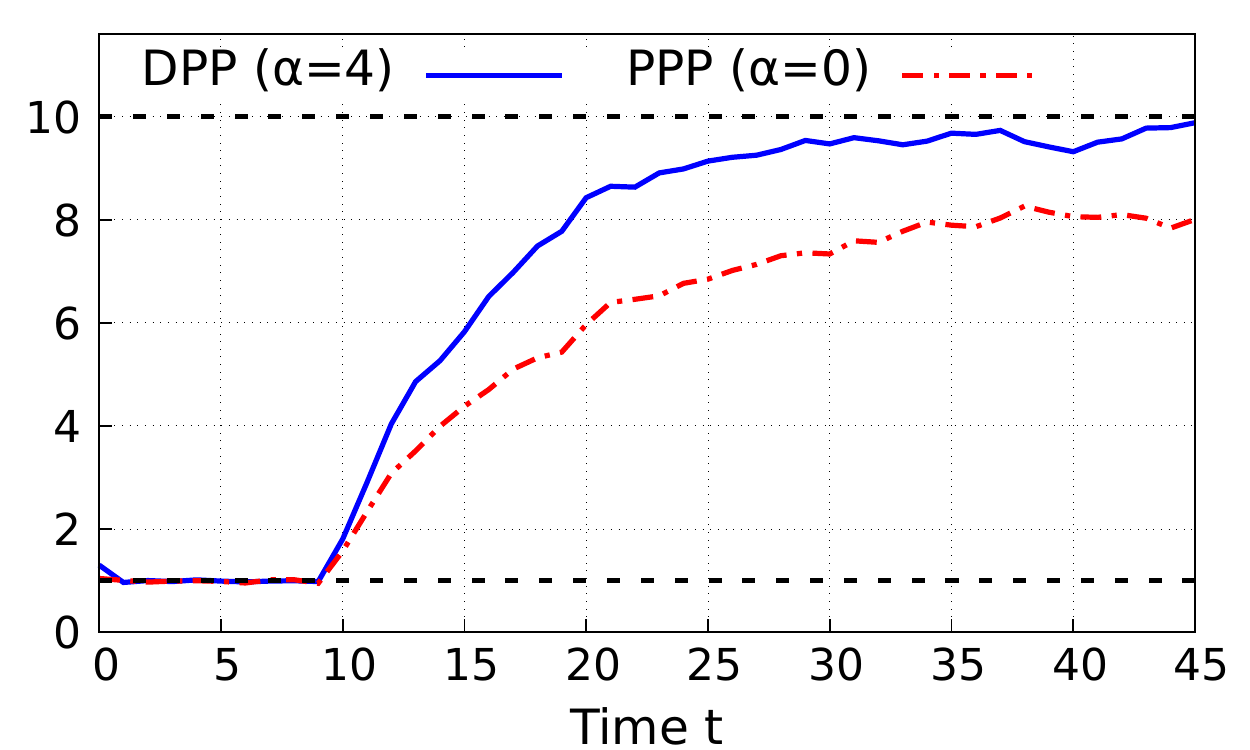}
  \caption{Graph with $P_p=40$ and $P_b=9$.}
\label{fig3a}
  \end{subfigure}
\hskip1.5cm
\begin{subfigure}{.4\textwidth}
  \centering
\includegraphics[width=7.8cm,height=4.5cm]{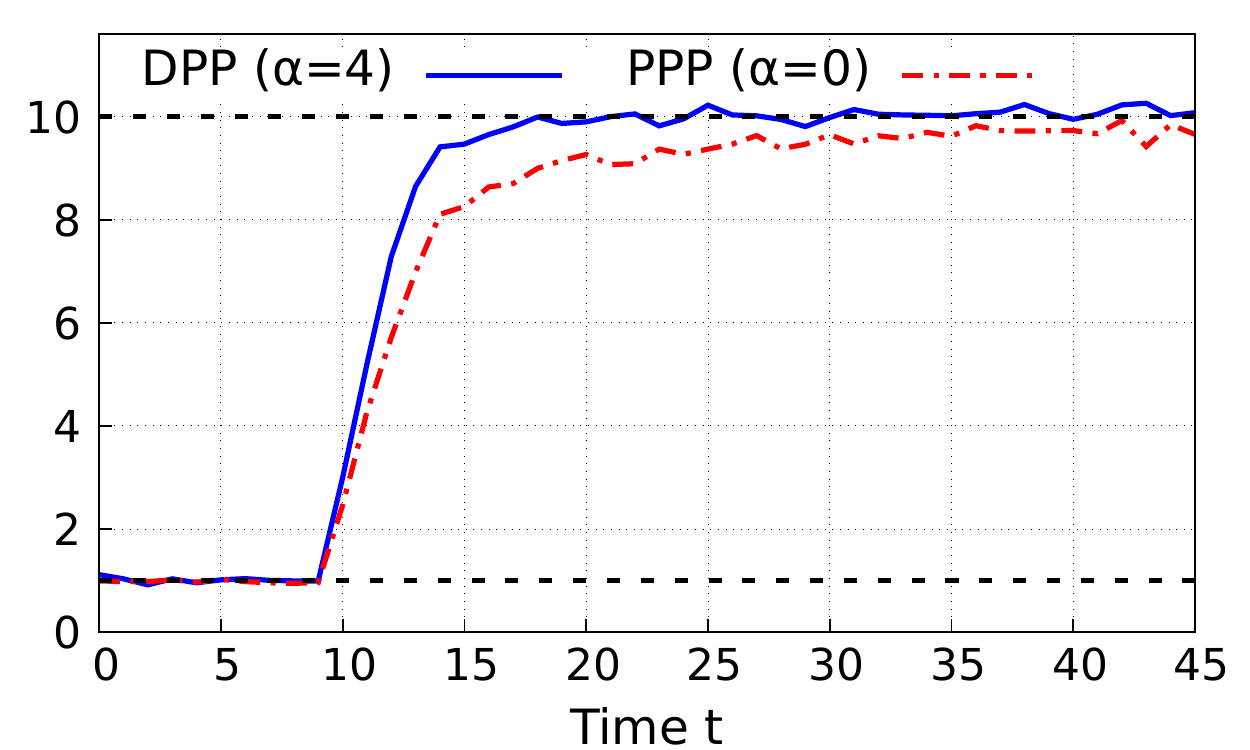}
  \caption{Graph with $P_p=50$ and $P_b=15$.}
\label{fig3b}
\end{subfigure}
\caption{Target count estimates from $0$ to $10$ targets, $p_d=0.9$ and $0$ to $5$ clutter points.}
\label{fig3}
\end{figure}

In Figure~\ref{fig3} we compare the robustness and performance recovery of the DPP
and PPP-PHD filters when subjected to a rapid birth in the number of targets in a single domain 
 of size $100$~m by $100$~m.
The experiment starts with a single target 
which survives throughout the $45$ time steps, without birth of new targets from time $t=0$ to time $t=9$.
At time $t=10$, $9$ new targets are born centrally distributed within the target space
and survive through the remaining time steps.

Each Monte Carlo run spans $45$ time steps,
with $400$ and $100$ Monte Carlo runs in the
experiments of Figures~\ref{fig3}-(\subref{fig3a}) and \ref{fig3}-(\subref{fig3b})
 respectively. 
In Figure~\ref{fig3} the spatial standard deviations (s.d.) 
$\sigma_{v_x} = \sigma_{v_y} = 1.0$~m/s$^2$, 
turn-rate noise s.d. $\sigma_{v_{\theta}} = \pi = 180$~rad/s, 
and bearing and range s.d. 
$\sigma_{w_{\omega}} = \pi = 180$~rad, 
$\sigma_{w_r} = \sqrt{2}$~m
are the same as in Figure~\ref{fig2}. 
The model generates measurement information from each target with a constant probability of detection $p_d=0.90$, mean clutter count at $0$ up to time
$t=9$ and then at $5$ afterwards
 for Figure~\ref{fig3}-(\subref{fig3a}), 
and mean clutter count at $0.05$ up to time
$t=9$ and then at $5$ afterwards
for Figure~\ref{fig3}-(\subref{fig3b}).
 We set $N_{\Phi,0}=300$ and $\gamma_{\Phi,0}=0.2$ at initialization in Figure~\ref{fig3}.
 Both DPP and PPP-PHD filter implementations use
 $P_b=40$ and $P_b=50$ resampled particles per target
 in Figure~\ref{fig3}-(\subref{fig3a})
 and Figure~\ref{fig3}-(\subref{fig3b}) respectively.
 For the target birth process we set $P_b=9$ and $P_b=15$ particles per birth target in
Figure~\ref{fig3}-(\subref{fig3a}) and Figure~\ref{fig3}-(\subref{fig3b}) respectively.

\appendix

\section{Appendix - Janossy density approximation} 
\label{s4}

Since the corrector terms
$l^{(1)}_{z_{1:m}} ( x )$,
$l^{(1)}_{z_{1:m}} (x;z)$,
$l^{(2)}_{z_{1:m}} (x,y)$,
$l^{(2)}_{z_{1:m}} (x,y;z)$,
$l^{(2)}_{z_{1:m}} (x,y;z,z')$
in \eqref{l1corr}, \eqref{l2corr}
and the kernel update formula \eqref{cjkdf}
have no closed form expression in the determinantal
setting, we propose to use the Janossy density approximations 
\begin{equation} 
\label{approx1} 
  j^{(n)}_\Phi (x_1,\ldots , x_{r-1},x,x_{r+1},\ldots , x_n)
 \simeq J_\Phi(x,x) j^{(n-1)}_\Phi (x_1,\ldots , x_{r-1},x_{r+1},\ldots , x_n)
\end{equation} 
$n\geq 1$, which corresponds to a (Poisson)
first-order approximation, and
\begin{eqnarray} 
  \label{approx2}
  \lefteqn{
    j^{(n)}_\Phi (x_1,\ldots , x_{r-1},x,x_{r+1},\ldots , x_{p-1},y,x_{p+1},\ldots , x_n)
  }
  \\
  \nonumber
  & \simeq & 
  (J_\Phi(x,x)J_\Phi(y,y)- ( J_\Phi(x,y))^2 )
  j^{(n-2)}_\Phi (x_1,\ldots , \hat{x}_r , \ldots , \hat{x}_p ,\ldots , x_n),
\end{eqnarray} 
$n\geq 2$, which corresponds to a 
second-order (determinant) approximation, obtained
from \eqref{eq:janossy} by assuming that 
the off-diagonal entries $J_\Phi(x_i,x_j)$, $i\not= j$, are small. 

\bigskip

This Janossy approximation is specially relevant to 
$\alpha$-determinantal Ginibre point processes (GPP)
which approximate a Poisson point process
when $\alpha \in [-1,0)$ tends to $0$, see \cite{shirai}. 

\begin{prop} 
\label{l1}
Under \eqref{approx1} we have the first-order Poisson approximations 
 $l^{(1)}_{z_{1:m}} (x) \simeq J_\Phi (x,x)$, 
  $m\geq 0$, and
  \begin{equation} 
 \nonumber 
 l^{(1)}_{z_{1:m}} (x;z)
 \simeq
 \frac{
   J_\Phi (x,x)
 }{
  l_c(z)
  +
  \int_\Lambda J_\Phi (u,u) \tilde{l}_d (z | u) \nu ( du)
}  ,
  \end{equation}
$z\in z_{1:m}$, 
  $x\in \Lambda$,
  $m\geq 1$. 
\end{prop}
\begin{Proof}
 By \eqref{djksds1} and \eqref{approx1} we have
\begin{eqnarray} 
\nonumber 
     \Upsilon^{(1)}_{z_{1:m}} (x)
& = &  
\displaystyle \sum_{S\subset \{1,\ldots , m\}}
\sum_{p\geq |S|}
 \frac{q_d^{p-|S|}}{(p-|S|)!}
 \prod_{j \notin S} l_c(z_j) 
 \int_{\Lambda^p}
 j^{(p+1)}_\Phi (x_{1:p} , x ) 
    \prod_{i \in S} \tilde{l}_d (z_i | x_i) 
\nu ( dx_{1:p})
\\
\nonumber
        & \simeq &   
 J_\Phi (x,x) 
    \sum_{p\geq 0}
     \displaystyle \sum_{S\subset \{1,\ldots , m\} \atop |S| \leq p} 
  \frac{q_d^{p-|S|}}{(p-|S|)!}
   \prod_{j \notin S} l_c(z_j) 
     \int_{\Lambda^p}
       j^{(p)}_\Phi (x_{1:p} ) 
    \prod_{i \in S} \tilde{l}_d (z_i | x_i) 
\nu ( dx_{1:p})
     \\
\label{cdsds}
        & = & 
     J_\Phi (x,x) 
     j^{(m)}_\Xi (z_1,\ldots , z_m)
, 
  \end{eqnarray} 
by \eqref{jd}, which yields the approximation
$l^{(1)}_{z_{1:m}} ( x ) \simeq      J_\Phi (x,x)$.
    On the other hand, for $r=1,\ldots ,m$,
using again \eqref{approx1} and \eqref{jd} we have 
  \begin{align}
   \nonumber 
 & j^{(m)}_\Xi (z_1,\ldots , z_m) = 
   \displaystyle
  \frac{\partial_{\delta_{z_1}} }{\partial g}
\cdots
\frac{\partial_{\delta_{z_m}} }{\partial g}
{\cal G}_{\Phi , \Xi } ( \ind,g)_{\mid g=0} 
  \\
   \nonumber 
   & = 
\displaystyle
\sum_{p\geq 0}
\sum_{S\subset \{1,\ldots , m\} \atop
|S|\leq p } 
     \frac{    q_d^{p-|S|}}{(p-|S|)!} 
  \prod_{j \notin S} l_c(z_j) 
    \int_{\Lambda^p} 
    j^{(p)}_\Phi (y_{1:p})
    \prod_{i \in S} \tilde{l}_d (z_i | y_i) 
           \nu (dy_{1:p} )
  \\
   \nonumber 
   & \simeq  
  l_c(z_r)
       \sum_{p\geq 0}
     \displaystyle \sum_{S\subset \{1,\ldots , m\} \setminus \{ r \}
       \atop |S| \leq p} 
  \frac{q_d^{p-|S|}}{(p-|S|)!}
  \prod_{j \notin S} l_c(z_j) 
     \int_{\Lambda^p}
       j^{(p)}_\Phi (x_{1:p} ) 
    \prod_{i \in S } \tilde{l}_d (z_i | x_i) 
\nu ( dx_{1:p})
\\
   \nonumber 
& 
\quad +
\int_\Lambda J_\Phi (x_r,x_r) \tilde{l}_d (z_r | x_r) \nu ( dx_r) 
\sum_{p\geq 0}
     \displaystyle
     \sum_{ S\subset \{1,\ldots , m\} \atop |S| \leq p+1 , r\in S } 
  \frac{q_d^{p+1-|S|}}{(p+1-|S|)!}
  \prod_{j \notin S} l_c(z_j) 
     \int_{\Lambda^p}
    j^{(p)}_\Phi (x_{1:p} ) 
    \prod_{i \in S \setminus \{ r \} } \tilde{l}_d (z_i | x_i) 
\nu ( dx_{1:p})
  \\
   \nonumber 
   & = 
  l_c(z_r)
  \sum_{p\geq 0}
     \displaystyle \sum_{S\subset \{1,\ldots , m\} \setminus \{ r \}
       \atop |S| \leq p} 
  \frac{q_d^{p-|S|}}{(p-|S|)!}
  \prod_{j \notin S} l_c(z_j) 
     \int_{\Lambda^p}
       j^{(p)}_\Phi (x_{1:p}) 
    \prod_{i \in S } \tilde{l}_d (z_i | x_i) 
\nu ( dx_{1:p})
\\
   \nonumber 
& 
\quad +
\int_\Lambda J_\Phi (u,u) \tilde{l}_d (z_r | u) \nu ( du) 
\sum_{p\geq 0}
     \displaystyle
     \sum_{ S\subset \{1,\ldots , m\} \setminus \{ r \} \atop |S| \leq p } 
  \frac{q_d^{p-|S|}}{(p-|S|)!}
  \prod_{j \notin S} l_c(z_j) 
     \int_{\Lambda^p}
    j^{(p)}_\Phi (x_{1:p} ) 
    \prod_{i \in S \setminus \{ r \} } \tilde{l}_d (z_i | x_i) 
\nu ( dx_{1:p})
  \\
   \nonumber 
   & = 
  \left(
  l_c(z_r)
  +
  \int_\Lambda J_\Phi (u,u) \tilde{l}_d (z_r | u) \nu ( du)
  \right)
  \sum_{p\geq 0}
     \displaystyle
     \sum_{ S\subset \{1,\ldots , m\} \setminus \{ r \} \atop |S| \leq p } 
  \frac{q_d^{p-|S|}}{(p-|S|)!}
  \prod_{j \notin S} l_c(z_j) 
     \int_{\Lambda^p}
    j^{(p)}_\Phi (x_{1:p} ) 
    \prod_{i \in S } \tilde{l}_d (z_i | x_i) 
\nu ( dx_{1:p})
  \\
 \label{djkd22}
  & = 
  \left(
  l_c(z_r)
  +
  \int_\Lambda J_\Phi (u,u) \tilde{l}_d (z_r | u) \nu ( du)
  \right)
  j^{(m-1)}_\Xi (z_1,\ldots , z_{r-1},z_{r+1}, \ldots , z_m). 
  \end{align}
 We conclude by taking $z_r=z$ and noting that
 by \eqref{l1corr} and \eqref{cdsds}-\eqref{djkd22} we have 
$$ 
 l^{(1)}_{z_{1:m}} (x;z)
 = 
 \frac{
 \Upsilon^{(1)}_{z_{1:m}\! \setminus z} (x)
 }{j^{(m)}_\Xi (z_{1:m})}
 \simeq
 J_\Phi (x,x)
 \frac{
 j^{(m-1)}_\Xi ( z_{1:m} \! \setminus z ) 
 }{j^{(m)}_\Xi (z_{1:m})}.
$$
\end{Proof} 
\begin{prop} 
\label{l2}
 Under \eqref{approx1}-\eqref{approx2} we have the second-order approximations 
$$
 l^{(2)}_{z_{1:m}} (x,y) \simeq J_\Phi(x,x)J_\Phi(y,y)-J_\Phi(x,y)^2,
\quad 
 l^{(2)}_{z_{1:m}} (x,y;z)
 \simeq
 \frac{
J_\Phi (x,x)J_\Phi (y,y)-J_\Phi (x,y)^2 
 }{ l_c(z) + \int_\Lambda J_\Phi (u,u) \tilde{l}_d (z | u) \nu ( du) 
}
     ,
 $$
 $z \in z_{1:m}$, $x,y\in \Lambda$, $m\geq 1$, and
\[
 l^{(2)}_{z_{1:m}} (x,y;z,z')
 \simeq 
 \frac{
 J_\Phi (x,x)J_\Phi (y,y)-J_\Phi (x,y)^2
 }{          s_c(z) s_c(z') 
     -     \int_{\Lambda^2} J_\Phi (u,v)^2 \tilde{l}_d (z | u) \tilde{l}_d (z' | v) \nu ( du) \nu ( dv) 
}
     ,
     \] 
     $z, z' \in z_{1:m}$, $z\not= z'$, $x,y\in \Lambda$, $m\geq 2$, where
\begin{equation} 
  \label{scz}
  s_c(z)
     : =
     l_c(z)
     + \int_\Lambda J_\Phi (v,v) \tilde{l}_d (z | v)\nu (dv),
     \qquad
     z\in \Lambda.
  \end{equation} 
     \end{prop}
\begin{Proof}
 By \eqref{djksds} and \eqref{approx2} we have 
\begin{align} 
  \nonumber
  &
    \Upsilon^{(2)}_{z_{1:m}} (x,y) 
  = 
  \sum_{p\geq 0}
 \displaystyle \sum_{S\subset \{1,\ldots , m\} \atop |S| \leq p} 
  \frac{q_d^{p-|S|}}{(p-|S|)!}
 \prod_{j \notin S} l_c(z_j) 
 \int_{\Lambda^p}
 j^{(p+2)}_\Phi (x_{1:p}, x, y ) 
    \prod_{i \in S} \tilde{l}_d (z_i | x_i) 
 \nu (dx_{1:p} )
  \\
     \nonumber
   & = 
  (J_\Phi (x,x)J_\Phi (y,y)-J_\Phi (x,y)^2 )
 \sum_{p\geq 0}
     \displaystyle \sum_{S\subset \{1,\ldots , m\} \atop |S| \leq p} 
  \frac{q_d^{p-|S|}}{(p-|S|)!}
 \prod_{j \notin S} l_c(z_j) 
 \int_{\Lambda^p}
 j^{(p)}_\Phi (x_{1:p} )
     \prod_{i \in S} \tilde{l}_d (z_i | x_i) 
  \nu (dx_{1:p} )
      \\
     \label{djklsf}
   & = 
 ( J_\Phi (x,x)J_\Phi (y,y)-J_\Phi (x,y)^2 )
 j^{(m)}_\Xi ( z_{1:m}), 
\end{align} 
and for $r,u=1,\ldots ,m$,
using \eqref{approx1}-\eqref{approx2} and \eqref{jd} we find 
 \begin{align}
   \nonumber 
 & j^{(m)}_\Xi (z_{1:m}) = 
\displaystyle
\sum_{n\geq 0}
\sum_{S\subset \{1,\ldots , m\} \atop
|S|\leq n } 
\prod_{j \notin S} l_c(z_j) 
      \frac{    q_d^{n-|S|}}{(n-|S|)!} 
    \int_{\Lambda^n} 
    j^{(n)}_\Phi (y_{1:n})
    \prod_{i \in S} \tilde{l}_d (z_i | y_i) 
           \nu (dy_{1:n} )
  \\
   \nonumber 
   & \simeq  
  l_c(z_r)
  l_c(z_u)
       \sum_{p\geq 0}
         \displaystyle \sum_{S\subset \{1,\ldots , m\} \setminus \{ r , u \}
       \atop |S| \leq p} 
  \frac{q_d^{p-|S|}}{(p-|S|)!}
 \prod_{j \notin S} l_c(z_j) 
 \int_{\Lambda^p}
 j^{(p)}_\Phi (x_{1:p} ) 
    \prod_{i \in S } \tilde{l}_d (z_i | x_i) 
\nu ( dx_{1:p})
\\
   \nonumber 
& \quad 
 +
  l_c(z_r)
  \int_\Lambda J_\Phi (v,v) \tilde{l}_d (z_u | v)\nu (dv)
  \\
  \nonumber
  & \qquad \times
  \sum_{p\geq 0}
     \int_{\Lambda^p}
     \displaystyle \sum_{S\subset \{1,\ldots , m\} \setminus \{ r \} \atop |S| \leq p+1,
   u \in S } 
  \frac{q_d^{p+1-|S|}}{(p+1-|S|)!}
  j^{(p)}_\Phi (x_{1:p} )
  \prod_{j \notin S } l_c(z_j) 
           \prod_{i \in S \setminus \{ u\} } \tilde{l}_d (z_i | x_i) 
\nu ( dx_{1:p})
\\
   \nonumber 
& 
\quad +
  l_c(z_u)
  \int_\Lambda J_\Phi (v,v) \tilde{l}_d (z_r | v)\nu (dv)
  \\
  \nonumber
  & \qquad \times
  \sum_{p\geq 0}
     \int_{\Lambda^p}
     \displaystyle \sum_{S\subset \{1,\ldots , m\} \setminus \{ u \} \atop |S| \leq p+1,
   r \in S } 
  \frac{q_d^{p+1-|S|}}{(p+1-|S|)!}
 j^{(p)}_\Phi (x_{1:p})
 \prod_{j \notin S } l_c(z_j) 
           \prod_{i \in S \setminus \{ r \} } \tilde{l}_d (z_i | x_i) 
\nu ( dx_{1:p})
\\
   \nonumber 
& 
\quad +
\int_\Lambda
 ( J_\Phi (x_r,x_r)J_\Phi (x_u,x_u)-J_\Phi (x_r,x_u)^2 )
 \tilde{l}_d (z_r | x_r) \tilde{l}_d (z_u | x_u) \nu ( dx_r) \nu ( dx_u)
      \\
  \nonumber
  & \qquad \times
\sum_{p\geq 0}
     \int_{\Lambda^p}
     \displaystyle
     \sum_{ S\subset \{1,\ldots , m\} \atop |S| \leq p+2 , r\in S } 
  \frac{q_d^{p+2-|S|}}{(p+2-|S|)!}
  j^{(p)}_\Phi (x_{1:p} )
  \prod_{j \notin S} l_c(z_j) 
        \prod_{i \in S \setminus \{ r , u \} } \tilde{l}_d (z_i | x_i) 
\nu ( dx_{1:p})
  \\
   \nonumber 
   & = 
  l_c(z_r)
  l_c(z_u)
  \sum_{p\geq 0}
     \int_{\Lambda^p}
     \displaystyle \sum_{S\subset \{1,\ldots , m\} \setminus \{ r , u\}
       \atop |S| \leq p} 
  \frac{q_d^{p-|S|}}{(p-|S|)!}
  j^{(p)}_\Phi (x_{1:p} )
  \prod_{j \notin S} l_c(z_j) 
           \prod_{i \in S } \tilde{l}_d (z_i | x_i) 
\nu ( dx_{1:p})
\\
   \nonumber 
& 
\quad +
  l_c(z_r)
  \int_\Lambda J_\Phi (v,v) \tilde{l}_d (z_u | v)\nu (dv)
  \sum_{p\geq 0}
     \int_{\Lambda^p}
     \displaystyle \sum_{S\subset \{1,\ldots , m\} \setminus \{ r , u \} \atop |S| \leq p,
   u \in S } 
  \frac{q_d^{p-|S|}}{(p-|S|)!}
 j^{(p)}_\Phi (x_{1:p} )
 \prod_{j \notin S } l_c(z_j) 
           \prod_{i \in S } \tilde{l}_d (z_i | x_i) 
\nu ( dx_{1:p})
\\
   \nonumber 
& 
\quad +
  l_c(z_u)
  \int_\Lambda J_\Phi (v,v) \tilde{l}_d (z_r | v)\nu (dv)
 \sum_{p\geq 0}
     \int_{\Lambda^p}
     \displaystyle \sum_{S\subset \{1,\ldots , m\} \setminus \{ r, u \} \atop |S| \leq p,
   r \in S } 
  \frac{q_d^{p-|S|}}{(p-|S|)!}
  j^{(p)}_\Phi (x_{1:p} )
 \prod_{j \notin S } l_c(z_j) 
          \prod_{i \in S } \tilde{l}_d (z_i | x_i) 
\nu ( dx_{1:p})
\\
   \nonumber 
& 
\quad +
\int_{\Lambda^2}
 ( J_\Phi (u,u)J_\Phi (v,v)-J_\Phi (u,v)^2 )
 \tilde{l}_d (z_r | u) \tilde{l}_d (z_u | v) \nu ( du) \nu ( dv) 
  \\
  \nonumber
  & \qquad \times
    \sum_{p\geq 0}
     \int_{\Lambda^p}
     \sum_{ S\subset \{1,\ldots , m\} \setminus \{ r , u \} \atop |S| \leq p } 
  \frac{q_d^{p-|S|}}{(p-|S|)!}
    j^{(p)}_\Phi (x_{1:p} )
   \prod_{j \notin S} l_c(z_j) 
       \prod_{i \in S } \tilde{l}_d (z_i | x_i) 
\nu ( dx_{1:p})
  \\
   \nonumber 
   & = 
  \left(
  l_c(z_r)
  l_c(z_u)
  +
    l_c(z_r)
  \int_\Lambda J_\Phi (v,v) \tilde{l}_d (z_u | v)\nu (dv)
 +  l_c(z_u)
  \int_\Lambda J_\Phi (v,v) \tilde{l}_d (z_r | v)\nu (dv)
  \right.
  \\
  \nonumber
  & \left.
  \quad +
\int_{\Lambda^2}
 (J_\Phi (u,u)J_\Phi (v,v)-J_\Phi (u,v)^2 )
 \tilde{l}_d (z_r | u) \tilde{l}_d (z_u | v) \nu ( du) \nu ( dv) 
 \right)
 \\
  \nonumber
  & \qquad \times
 \sum_{p\geq 0}
     \int_{\Lambda^p}
     \displaystyle
     \sum_{ S\subset \{1,\ldots , m\} \setminus \{ r , u \} \atop |S| \leq p } 
  \frac{q_d^{p-|S|}}{(p-|S|)!}
  j^{(p)}_\Phi (x_{1:p} )
  \prod_{j \notin S} l_c(z_j) 
        \prod_{i \in S } \tilde{l}_d (z_i | x_i) 
\nu ( dx_{1:p})
  \\
  \nonumber 
  & = 
  \left(
  l_c(z_r)
  l_c(z_u)
  +
    l_c(z_r)
  \int_\Lambda J_\Phi (v,v) \tilde{l}_d (z_u | v)\nu (dv)
  \right.
  \\
  \nonumber
  &   \left.
  \quad +  l_c(z_u)
  \int_\Lambda J_\Phi (v,v) \tilde{l}_d (z_r | v)\nu (dv)
  +
\int_{\Lambda^2}
 ( J_\Phi (u,u)J_\Phi (v,v)-J_\Phi (u,v)^2 )
 \tilde{l}_d (z_r | u) \tilde{l}_d (z_u | v) \nu ( du) \nu ( dv) 
 \right)
  \\
  \nonumber
  & \qquad \times
 j^{(m-2)}_\Xi (z_1,\ldots , z_{r-1},z_{r+1}, \ldots , z_{u-1},z_{u+1}, \ldots , z_m)
  \\
  \label{dkjdsf} 
& \quad
= 
  \left(
   s_c(z_r)
  s_c(z_u)
-
  \int_{\Lambda^2} J_\Phi (u,v)^2 \tilde{l}_d (z_r | u) \tilde{l}_d (z_u | v) \nu ( du) \nu ( dv) 
 \right)
  \\
  \nonumber
  & \qquad \times
  j^{(m-2)}_\Xi (z_1,\ldots , z_{r-1},z_{r+1}, \ldots , z_{u-1},z_{u+1}, \ldots , z_m)
 .
\end{align} 
 We conclude by taking $(z_r,z_u)=(z,z')$ and noting that by
 \eqref{l2corr.2} and \eqref{djklsf}-\eqref{dkjdsf} we have 
\begin{eqnarray*} 
     l^{(2)}_{z_{1:m}} (x,y,z,z')
 & = & 
 \frac{\Upsilon^{(2)}_{z_{1:m}\! \setminus \{ z,z'\}} (x,y)}{j^{(m)}_\Xi (z_{1:m})}
  \\
 & \simeq & 
 ( J_\Phi (x,x)J_\Phi (y,y)-J_\Phi (x,y)^2 )
 \frac{
 j^{(m-2)}_\Xi ( z_{1:m}\! \setminus \{ z,z' \} ) 
 }{j^{(m)}_\Xi (z_{1:m})}
\\
 & \simeq & 
 \frac{
 ( J_\Phi (x,x)J_\Phi (y,y)-J_\Phi (x,y)^2 )
 }{          s_c(z) s_c(z') 
     -     \int_{\Lambda^2} J_\Phi (u,v)^2 \tilde{l}_d (z | u) \tilde{l}_d (z' | v) \nu ( du) \nu ( dv) 
}
     ,
\end{eqnarray*} 
 $z, z' \in z_{1:m}$, $z\not= z'$, $m\geq 2$.
 Similarly, by \eqref{l2corr} and \eqref{djkd22}, \eqref{djklsf} we also have
\begin{eqnarray*} 
      l^{(2)}_{z_{1:m}} (x,y;z)
  & = &  
 \frac{\Upsilon^{(2)}_{z_{1:m}\! \setminus z } (x,y)}{j^{(m)}_\Xi (z_{1:m})}
  \\
 & \simeq & 
 ( J_\Phi (x,x)J_\Phi (y,y)-J_\Phi (x,y)^2 )
 \frac{
 j^{(m-1)}_\Xi ( z_{1:m}\! \setminus z ) 
 }{j^{(m)}_\Xi (z_{1:m})}
\\
 & \simeq & 
 \frac{
J_\Phi (x,x)J_\Phi (y,y)-J_\Phi (x,y)^2 
 }{l_c(z) + \int_\Lambda J_\Phi (u,u) \tilde{l}_d (z | u) \nu ( du) 
},
 \qquad
 z \in z_{1:m}, \quad
 m\geq 1.
\end{eqnarray*} 
\end{Proof}
As a consequence of \eqref{rho22} and Proposition~\ref{l2},
the second-order conditional factorial moment density of 
$\Phi$ given that $\Xi =z_{1:m}=(z_1,\ldots , z_m)$ 
will be approximated as 
\begin{eqnarray} 
\label{djkls2}
\lefteqn{
  \rho^{(2)}_{\Phi \mid \Xi =z_{1:m}} (x,y)  
    \simeq 
    q_d^2 ( J_\Phi (x,x)J_\Phi (y,y)-J_\Phi (x,y)^2 )
}
\\
    \nonumber
 &   & +     q_d  \sum_{z\in z_{1:m}} 
    \frac{
         ( J_\Phi (x,x)J_\Phi (y,y)-J_\Phi (x,y)^2 )
 \big( \tilde{l}_d (z | x)
      +
 \tilde{l}_d (z | y) \big)}{
   s_c(z)
    }
\\
     \nonumber
  &   &
        + 
     \sum_{z,z'\in z_{1:m} \atop z \not= z'}
     \frac{  
      ( J_\Phi (x,x)J_\Phi (y,y)-J_\Phi (x,y)^2 )
 \tilde{l}_d (z | x)
 \tilde{l}_d (z' | y)
}{
         s_c(z) s_c(z')
         -
  \int_{\Lambda^2} J_\Phi (u,v)^2 \tilde{l}_d (z | u) \tilde{l}_d (z' | v) \nu (du) \nu (dv) 
 }
     , 
\end{eqnarray} 
 $m\geq 0$,  
 with $\rho^{(2)}_{\Phi \mid \Xi =z_{1:m}} (x,x):=0$, $x\in \Lambda$.  

\begin{prop}
\label{thasd}
  The (approximate) kernel update formula is given by 
    \begin{eqnarray} 
\nonumber 
 K_{\Phi \mid \Xi =z_{1:m}} (x,y)^2
  &     \simeq
       &
       q^2_d J_\Phi (x,y)^2 
    +   q_d  J_\Phi (x,y)^2 
 \sum_{z\in z_{1:m}}
    \frac{  \big( \tilde{l}_d (z | x) + \tilde{l}_d (z | y) \big) }{
          s_c ( z ) 
}
     \displaystyle
      \\
      \nonumber
      & 
 &     \displaystyle
     +
J_\Phi ( x , x ) J_\Phi ( y , y )      \sum_{z,z'\in z_{1:m}}
    \frac{ \tilde{l}_d (z | x)  \tilde{l}_d (z' | y)}{
           s_c(z) s_c(z')
      }
      \\
      \nonumber
      & 
  & + \sum_{z,z'\in z_{1:m} \atop z \not= z'} 
 \frac{      (J_\Phi (x,y)^2 - J_\Phi (x,x)J_\Phi (y,y) )   \tilde{l}_d (z | x)
 \tilde{l}_d (z' | y)
}{
      s_c(z) s_c(z')
  -
  \int_{\Lambda^2} J_\Phi (u,v)^2 \tilde{l}_d (z | u) \tilde{l}_d (z' | v) \nu (du) \nu (dv) 
 }
   ,
\end{eqnarray} 
    $m\geq 0$,
    $x,y\in \Lambda$.
\end{prop}
\begin{Proof}
 By \eqref{rho11} and Proposition~\ref{l1},
 we have the approximation
\begin{eqnarray} 
  \label{djkls1}
  \mu^{(1)}_{\Phi \mid \Xi =z_{1:m}} (x) 
& = &  
  q_d l^{(1)}_{z_{1:m}} ( x )
  +
  \sum_{z\in z_{1:m}}
  \tilde{l}_d(x\mid z)
  l^{(1)}_{z_{1:m}} (x;z)
\\
\nonumber 
  &   \simeq & 
q_d J_\Phi ( x , x )
 + \sum_{z\in z_{1:m}}
    \frac{ J_\Phi ( x , x ) \tilde{l}_d (z | x)}{
     l_c (z ) + \int_{\Lambda} \tilde{l}_d (z | u ) J_\Phi ( u,u) \nu ( du )  }
    , \qquad
m\geq 0,  
\end{eqnarray} 
 hence by \eqref{djkls2} and \eqref{djkls1}, we find 
\begin{align} 
  \nonumber
&  \rho^{(2)}_{\Phi \mid \Xi =z_{1:m}} (x,y)
  -
       \mu^{(1)}_{\Phi, \Xi =z_{1:m}} (x) 
     \mu^{(1)}_{\Phi, \Xi =z_{1:m}} (y) 
\\
\nonumber
        & \simeq     
 -   q_d J_\Phi ( x , x )J_\Phi ( y , y )
    \bigg( 
    q_d 
 + \sum_{z\in z_{1:m}}
    \frac{ \tilde{l}_d (z | x) + \tilde{l}_d (z | y)}{
 s_c(z) } 
\bigg) 
- 
J_\Phi ( x , x )J_\Phi ( y , y )
\sum_{z,z'\in z_{1:m}}
    \frac{     \tilde{l}_d (z | x)\tilde{l}_d (z' | y)}{
      s_c (z ) 
     s_c (z' ) 
    }
\\
\nonumber
&  
\quad +
 q_d^2    ( J_\Phi (x,x)J_\Phi (y,y)-J_\Phi (x,y)^2 )
   + q_d 
 ( J_\Phi (x,x)J_\Phi (y,y)-J_\Phi (x,y)^2 )
 \sum_{z\in z_{1:m}} 
 \frac{
 \tilde{l}_d (z | x) + \tilde{l}_d (z | y)}{
  s_c(z) 
   }
\\
     \nonumber
& \quad + 
 \sum_{z,z'\in z_{1:m} \atop z\not= z'} 
 \frac{      ( J_\Phi (x,x)J_\Phi (y,y)-J_\Phi (x,y)^2 )
 \tilde{l}_d (z | x)
 \tilde{l}_d (z' | y)
}{
    s_c(z)
  s_c(z')
   -
  \int_{\Lambda^2} J_\Phi (u,v)^2 \tilde{l}_d (z | u) \tilde{l}_d (z' | v) \nu (du) \nu (dv) 
 }
 \\
\nonumber
        & =  - q^2_d J_\Phi (x,y)^2 
  -   q_d J_\Phi (x,y)^2
 \sum_{z\in z_{1:m}}
    \frac{ \tilde{l}_d (z | x) + \tilde{l}_d (z | y) }{
       s_c(z) 
} 
\\
\nonumber
& \quad  - J_\Phi ( x , x ) J_\Phi ( y , y ) 
     \sum_{z,z'\in z_{1:m}}
    \frac{ \tilde{l}_d (z | x)  \tilde{l}_d (z' | y)}{
                 s_c (z ) 
    s_c (z' ) 
    }
  + \sum_{z,z'\in z_{1:m} \atop z\not= z'} 
 \frac{     ( J_\Phi (x,x)J_\Phi (y,y)-J_\Phi (x,y)^2 )
 \tilde{l}_d (z | x)
 \tilde{l}_d (z' | y)
}{
     s_c(z)
    s_c(z')
  -
  \int_{\Lambda^2} J_\Phi (u,v)^2 \tilde{l}_d (z | u) \tilde{l}_d (z' | v) \nu (du) \nu (dv)  
 }
 ,
\end{align} 
$m\geq 0$, and we conclude by \eqref{kxy}, i.e. 
$$ 
  ( K_{\Phi \mid \Xi =z_{1:m}} (x,y) )^2 =
     \mu^{(1)}_{\Phi \mid \Xi =z_{1:m}} (x) 
     \mu^{(1)}_{\Phi \mid \Xi =z_{1:m}} (y) 
-
  \rho^{(2)}_{\Phi \mid \Xi =z_{1:m}} (x,y)
$$ 
  and \eqref{cnjdf}. 
\end{Proof}
The next result, which provides
an approximation formula for the
posterior covariance of Proposition~\ref{djlksd}, 
is a consequence of Proposition~\ref{thasd}
and \eqref{djkls1}. 
  
\begin{corollary} 
  \label{fdsfsfd}
   Under \eqref{approx1}-\eqref{approx2} the posterior covariance
 of $\Phi$ given that $\Xi =z_{1:m}=(z_1,\ldots , z_m)$ is approximated as 
\begin{eqnarray}
 \label{cnjdf} 
 \lefteqn{
   c^{(2)}_{\Phi \mid \Xi =z_{1:m}} (A,B)  
  \simeq 
q_d \int_{A\cap B} J_\Phi ( x , x ) \nu ( dx ) 
    - q^2_d \int_{A\times B} J_\Phi (x,y)^2 \nu ( dx ) \nu ( dy )
 }
 \\
\nonumber
      & & 
      -   q_d 
 \sum_{z\in z_{1:m}}
 \frac{ 1 }{
       s_c (z )  }
 \int_{A\times B} J_\Phi (x,y)^2 \big( \tilde{l}_d (z | x) + \tilde{l}_d (z | y) \big)
   \nu ( dx ) \nu ( dy )
     \displaystyle
      \\
\nonumber
      & & 
     \displaystyle
 + \sum_{z\in z_{1:m}}
 \frac{1}{ s_c (z ) } 
 \bigg(
 \int_{A\cap B} \tilde{l}_d (z | x) J_\Phi ( x , x )  \nu (dx)
 -  
 \frac{ \int_A \tilde{l}_d (z | x)  J_\Phi ( x , x ) \nu (dx )
   \int_B \tilde{l}_d (z | y) J_\Phi ( y , y )  \nu (dy)
 }{
           s_c (z ) 
    }
 \bigg)
    \\
      \nonumber
 & & 
      \hskip-0.1cm
      + \hskip-0.3cm
      \sum_{z,z'\in z_{1:m} \atop z\not= z'} 
      \hskip-0.4cm
      \frac{        \int_{\Lambda^2} ( J_\Phi (x,x)J_\Phi (y,y)-J_\Phi (x,y)^2 )
 \tilde{l}_d (z | x)
 \tilde{l}_d (z' | y) \nu(dx) \nu (dy) 
}{
     s_c(z)
    s_c(z')
  -
  \int_{\Lambda^2} J_\Phi (u,v)^2 \tilde{l}_d (z | u) \tilde{l}_d (z' | v) \nu (du) \nu (dv)  
 }, \qquad m \geq 0.
\end{eqnarray} 
\end{corollary} 

\subsubsection*{Conclusion}
Our observations have shown that the performance of the
 multi-target tracking PPP-based
 standard PHD filter
 is degraded in the presence of target interaction such as repulsion.
 To address this issue, we have constructed 
 a second-order DPP-based PHD filter based on Determinantal Point Processes 
 which are able to model repulsion between targets, 
 and can propagate variance and covariance information
 in addition to first-order target count estimates. 
 We have derived posterior moment formulas for the estimation
 of DPPs after thinning and superposition with a Poisson Point
 Process (PPP), based on suitable approximation formulas.
 Our numerical experiments include
 an assessment of the spooky effect on disjoint domains,
 with negative correlation estimates which are consistent with
 the nature of DPPs.
 We have also compared the robustness and performance recovery of the DPP and PPP-PHD
 filters when subjected to sudden changes in target numbers.

\footnotesize 
 
\setcitestyle{numbers}

\def\cprime{$'$} \def\polhk#1{\setbox0=\hbox{#1}{\ooalign{\hidewidth
  \lower1.5ex\hbox{`}\hidewidth\crcr\unhbox0}}}
  \def\polhk#1{\setbox0=\hbox{#1}{\ooalign{\hidewidth
  \lower1.5ex\hbox{`}\hidewidth\crcr\unhbox0}}} \def\cprime{$'$}


\begin{thebibliography}{}

\bibitem[\protect\astroncite{Brezis}{1983}]{brezis2}
Brezis, H. (1983).
\newblock {\em Analyse fonctionnelle}.
\newblock Collection Math\'ematiques Appliqu\'ees pour la Ma\^\i trise.
  [Collection of Applied Mathematics for the Master's Degree]. Masson, Paris.

\bibitem[\protect\astroncite{Clark and de~Melo}{2018}]{demelo}
Clark, D. and de~Melo, F. (2018).
\newblock A linear-complexity second-order multi-object filter via factorial
  cumulants.
\newblock In {\em 2018 21st International Conference on Information Fusion
  (FUSION)}, pages 1250--1259.

\bibitem[\protect\astroncite{Clark et~al.}{2016}]{cdh}
Clark, D., Delande, E., and Houssineau, J. (2016).
\newblock Basic concepts for multi-object estimation.
\newblock Lecture notes, Heriot-Watt University.

\bibitem[\protect\astroncite{Clark and Houssineau}{2012}]{Faa}
Clark, D. and Houssineau, J. (2012).
\newblock Faa di {B}runo's formula for {G}ateaux differentials and interacting
  stochastic population processes.
\newblock Preprint arXiv:1202.0264v4.

\bibitem[\protect\astroncite{Daley and Vere-Jones}{2003}]{daley}
Daley, D.~J. and Vere-Jones, D. (2003).
\newblock {\em An introduction to the theory of point processes. {V}ol. {I}}.
\newblock Probability and its Applications. Springer-Verlag, New York.

\bibitem[\protect\astroncite{de~Melo and Maskell}{2019}]{maskell}
de~Melo, F. and Maskell, S. (2019).
\newblock A {CPHD} approximation based on a discrete-gamma cardinality model.
\newblock {\em IEEE Trans. Signal Processing}, 67(2):336--350.

\bibitem[\protect\astroncite{Decreusefond et~al.}{2016}]{dfpt2}
Decreusefond, L., Flint, I., Privault, N., and Torrisi, G. (2016).
\newblock Determinantal point processes.
\newblock In Peccati, G. and Reitzner, M., editors, {\em Stochastic Analysis
  for {P}oisson Point Processes: {M}alliavin Calculus, {W}iener-{I}t{\^o} Chaos
  Expansions and Stochastic Geometry}, volume~7 of {\em Bocconi \& Springer
  Series}, pages 311--342, Berlin. Springer.

\bibitem[\protect\astroncite{Delande et~al.}{2014}]{duhc}
Delande, E., {\"U}ney, M., Houssineau, J., and Clark, D. (2014).
\newblock Regional variance for multi-object filtering.
\newblock {\em IEEE Trans. Signal Processing}, 62(13):3415--3428.

\bibitem[\protect\astroncite{Fr{\"a}nken et~al.}{2009}]{franken}
Fr{\"a}nken, D., Schmidt, M., and Ulmke, M. (2009).
\newblock Spooky action at a distance in the cardinalized probability
  hypothesis density filter.
\newblock {\em IEEE Transactions on Aerospace and Electronic Systems},
  45(4):1657--1664.

\bibitem[\protect\astroncite{Georgii and Yoo}{2005}]{georgiiyoo}
Georgii, H. and Yoo, H. (2005).
\newblock Conditional intensity and {G}ibbsianness of determinantal point
  processes.
\newblock {\em J. Stat. Phys.}, 118(1-2):55--84.

\bibitem[\protect\astroncite{Hoffman and Mahler}{2004}]{hoffman}
Hoffman, J. and Mahler, R. (2004).
\newblock Multitarget {B}ayes filtering via first-order multitarget moments.
\newblock {\em IEEE Transactions on Systems, Man, and Cybernetics - Part A:
  Systems and Humans}, 34(3):327--336.

\bibitem[\protect\astroncite{Hough et~al.}{2009}]{hough}
Hough, J.-B., Krishnapur, M., Peres, Y., and Vir{\'a}g, B. (2009).
\newblock {\em Zeros of {G}aussian analytic functions and determinantal point
  processes}, volume~51 of {\em University Lecture Series}.
\newblock American Mathematical Society, Providence, RI.

\bibitem[\protect\astroncite{Jorquera et~al.}{2018}]{jorquera}
Jorquera, F., Hern{\'a}ndez, S., and Vergara, D. (2018).
\newblock Multi target tracking using determinantal point processes.
\newblock In {\em Progress in Pattern Recognition, Image Analysis, Computer
  Vision, and Applications}, volume 10657 of {\em Lecture Notes in Computer
  Science}, pages 323--330. Springer.

\bibitem[\protect\astroncite{Jorquera et~al.}{2019}]{jorquera2}
Jorquera, F., Hern{\'a}ndez, S., and Vergara, D. (2019).
\newblock Probability hypothesis density filter using determinantal point
  processes for multi object tracking.
\newblock {\em Computer Vision and Image Understanding}, 183:33--41.

\bibitem[\protect\astroncite{Koch}{2018}]{koch}
Koch, W. (2018).
\newblock On anti-symmetry in multiple target tracking.
\newblock In {\em 2018 21st International Conference on Information Fusion
  (FUSION)}, pages 957--964.

\bibitem[\protect\astroncite{Li et~al.}{2017}]{litiancheng}
Li, T., Corchado, J., Sun, S., and Fan, H. (2017).
\newblock Multi-{EAP}: Extended {EAP} for multi-estimate extraction for
  {SMC-PHD} filter.
\newblock {\em Chinese Journal of Aeronautics}, 30(1):368--379.

\bibitem[\protect\astroncite{Li et~al.}{2013}]{sattar}
Li, T., Sattar, T.~P., Han, Q., and Sun, S. (2013).
\newblock Roughening methods to prevent sample impoverishment in the particle
  {PHD} filter.
\newblock In {\em Proceedings of the 16th International Conference on
  Information Fusion}, pages 17--22. IEEE, Istanbul.

\bibitem[\protect\astroncite{Lund and Rudemo}{2000}]{lund}
Lund, J. and Rudemo, M. (2000).
\newblock Models for point processes observed with noise.
\newblock {\em Biometrika}, 87(2):235--249.

\bibitem[\protect\astroncite{Macchi}{1975}]{macchi}
Macchi, O. (1975).
\newblock The coincidence approach to stochastic point processes.
\newblock {\em Advances in Appl. Probability}, 7:83--122.

\bibitem[\protect\astroncite{Mahler}{2003}]{mahler2003}
Mahler, R. (2003).
\newblock Multitarget bayes filtering via first-order multitarget moments.
\newblock {\em IEEE Transactions on Aerospace and Electronic Systems},
  39(4):1152--1178.

\bibitem[\protect\astroncite{Mahler}{2007}]{mahler2007}
Mahler, R. (2007).
\newblock {PHD} filters of higher order in target number.
\newblock {\em IEEE Transactions on Aerospace and Electronic Systems},
  43(4):1523--1543.

\bibitem[\protect\astroncite{Mahler}{2015}]{mahlerpermanental}
Mahler, R. (2015).
\newblock Tracking ``bunching'' multitarget correlations.
\newblock In {\em IEEE International Conference on Multisensor Fusion and
  lntegration for Intelligent Systems (MFI)}, pages 102--109.

\bibitem[\protect\astroncite{Mori}{1997}]{mori}
Mori, S. (1997).
\newblock Random sets in data fusion. {M}ulti-object state-estimation as a
  foundation of data fusion theory.
\newblock In {\em Random sets ({M}inneapolis, {MN}, 1996)}, volume~97 of {\em
  IMA Vol. Math. Appl.}, pages 185--207. Springer, New York.

\bibitem[\protect\astroncite{Moyal}{1964}]{moyal2}
Moyal, J. (1964).
\newblock Multiplicative population processes.
\newblock {\em J. Appl. Probability}, 1:267--283.

\bibitem[\protect\astroncite{Moyal}{1962}]{moyal}
Moyal, J.~E. (1962).
\newblock The general theory of stochastic population processes.
\newblock {\em Acta Math.}, 108:1--31.

\bibitem[\protect\astroncite{Osada}{2013}]{osada2}
Osada, H. (2013).
\newblock Interacting {B}rownian motions in infinite dimensions with
  logarithmic interaction potentials.
\newblock {\em Ann. Probab.}, 41(1):1--49.

\bibitem[\protect\astroncite{Portenko et~al.}{1997}]{portenko}
Portenko, N., Salehi, H., and Skorokhod, A. (1997).
\newblock On optimal filtering of multitarget tracking systems based on point
  processes observations.
\newblock {\em Random Oper. Stoch. Equ.}, 5(1):1--34.

\bibitem[\protect\astroncite{Schlangen et~al.}{2018}]{sdhc}
Schlangen, I., Delande, E., Houssineau, J., and Clark, D. (2018).
\newblock A second-order {PHD} filter with mean and variance in target number.
\newblock {\em IEEE Trans. Signal Processing}, 66(1):48--63.

\bibitem[\protect\astroncite{Schumacher et~al.}{2008}]{schumacher}
Schumacher, D., Vo, B.-T., and Vo, B.-N. (2008).
\newblock A consistent metric for performance evaluation of multi-object
  filters.
\newblock {\em IEEE Trans. Signal Processing}, 56(8):3447--3457.

\bibitem[\protect\astroncite{Shirai and Takahashi}{2003}]{shirai}
Shirai, T. and Takahashi, Y. (2003).
\newblock Random point fields associated with certain {F}redholm determinants.
  {I}. {F}ermion, {P}oisson and boson point processes.
\newblock {\em J. Funct. Anal.}, 205(2):414--463.

\bibitem[\protect\astroncite{Singh et~al.}{2009}]{singh}
Singh, S., Vo, B.-N., Baddeley, A., and Zuyez, S. (2009).
\newblock Filters for spatial point processes.
\newblock {\em SIAM J. Control Optim.}, 48(4):2275--2295.

\bibitem[\protect\astroncite{Soshnikov}{2000}]{soshnikov}
Soshnikov, A. (2000).
\newblock Determinantal random point fields.
\newblock {\em Uspekhi Mat. Nauk}, 55(5(335)):107--160.

\bibitem[\protect\astroncite{van Lieshout}{1995}]{vanlieshout2}
van Lieshout, M. N.~M. (1995).
\newblock {\em Stochastic geometry models in image analysis and spatial
  statistics}, volume 108 of {\em CWI Tract}.
\newblock Stichting Mathematisch Centrum, Centrum voor Wiskunde en Informatica,
  Amsterdam.

\bibitem[\protect\astroncite{Vo and Ma}{2006}]{vo-ma}
Vo, B.-N. and Ma, W.-K. (2006).
\newblock The {G}aussian mixture probability hypothesis density filter.
\newblock {\em IEEE Transactions on Aerospace and Electronic Systems},
  54(11):4091--4104.

\bibitem[\protect\astroncite{Vo et~al.}{2005}]{vo-singh}
Vo, B.-N., Singh, S.~S., and Doucet, A. (2005).
\newblock Sequential {M}onte {C}arlo methods for multitarget filtering with
  random finite sets.
\newblock {\em IEEE Transactions on Aerospace and Electronic Systems},
  41(4):1224--1245.

\bibitem[\protect\astroncite{Vo et~al.}{2007}]{batuongvo}
Vo, B.-T., Vo, B.-N., and Cantoni, A. (2007).
\newblock Analytic implementations of the cardinalized probability hypothesis
  density filter.
\newblock {\em IEEE Trans. Signal Processing}, 55(7):3553--3567.

\bibitem[\protect\astroncite{Vo et~al.}{2009}]{cantoni}
Vo, B.-T., Vo, B.-N., and Cantoni, A. (2009).
\newblock The cardinality balanced multi-target multi-{B}ernoulli filter and
  its implementations.
\newblock {\em IEEE Trans. Signal Processing}, 57(2):409--423.

\end{thebibliography}
\end{document}